\documentclass[10pt,a4paper]{amsart}

\usepackage[T1]{fontenc}
\usepackage[utf8]{inputenc}
\usepackage[english]{babel}
\usepackage{amsmath,amsthm,amssymb,geometry,fancyhdr,enumerate,bbm,mathrsfs}
\usepackage[all]{xy}
\usepackage[active]{srcltx}
\usepackage[hypertexnames=false]{hyperref}

\date{\today}

\theoremstyle{plain}
\newtheorem{THM}{Theorem}[section]

\newtheorem{COR}[THM]{Corollary}
\newtheorem{thm}{Theorem}[subsection]
\newtheorem{cor}[thm]{Corollary}
\newtheorem{lem}[thm]{Lemma}
\newtheorem{prop}[thm]{Proposition}
\theoremstyle{remark}
\newtheorem{ex}[thm]{Example}
\newtheorem{rem}[thm]{Remark}
\theoremstyle{definition}
\newtheorem{definition}[thm]{Definition}
\newtheorem{setting}[THM]{Setting}

\numberwithin{equation}{subsection}

\begin{document}

\title{Crossed products of Calabi-Yau algebras by finite groups}

\author{Patrick Le Meur}

\thanks{This is the accepted manuscript of the article published in
  the \emph{Journal of Pure and Applied Algebra} available at
  \texttt{https://doi.org/10.1016/j.jpaa.2020.106394}. It is made
  available under the CC-BY-NC-ND license}

\address{CMLA, ENS Cachan, CNRS, UniverSud, 61 Avenue du President
  Wilson, F-94230 Cachan}

\curraddr{Universit\'e de Paris, Sorbonne Universit\'e, CNRS, Institut
  de Math\'ematiques de Jussieu-Paris Rive Gauche, F-75013 Paris,
  France}

\email{patrick.le-meur@imj-prg.fr}

\subjclass[2010]{18G20; 16W22; 18E30; 16G20; 16G70}

\keywords{Calabi-Yau algebras; skew group algebras; quivers with
  potentials; Ginzburg algebras; cluster categories; higher
  Auslander-Reiten theory}

\begin{abstract}
  Let a finite group $G$ act on a differential graded algebra
  $A$. This article presents necessary conditions and sufficient
  conditions for the skew group algebra $A*G$ to be Calabi-Yau. In
  particular, when $A$ is the Ginzburg dg algebra of a quiver with an
  invariant potential, then $A*G$ is Calabi-Yau and Morita equivalent to a
  Ginzburg dg algebra. Some applications of these results are derived
  to compare the generalised cluster categories of $A$ and $A*G$ when
  they are defined and to compare the higher Auslander-Reiten
  theories of $A$ and $A*G$ when $A$ is a finite dimensional algebra.
\end{abstract}

\maketitle

\setcounter{tocdepth}{1}
\tableofcontents

\section*{Introduction}

Calabi-Yau algebras occur in many parts of mathematics and play a
special role in several areas related to representation theory. In the
non commutative algebraic geometry initiated by Artin and Schelter
(\cite{MR917738}), they appear naturally among the regular algebras
used as non commutative analogues of polynomial algebras. They also
serve to define generalised cluster categories (\cite{A08}), which are
used to categorify the cluster algebras of Fomin and Zelevinsky. When
a cluster algebra is given by a marked surface, quivers with
potentials are associated to the triangulations of the surface
(\cite{MR2592019,MR2500873}); recall that quivers with potentials are
the base ingredients in the definition of Ginzburg differential graded
(dg) algebras, which are Calabi-Yau. Finally, in the higher
Auslander-Reiten theory developed by Iyama et al., the higher
representation infinite algebras (\cite{MR3144232}) are defined using
higher derived preprojective algebras, which are Calabi-Yau
completions of finite dimensional algebras of finite global dimension
(\cite{MR2795754}).

Recall that a dg algebra $A$ over a field ${\mathbbm{k}}$ is called $d$-Calabi-Yau, where $d$ is
an integer, when
\begin{itemize}
\item $A$ is homologically smooth, that is $A\in \mathrm{per}(A^e)$, and
\item $\Sigma^d\mathrm{RHom}_{A^e}(A,A^e)\simeq A$ in the derived category
  $\mathcal{D}(A^e)$ of $A^e$.
\end{itemize}
Here, $\Sigma$ denotes the suspension for all triangulated categories,
$A^e$ denotes the enveloping dg algebra $A\otimes_{\mathbbm{k}} A^{\mathrm{op}}$ and
$\mathrm{per}(A^e)$ denotes the perfect derived category of $A^e$.

This article investigates the Calabi-Yau duality of skew group algebras
$A*G$ of dg algebras $A$ under the action of a finite group $G$ with a
particular interest on the dg algebras which are involved in the above
mentioned theories. This is motivated by several developments some of
which being recent and others being older.

On one hand, some of the early developments of non commutative
algebraic geometry focused to the regularity properties of
invariant subrings under actions of finite groups. And this has led to
the study of skew group algebras from the viewpoint of
regularity. Recently, many articles have studied when the skew group
algebra of an ordinary (not dg) algebra is Calabi-Yau and fairly complete
answers are known when the original algebra is Koszul
(\cite{MR2813562}) connected graded (\cite{MR3250287}), or is
the universal enveloping algebra of a finite dimensional Lie algebra
(\cite{MR2678828}).

On the other hand, the recent works \cite{AP} and \cite{MR3899030} use skew
group algebras of Ginzburg dg algebras and the associated generalised
cluster categories for marked surfaces. While the former uses actions
of $\mathbb Z/2\mathbb Z$ to relate punctured surfaces to unpunctured
surfaces, the latter uses free actions of finite groups for
orbifolds. In both cases, the considered skew group algebra happens to
be a Ginzburg dg algebra up to a Morita equivalence. Recall that, in
the context of a generalised cluster category $\mathcal C$, important
information is provided by the cluster tilting objects and the cluster
tilting subcategories. A functorially finite subcategory
$\mathcal T\subseteq \mathcal C$ is called \emph{cluster tilting} if
\[
\begin{array}{rcl}
  \mathcal T & = & \{X\in \mathcal C \ |\ (\forall M\in \mathcal T)\
                   \mathcal C(X,\Sigma M) = 0 \}\,,
\end{array}
\]
and an object $T\in \mathcal C$ is called cluster tilting if its
closure $\mathrm{add}(T)$ under direct sums and direct summands is cluster
tilting.

The present article therefore aims at giving general answers to the
following questions for a given dg algebra $A$ acted upon
by a finite group $G$ by dg automorphisms.\footnote{In the original
  version of this article, \cite{lemeur2}, only the first two of these
  questions are treated.}
\begin{itemize}
\item Under which conditions (necessary and/or sufficient) is $A*G$
  Calabi-Yau?
\item To what extent are constructions of Calabi-Yau algebras
  such as Ginzburg dg algebras or deformed Calabi-Yau completions
  compatible with taking skew group algebras?
\item Assuming that $A$ and $A*G$ are Calabi-Yau and such that their
  generalised cluster categories are defined, what kind of relation is
  there between these categories, and what kind of relation is there
  between their respective cluster tilting objects?
\item When $A$ is a finite dimensional algebra, is there a connection
  between $A$ being $d$-representation infinite and $A*G$ being
  $n$-representation infinite?
\end{itemize}
Note that the dg algebras considered in this article do not satisfy
particular conditions such as being Koszul, connected graded or
Artin-Schelter.

\section{Main results and structure of the article}
\label{sec:main-results}

The main results assume the following setting
\begin{setting}
  \label{sec:main-results-1}
  Let $d$ be an integer. Let ${\mathbbm{k}}$ be a field. Let $A$ be a dg algebra
  over ${\mathbbm{k}}$. Let $G$ be a finite group acting on $A$ by dg
  automorphisms $(g,a)\mapsto \,^ga$. Denote by $\Lambda$ the skew
  group dg algebra $A*G$.
\end{setting}

Recall that $\Lambda$ is defined as follows. As a complex it is equal
to $A\otimes {\mathbbm{k}} G$, where ${\mathbbm{k}} G$ is a complex concentrated in degree
$0$, and the product is such that $(a*g)\cdot (b*h) = a\,^gb*gh$ for
all $a,b\in A$ and $g,h\in G$, where $a*g$ stands for $a\otimes g$.

The investigation of the Calabi-Yau duality for $\Lambda$ makes use of an
auxiliary dg subalgebra of $\Lambda^e$ denoted by $\Delta$,
\setcounter{subsection}{1}
\begin{equation}
  \label{eq:17}
  \Delta = \bigoplus\limits_{g\in G}(A*g)\otimes (A*g^{-1})\,.
\end{equation}
It features the following properties.
\begin{itemize}
\item $A$ is a dg $\Delta$-module and $A\otimes^{\mathrm L}_\Delta \Lambda^e \simeq
  \Lambda$ in $\mathcal{D}(\Lambda^e)$.
\item The functor ${\mathrm{Hom}}_{A^e}(-,A^e)\colon \mathrm{Mod}(A^e) \to \mathrm{Mod}(A^e)$
  induces a functor $\mathrm{Mod}(\Delta)\to \mathrm{Mod}(\Delta)$ also denoted by
  $\mathrm{Hom}_{A^e}(-,A^e)$ and whose total derived functor
  $\mathrm{RHom}_{A^e}(-,A^e)$ is such that
  $\mathrm{RHom}_{A^e}(-,A^e)\otimes^{\mathrm L}_\Delta \Lambda^e \simeq
  \mathrm{RHom}_{\Lambda^e}(-\otimes^{\mathrm L}_\Delta\Lambda^e,\Lambda^e)$.
\end{itemize}
These are used to prove the first main result of this article.
\begin{THM}[Theorem~\ref{sec:main-theorem-1}]
\label{sec:main-results-3}
Assume setting~\ref{sec:main-results-1}. Assume that
$\mathrm{char}({\mathbbm{k}})\nmid\mathrm{Card}(G)$, that $A$ is homologically smooth and that
$\Sigma^d\mathrm{RHom}_{A^e}(A,A^e)\simeq A$ in $\mathcal{D}(\Delta)$, then $A$ and
$\Lambda$ are $d$-Calabi-Yau.
\end{THM}

See Proposition~\ref{sec:main-theorem} for a necessary condition for
$\Lambda$ to be $d$-Calabi-Yau expressed in terms of $\mathrm{RHom}_{A^e}(A,A^e)$
when the Hochschild cohomology group $\mathrm{HH}^0(A)$ is a local algebra.

Several constructions of Calabi-Yau dg algebras behave nicely with respect
to taking skew group algebras. This article concentrates on two
constructions, namely, the Calabi-Yau completions (and their deformations)
of \cite{MR2795754} and the Ginzburg dg algebras of
\cite{G06,MR2795754}. When $A$ is homologically smooth, its $d$-Calabi-Yau
completion $\mathbf{\Pi}_d(A)$ is defined as $T_A(\Sigma^{d-1}\Theta_A)$ where
$\Theta_A$ is any cofibrant replacement of $\mathrm{RHom}_{A^e}(A,A^e)$ also
called an \emph{inverse dualising complex} of $A$. It is
$d$-Calabi-Yau. Any Hochschild homology class $c\in \mathrm{HH}_{d-2}(A)$ yields a
deformation of the differential of $T_A(\Sigma^{d-1}\Theta_A)$, thus
giving rise to the deformed Calabi-Yau completion $\mathbf{\Pi}_d(A,c)$. If $c$ lies
in the image of the natural mapping
$\mathrm{HN}_{d-2}(A) \to \mathrm{HH}_{n-2}(A)$, where $\mathrm{HN}$
is the negative cyclic homology, then $\mathbf{\Pi}_d(A,c)$ is $d$-Calabi-Yau
(\cite{yeung}). Both constructions depend on the choice of $\Theta_A$
and on the choice of a representative of $c$, but they are independent
of these choices up to a quasi-isomorphism of dg algebras. Now, here
is how these constructions are compatible with taking skew group
algebras.
\begin{THM}[Theorem~\ref{sec:skew-group-algebras-2}]
\label{sec:main-results-2}
Assume setting~\ref{sec:main-results-1}. Assume that $A$ is
homologically smooth and that $\mathrm{char}({\mathbbm{k}})\nmid\mathrm{Card}(G)$.
  \begin{enumerate}
  \item $A*G$ is homologically smooth and, up to appropriate choices
    of inverse dualising complexes of $A$ and $A*G$, the group $G$
    acts on $\mathbf{\Pi}_d(A)$ by dg automorphisms and
    $\mathbf{\Pi}_d(A)*G\simeq \mathbf{\Pi}_d(A*G)$ as dg algebras. In particular,
    $\mathbf{\Pi}_d(A)*G$ is $d$-Calabi-Yau.
  \item Let $c\in \mathrm{HH}_{d-2}(A)^G$. Let $\overline c$ be the image of
    $c$ under the mapping $\mathrm{HH}_{d-2}(A)\to \mathrm{HH}_{d-2}(A*G)$ induced by
    $A \to A*G,\,a \mapsto a*e$. Up to appropriate choices as in (1)
    and up to appropriate choices of representatives of $c$ and
    $\overline c$, the group $G$ acts on $\mathbf{\Pi}_d(A,c)$ by dg
    automorphisms and $\mathbf{\Pi}_d(A,c)*G \simeq \mathbf{\Pi}_d(A*G,\overline c)$ as dg
    algebras. In particular, if $c$ lifts to $\mathrm{HN}_{d-2}(A)$
    (so that $\mathbf{\Pi}_d(A,c)$ is $d$-Calabi-Yau), then $\mathbf{\Pi}_d(A,c)*G$ is $d$-Calabi-Yau.
  \end{enumerate}  
\end{THM}

Note that, if $A$ is a finite dimensional algebra of global
dimension $d$, then $H^0(\mathbf{\Pi}_{d+1}(A))$ is isomorphic to the
$(d+1)$-preprojective algebra of $A$ in the sense of \cite[Definition
2.11]{MR3077865}. Hence, part (1) of Theorem~\ref{sec:main-results-2}
entails that, in this case, $G$ acts on the $(d+1)$-preprojective
algebra of $A$ and that the resulting skew group algebra is isomorphic
to the $(d+1)$-preprojective algebra of $A*G$ (see
Proposition~\ref{sec:texorpdfstr-repr-inf-1}). In the more particular
case where $d=1$, these are classical preprojective algebras of
hereditary algebras and this result is proved in \cite{MR2578593}.

The Ginzburg dg algebras are particular instances of deformed Calabi-Yau
completions. They are Calabi-Yau. Here, they are considered using the
following setting.
\begin{setting}
  \label{sec:appl-ginzb-dg}
  Let $Q,G,W$ be as follows.
  \begin{itemize}
  \item $Q$ is a finite graded ${\mathbbm{k}}$-quiver.
\item $G$ is a finite group with order not divisible by $\mathrm{char}({\mathbbm{k}})$ and
  acting on the path algebra ${\mathbbm{k}} Q$ by degree preserving automorphisms
  in such a way that both the set of (idempotents associated with) vertices
  and the vector space generated by the arrows of $Q$ are stabilised
  by the action.
\item $W$ is a homogeneous potential of degree $d-3$ on $Q$.
  \end{itemize}
\end{setting}
The Ginzburg dg algebra $\mathcal A(Q,W)$ is quasi-isomorphic to the
deformed Calabi-Yau completion $\mathbf{\Pi}_d({\mathbbm{k}} Q,c)$, where the graded path
algebra ${\mathbbm{k}} Q$ has zero differential and $c$ is the image of $W$ under
Connes' boundary $\mathrm{HH}_{d-3}({\mathbbm{k}} Q) \to \mathrm{HH}_{d-2}({\mathbbm{k}}
Q)$. Theorem~\ref{sec:main-results-2} applies as follows.
\begin{COR}[Corollary~\ref{sec:skew-group-algebras-3}]
  \label{sec:main-results-4}
  Assume setting~\ref{sec:appl-ginzb-dg}. Assume that $W$ is
  $G$-invariant up to cyclic permutation.
  \begin{enumerate}
  \item The action of $G$ on ${\mathbbm{k}} Q$ extends to an action of $G$ on
    $\mathcal A(Q,W)$ by dg automorphisms and $\mathcal A(Q,W)*G$ is
    $d$-Calabi-Yau.
  \item For all graded ${\mathbbm{k}}$-quivers $Q'$ and for all (non unital)
    graded algebra homomorphisms ${\mathbbm{k}} Q'\to {\mathbbm{k}} Q*G$ whose
    restriction-of-scalars functor is an equivalence from $\mathrm{Mod}({\mathbbm{k}}
      Q*G)$ to $\mathrm{Mod}({\mathbbm{k}} Q')$ (see \cite{MR2578593}), there exists a
    homogeneous of degree $d-3$ potential $W'$ on $Q'$ such that ${\mathbbm{k}}
    Q'\to {\mathbbm{k}} Q*G$ extends to a (non unital) dg algebra homomorphism
    \[
    \mathcal A(Q',W')\to \mathcal A(Q,W)*G
    \]
    whose restriction-of-scalars functor induces an equivalence
    \[
    \mathcal{D}(\mathcal A(Q,W)*G) \xrightarrow{\sim} \mathcal{D}(\mathcal A(Q',W'))\,.
    \]
  \end{enumerate}
\end{COR}
More precisely, ${\mathbbm{k}} Q'\to {\mathbbm{k}} Q*G$ induces an isomorphism
$\mathrm{HH}_{d-3}({\mathbbm{k}} Q')\to \mathrm{HH}_{d-3}({\mathbbm{k}} Q*G)$. Using it as an
identification, $W'$ may be taken equal to the image of $W$ under the
mapping $\mathrm{HH}_{d-3}({\mathbbm{k}} Q) \to \mathrm{HH}_{d-3}({\mathbbm{k}} Q*G)$ induced by the natural
embedding ${\mathbbm{k}} Q\to {\mathbbm{k}} Q*G,\,a\mapsto a*e$. When $d=3$,
Corollary~\ref{sec:main-results-4} is proved in \cite[Theorem 2.6,
Corollary 2.7]{AP} assuming that $G=\mathbb Z/2\mathbb Z$ and the
action of $G$ on ${\mathbbm{k}} Q$ is induced by an action on $Q$ by quiver
automorphisms.

As mentioned in the introduction, $3$-Calabi-Yau algebras are base
ingredients in the construction of generalised cluster
categories. Theorem~\ref{sec:main-results-3} and
Corollary~\ref{sec:main-results-4} raise the question of the
comparison of generalised cluster categories of $A$ and $A*G$ when
they are defined. The setting in which this article addresses this
question is the following.
\begin{setting}
  \label{sec:main-results-5}
  Let ${\mathbbm{k}}$ be a field. Let $G$ be a finite group. Let $A$ be dg
  ${\mathbbm{k}}$-algebra acted upon by $G$ by dg automorphisms. Assume that
\begin{itemize}
\item $A$ is concentrated in non positive degrees and $H^0(A)$ is
  finite dimensional,
\item $\mathrm{char}({\mathbbm{k}})$ does not divide $\mathrm{Card}(G)$,
\item and both $A$ and $A*G$ are $3$-Calabi-Yau.
\end{itemize}
\end{setting}
This ensures the existence of the generalised cluster categories ${\mathcal{C}_A}$
and ${\mathcal{C}_{A*G}}$ of $A$ and $A*G$, respectively. In particular, they have 
distinguished cluster tilting objects $\overline A$ and
$\overline{A*G}$, respectively. In this setting, the group $G$ acts
(strictly) on ${\mathcal{C}_A}$ by (strict) automorphisms of triangulated
categories. Also, the extension-of-scalars functor $\mathcal{D}(A)\to \mathcal{D}(A*G)$
and the restriction-of-scalars functor $\mathcal{D}(A*G)\to \mathcal{D}(A)$ induce
triangle functors
\[
\xymatrix{
  F_\lambda \colon {\mathcal{C}_A}
  \ar@<+0.5ex>[rr] &&
  {\mathcal{C}_{A*G}} \colon F_\cdot
  \ar@<+0.5ex>[ll]
}
\]
which form a biadjoint pair with split units (see
Proposition~\ref{sec:adjunct-betw-ca-2}). These functors are the
basement of the comparison between the cluster tilting objects
of ${\mathcal{C}_A}$ and those of ${\mathcal{C}_{A*G}}$.
\begin{THM}[Corollary~\ref{sec:clust-tilt-objects-3} and
  Proposition~\ref{sec:clust-tilt-objects-4}]
  \label{sec:main-results-6}
  Assume setting~\ref{sec:main-results-5}.
  \begin{enumerate}
  \item The assignments $\mathcal T \mapsto \mathrm{add}(F_\lambda \mathcal
      T)$ and $\mathcal T'\mapsto \mathrm{add}(F_\cdot \mathcal T')$ induce
    mutually inverse bijections between the set of $G$-stable cluster
    tilting subcategories of ${\mathcal{C}_A}$ and the set of cluster tilting
    subcategories of ${\mathcal{C}_{A*G}}$ which are stable under $F_\lambda
    F_\cdot$.
  \item For all $G$-stable cluster tilting subcategories $\mathcal T$
    of ${\mathcal{C}_A}$, there exists a cluster tilting object
    $T\in {\mathcal{C}_A}$ such that $\mathcal T=\mathrm{add}(T)$ and
    $T^g=T$ for all $g\in G$. Then,
    ($\mathrm{add}(F_\lambda \mathcal T) = \mathrm{add}(F_\lambda
    T)$), the group $G$ acts on $\mathrm{End}(T)$ by algebra
    automorphisms and $\mathrm{End}(T)*G$ is Morita equivalent to
    $\mathrm{End}(F_\lambda T)$.
  \end{enumerate}  
\end{THM}
Part (1) of Theorem~\ref{sec:main-results-6} is similar to some of the
results of \cite[Section 3]{MR2844757} where a correspondence is
established between the $G$-stable cluster tilting objects of a given
$\mathrm{Hom}$-finite Krull-Schmidt and stably $2$-Calabi-Yau Frobenius category
endowed with a suitable action of $G$ and the ${\mathbbm{k}}[G]$-stable cluster
tilting objects of the associated $G$-equivariant category.

With additional assumptions, part (1) of
Theorem~\ref{sec:main-results-6} has a more precise
formulation. Indeed, assume that the action of $G$ on the isomorphism
classes of indecomposable projective $H^0(A)$-modules is free. Then,
Proposition~\ref{sec:clust-tilt-objects-6} asserts that, for all
cluster tilting objects $T\in {\mathcal{C}_{A*G}}$ in the mutation class of
$\overline{A*G}$, every indecomposable direct summand of $T$ lies in
the essential image of $F_\lambda$; moreover, $T$ is isomorphic to
$F_\lambda \widetilde T$ for some cluster tilting object
$\widetilde T\in{\mathcal{C}_A}$ such that $\widetilde T^g=T$ for all $g\in G$.
When $A$ is given by a quiver with potential like in
Corollary~\ref{sec:main-results-4}, the above condition relative to
$H^0(A)$ is equivalent to the action of $G$ on the set of vertices of
$Q$ being free.

These results have particular implications in the acyclic case. A
$\mathrm{Hom}$-finite, Krull-Schmidt and $2$-Calabi-Yau category is called
\emph{acyclic} when it is equivalent as a triangulated category to the
cluster category, in the sense of \cite{BMRRT06}, of a hereditary
algebra. Following \cite{KR08}, this is the case if and only if the
category is algebraic and contains a cluster tilting object whose
endomorphism algebra is hereditary.
\begin{COR}[Corollary~\ref{sec:acyclic-case}]
  \label{sec:main-results-7}
  Assume setting~\ref{sec:main-results-5}.
  \begin{enumerate}
  \item The following assertions are equivalent.
    \begin{enumerate}[(i)]
    \item ${\mathcal{C}_A}$ is equivalent to the cluster category of a Dynkin
      quiver.
    \item ${\mathcal{C}_{A*G}}$ is equivalent to the cluster category of a Dynkin quiver.
    \end{enumerate}
  \item If ${\mathcal{C}_A}$ is acyclic and has infinitely many isomorphism
    classes of indecomposable objects, then
    \begin{enumerate}[(a)]
    \item there exists a cluster tilting object $T\in {\mathcal{C}_A}$ such that
      $\mathrm{End}(T)$ is hereditary and $T^g = T$ for all $g\in G$ and
    \item for any such $T$, there exists a cluster tilting object
      $T'\in {\mathcal{C}_{A*G}}$ such that the algebras $\mathrm{End}(T)*G$ and $\mathrm{End}(T')$
      are Morita equivalent.
  \end{enumerate}
  Consequently, ${\mathcal{C}_{A*G}}$ is acyclic.
  \item Assume that $G$ acts freely on the set of isomorphism classes
    of indecomposable projective $H^0(A)$-modules. If ${\mathcal{C}_{A*G}}$ is
    acyclic, then so is ${\mathcal{C}_A}$.
  \end{enumerate}
\end{COR}
Note that part (1) is proved in \cite[Proposition 7.12]{MR3899030} when $A$
is given by a quiver with potential $(Q,W)$ and the action of $G$ on
$A$ is induced by a free action on $Q$ by quiver automorphisms in such
a way that $W$ is $G$-invariant up to cyclic permutation.

The last main result of this article is an application of
Theorem~\ref{sec:main-results-2} to higher Auslander-Reiten
theory. Indeed, a finite dimensional algebra of global dimension $d$
is $d$-representation infinite, in the sense of \cite{MR3144232}, if
and only if the cohomology of its derived $(d+1)$-preprojective
algebra is concentrated in degree $0$. Since the derived
$(d+1)$-preprojective algebra is equal to the $(d+1)$-Calabi-Yau
completion, Theorem~\ref{sec:main-results-2} applies to prove the
following.
\begin{COR}[Propositions~\ref{sec:texorpdfstr-repr-fin} and
  \ref{sec:texorpdfstr-repr-inf-1}]
  \label{sec:main-results-8}
  Let $A$ be a finite dimensional algebra of global dimension $d$. Let
  $G$ be a finite group acting on $A$ by algebra automorphisms. Assume
  that $\mathrm{char}({\mathbbm{k}})$ does not divide $\mathrm{Card}(G)$. 
  \begin{enumerate}
  \item $A$ is $d$-representation finite if and only if so is $A*G$.
  \item $A$ is $d$-representation infinite if and only if so is $A*G$.
  \end{enumerate}
\end{COR}

The article is organised as follows. Basic material is recalled in
section~\ref{sec:not}. Section~\ref{sec:s2} develops the needed
properties of the dg algebra $\Delta$ and proves
Theorems~\ref{sec:main-results-3} and
\ref{sec:main-results-2}. Section~\ref{sec:smash-prod-ginzb} is
devoted to the application to Ginzburg dg algebras, it proves
Corollary~\ref{sec:main-results-4}. The comparison between the
generalised cluster categories of $A$ and $A*G$, when these are
defined, is made in section~\ref{sec:smash-prod-ginzb}; this includes
proofs of Theorem~\ref{sec:main-results-6} and
Corollary~\ref{sec:main-results-7}. Finally,
Corollary~\ref{sec:main-results-8} on higher representation (in)finite
algebras is proved in section \ref{sec:appl-high-ausl}.

\section{Definitions and notation}
\label{sec:not}
Throughout the text, ${\mathbbm{k}}$ denotes a field. The tensor product over
${\mathbbm{k}}$ is denoted by $\otimes$. Let $A$ be a dg algebra over ${\mathbbm{k}}$.
All dg algebras are over ${\mathbbm{k}}$ and with differential of degree
$1$. Denote by $A^{\mathrm{op}}$ the \emph{opposite dg algebra} of $A$. Denote
by $A^e$ the \emph{enveloping dg algebra} $A\otimes A^{\mathrm{op}}$ of $A$,
with product given by
$(a\otimes b) \cdot (c\otimes d) = (-1)^{\mathrm{deg}(b)\cdot (\mathrm{deg}(c)+\mathrm{deg}(d))}
ac \otimes db$ for all (homogeneous) $a,b,c,d\in A$.

The dg category of dg (right) $A$-modules is denoted by $\mathrm{Mod}(A)$.  See
\cite{MR2275593} for a background. In all places, the suspension
functor is denoted by $\Sigma$. Given dg $A$-modules $M,N$ the complex
of morphisms $M\to N$ in $\mathrm{Mod}(A)$ is denoted by $\mathrm{Hom}_A(M,N)$; in
particular $Z^0\mathrm{Hom}_A(M,N)$ is the vector space of morphisms of dg
$A$-modules, that is, with degree $0$ and commuting with the
differentials, and $H^0\mathrm{Hom}_A(M,N)$ is the vector space of morphisms
of dg $A$-modules up to homotopy.  The \emph{derived category} of $A$
is denoted by $\mathcal{D}(A)$, this is the localisation of the graded category
$H^0\mathrm{Mod}(A)$ at quasi-isomorphisms. It contains two full subcategories
used here,
\begin{itemize}
\item the \emph{perfect} derived category denoted by $\mathrm{per}(A)$, this is
  the thick triangulated subcategory of $\mathcal{D}(A)$ generated by $A$ and
\item the full subcategory $\mathcal{D}_{\mathrm{fd}}(A)$ of objects $X\in \mathcal{D}(A)$ such that
  $\mathrm{dim}(\oplus_{i\in \mathbb Z}H^i(X))<\infty$.
\end{itemize}
In particular, $A$ is a dg $A^e$-module for the following
structure map
\[
\begin{array}{rcl}
  A \otimes A^e & \to & A^e \\
  a \otimes (b \otimes c) & \mapsto & (-1)^{\mathrm{deg}(c) \cdot(\mathrm{deg}(b) +\mathrm{deg}(a))} cba
\end{array}
\]
and $A$ is called \emph{homologically smooth} if $A\in \mathrm{per}(A^e)$.

\subsection{Calabi-Yau algebras} 

For all $X\in \mathrm{Mod}(A)$, the complex $\mathrm{Hom}_{A^e}(X,A^e)$ is viewed as a
dg $A^e$-module such that
$(f\cdot (a \otimes b))(x) = (-1)^{\mathrm{deg}(f)\cdot(\mathrm{deg}(a) +\mathrm{deg}(b)) + \mathrm{deg}(
  b) \cdot \mathrm{deg}(a)} (b \otimes a) \cdot f(x)$
for all $f\in \mathrm{Hom}_{A^e}(X,A^e)$, $x\in X$ and $a,b\in A$.  When $A$
is homologically smooth, $\mathrm{RHom}_{A^e}(-,A^e)$ induces a duality
$\mathrm{RHom}_{A^e}(-,A^e)\colon \mathrm{per}(A^e)^{\mathrm{op}}\to \mathrm{per}(A^e)$.  Let $d$ be an
integer. Following \cite[3.2.3]{G06}, the dg algebra $A$ is called
$d$-Calabi-Yau if the two following assertions hold true:
\begin{enumerate}
\item $A$ is homologically smooth,
\item $\Sigma^d\mathrm{RHom}_{A^e}(A,A^e)\simeq A$ in $\mathcal{D}(A^e)$.
\end{enumerate}

\subsection{Skew group algebras}
\label{sec:smash-products}

Let $G$ be a finite group and $A$ be a dg algebra. In all places, $e$
denotes the neutral element of $G$. By an \emph{action of $G$ on $A$
  by dg automorphism} is meant a morphism of complexes
${\mathbbm{k}} G \otimes A \to A$, where ${\mathbbm{k}} G$ is a
complex in degree $0$ and with zero differential, such that, denoting
by $^ga$ the image of $g\otimes a$ for all $g\in G$ and $a\in A$, then
\begin{itemize}
\item $^g(ab) = \,^ga\,^gb$ for all $a,b\in A$ and $g\in G$ and
\item $^ea = a$ and $^{gh}a = \,^g(\,^h a)$ for all $a\in A$ and
  $g,h\in G$.
\end{itemize}
This defines the \emph{skew group} dg algebra $A*G$ with underlying
complex of vector spaces equal to $A\otimes {\mathbbm{k}} G$, where any tensor
$a\otimes g$ with $a\in A$ and $g\in G$ is denoted by $a*g$ and with
product such that $(a*g) \cdot (b *h) = a\,^gb*gh$ for all $a,b\in A$
and $g,h\in G$.

For all $g\in G$, denote by $A*g$ the dg $A^e$-module
$A*g=\{a*g\ |\ a\in A\}$. It is isomorphic to $A$ as a dg
$A$-module. Besides, $A*G = \oplus_{g\in G} A*g$ as a dg
$A^e$-module.

\subsection{Adjunctions}
\label{sec:adjunctions}

Let $G$ be a finite group acting on a dg algebra $A$ by dg automorphism.
For all $g\in G$ denote the restriction of scalars along the dg
automorphism $A\to A\,,a\mapsto \,^ga$ as follows
\begin{equation}
  \label{eq:29}
\begin{array}{rcl}
  \mathrm{Mod}(A)& \to & \mathrm{Mod}(A) \\
  M & \mapsto & M^g \\
  M\xrightarrow{f}M' & \mapsto & M^g\xrightarrow{f^g} M'^g\,.
\end{array}
\end{equation}
This defines a \emph{strict} action of $G$ on $\mathrm{Mod}(A)$
on the right by automorphisms of dg categories. It preserves
both projective and injective model structures of $\mathrm{Mod}(A)$, and hence
defines a strict action of $G$ on $\mathcal{D}(A)$ by strict automorphisms of
triangulated category.

Since $A*g \simeq A$ as dg $A$-modules for all $g\in G$, the
extension-of-scalars functor
$- \otimes_A (A*G) \colon \mathrm{Mod}(A) \to \mathrm{Mod}(A*G)$ preserves
quasi-isomorphisms and so does the restriction-of-scalars functor
$\mathrm{res}\colon \mathrm{Mod}(A*G) \to \mathrm{Mod}(A)$. Consider the adjoint pair of
triangle functors
\[
\xymatrix{
  - \otimes^{\mathrm L}_A (A*G) \colon \mathcal{D}(A) \ar@<+0.5ex>[rr] &&
  \mathcal{D}(A*G) \colon \mathrm{Rres}  \ar@<+0.5ex>[ll] }
\]

\begin{lem}
  \label{sec:adjunctions-1}
  The functors $- \otimes^{\mathrm L}_A A*G$ and $\mathrm{Rres}$ have
  the following properties.
  \begin{enumerate}
  \item The pairs $(\mathrm{Rres}\,,-\otimes^{\mathrm L}_A (A*G))$
    and $(-\otimes^{\mathrm L}_A (A*G),\mathrm{Rres}\,)$ are
    adjoint.
  \item The unit of the adjoint pair $(-\otimes^{\mathrm L}_A (A*G),\mathrm{Rres}\,)$ splits.
  \item If $\mathrm{char}({\mathbbm{k}})\nmid \mathrm{Card}(G)$, then the unit morphism
    $N \to (\mathrm{Rres}\,N) \otimes^{\mathrm L}_A (A*G)$ of the adjoint pair
    $(\mathrm{Rres}\,,-\otimes^{\mathrm L}_A (A*G))$ splits functorially for all
    $N\in \mathcal{D}(A*G)$.
  \item For all $M\in \mathcal{D}(A)$ there is a functorial isomorphism in $\mathcal{D}(A)$
    \begin{equation}
      \label{eq:30}
      \mathrm{Rres}\,( M\otimes^{\mathrm L}_A (A*G)) \simeq \bigoplus_{g\in G}M^g\,.
    \end{equation}
  \item For all $M,M'\in \mathcal{D}(A)$ there is a bifunctorial isomorphism
    \[
    \bigoplus\limits_{g\in G} \mathcal{D}(A)(M,M'^g) \xrightarrow{\sim}
    \mathcal{D}(A*G)(M\otimes^{\mathrm L}_A(A*G),M' \otimes^{\mathrm L}_A(A*G))\,.
    \]
  \item There exists an isomorphism of functors
    $\lambda_g \colon -^g\otimes^{\mathrm L}_A(A*G) \xrightarrow{\sim}
    -\otimes^{\mathrm L}_A (A*G)$, for all $g\in G$, such that the following
    equalities hold true, for all $M\in \mathcal{D}(A)$ and
    $g,h\in G$,
    \[
    \left\{
      \begin{array}{rcl}
        (\lambda_e)_M & = & \mathbbm 1_M \\
        (\lambda_g)_M \circ (\lambda_h)_{M^g} & = & (\lambda_{gh})_M\,.
      \end{array}\right.
    \]
  \item $\mathrm{Rres} \simeq (\mathrm{Rres})^g$ for all $g\in G$.
  \end{enumerate}
\end{lem}
\begin{proof} This is known for finitely generated modules over Artin
  algebras. Here are some details.

  (1) See \cite[Theorem 1.4]{RR85}. Given $X\in \mathrm{Mod}(A*G)$ and
  $Y\in \mathrm{Mod}(A)$, the following is an isomorphism of complexes
  \[
  \Phi \colon \mathrm{Hom}_A(\mathrm{res}\,X,Y) \to \mathrm{Hom}_{A*G}(X,Y \otimes_A (A*G))\,,
  \]
  where $\Phi(f)$ is defined as follows, for all $f\in\mathrm{Hom}_A(\mathrm{res}\,X,Y)$,
  \[
  \begin{array}{crcl}
    \Phi(f) \colon & X & \to & Y\otimes_A (A*G) \\
    & x & \mapsto & \sum_{g\in G} f(x \cdot (1*g^{-1})) \otimes (1*g)\,.
  \end{array}
  \]
  The inverse of $\Phi$ maps any
  $f'\in \mathrm{Hom}_{A*G}(X,Y\otimes_A (A*G))$ to its composition with the
  mapping $Y\otimes_A A*G \to Y$ induced by the counit ${\mathbbm{k}} G\to {\mathbbm{k}}$.
  Hence, the pair $(\mathrm{res},-\otimes_A (A*G))$ of functors between
  $\mathrm{Mod}(A)$ and $\mathrm{Mod}(A*G)$ is adjoint.  This is a Quillen adjunction,
  which yields the adjunction at the level of derived categories (note
  that both $\mathrm{res}$ and $- \otimes_A (A*G)$ preserve
  quasi-isomorphisms).

  (2) See \cite[Theorem 1.4]{RR85}. The unit morphism
  $M\to \mathrm{Rres}\,(M\otimes^{\mathrm L}_A (A*G))$ is induced by
  \[
  \begin{array}{rcl}
    M & \to & M \otimes_A (A*G) \\
    m & \mapsto & m \otimes (1*e)\,.
  \end{array}
  \]
  A splitting in $\mathcal{D}(A)$ is induced by the counit ${\mathbbm{k}} G\to {\mathbbm{k}}$.

  (3) The unit morphism
  $N\to (\mathrm{Rres}\,N) \otimes^{\mathrm L}_A (A*G)$ is induced by
  \[
  \begin{array}{rcl}
    N & \to &  N\otimes_A (A*G) \\
    n & \mapsto & \sum_{g\in G} (n \cdot (1*g^{-1}))
                  \otimes (1*g)\,.
  \end{array}
  \]
  A splitting in $\mathcal{D}(A*G)$ is induced by $n \otimes (a*g)
  \mapsto \frac{1}{\mathrm{Card}(G)}\cdot n \cdot (a*g)$.

  (4) See \cite[Proposition 1.8]{RR85}. The isomorphism is induced by
  \[
  \begin{array}{rcl}
    \bigoplus_{g\in G}M^g & \to & M \otimes_A (A*G) \\
    \sum_{g\in G}m_g & \mapsto & \sum_{g\in G} m_g \otimes (1*g)\,.
  \end{array}
  \]

  (5) follows from (1) and (4).

  (6) There exists an isomorphism of dg $A*G$-modules, for all
  $M\in \mathrm{Mod}(A)$ and $g\in G$,
  \[
  \begin{array}{crcl}
    \ell_{g,M}\colon & M^g \otimes_A (A*G) & \to & M \otimes_A (A*G) \\
                     & m \otimes (a * h) & \mapsto & m \otimes(\,^ga \otimes h)\,.
  \end{array}
  \]
  And $\ell_{e,M} = \mathbbm 1_M$ and
  $\ell_{h,M^g}\circ \ell_{g,M}=\ell_{gh,M}$ for all $M\in \mathrm{Mod}(A)$ and
  $g,h\in G$. This proves (6).

  (7) $N\to N^g,\,n \mapsto n \cdot(1*g^{-1})$ is an isomorphism of dg
  $A$-modules, for all $N\in \mathrm{Mod}(A*G)$.
\end{proof}

\section{The main theorems}
\label{sec:s2}

This section assumes setting~\ref{sec:main-results-1}. This setting is
restated below for convenience.
\begin{setting}
  \label{sec:main-results-1-bis}
  Let $d$ be an integer. Let ${\mathbbm{k}}$ be a field. Let $A$ be a dg algebra
  over ${\mathbbm{k}}$. Let $G$ be a finite group acting on $A$ by dg
  automorphisms $(g,a)\mapsto \,^ga$. Denote by $\Lambda$ the dg
  algebra $A*G$.
\end{setting}
The section aims at proving Theorem~\ref{sec:main-results-3} which
gives sufficient conditions for $\Lambda$ to be $d$-Calabi-Yau and
Theorem~\ref{sec:main-results-2} which explains when $G$ acts on the
(deformed) Calabi-Yau-completions of $A$ so that the resulting skew group
algebras are Calabi-Yau.

For this purpose, an auxiliary dg subalgebra $\Delta$ of $\Lambda^e$
such that $A\otimes_\Delta \Lambda^e\simeq \Lambda$ is introduced in
subsection~\ref{sec:s2.1}. This allows to relate $A$ being
homologically smooth to $\Lambda$ being homologically smooth in
subsection~\ref{sec:homol-smoothn-a}. The proofs of these results are
based on a relationship between the inverse dualising complexes of $A$
and $\Lambda$ established in
subsection~\ref{sec:s2.3}. Theorem~\ref{sec:main-results-3} is proved
in subsection~\ref{sec:s2.4} and Theorem~\ref{sec:main-results-2} is
proved in subsection~\ref{sec:smash-prod-deform}.

\subsection{The dg algebra \texorpdfstring{$\Delta$}{Δ} and its dg modules}
\label{sec:s2.1}

Recall that $\Delta$ is defined in \eqref{eq:17}.  The group $G$ acts
on $A^e$ by dg automorphisms as follows
\[
\begin{array}{rcl}
  G \times A^e & \to & A^e \\
  (g,a\otimes b) & \mapsto &\,^g(a\otimes b) = \,^ga \otimes \,^gb\,,
\end{array}
\]
and the resulting skew group algebra is isomorphic to $\Delta$
\emph{via} the following isomorphism
\begin{equation}
  \label{eq:82}
\begin{array}{rcl}
  \Delta & \to & A^e * G \\
  (a*g) \otimes (b*g^{-1}) & \mapsto & (a\otimes\,^{g^{-1}}b) * g\,.
\end{array}
\end{equation}

To be given a dg $\Delta$-module is equivalent to be given a dg
$A^e$-module $M$ together with a \emph{compatible} action of $G$ on
$M$ by automorphisms of complexes of vector spaces
$G \times M \to M\,,(g,m) \mapsto \,^gm$. Here, by ``compatible'', is
meant that, for all $g\in G$, $a,b\in A$ and $m\in M$,
\begin{equation}
  \label{eq:36}
^g(m\cdot (a\otimes b)) = \,^gm \cdot \,^g(a\otimes b)\,,
\end{equation}
the structure of dg $\Delta$-module is then such that
$m\cdot ((a*g) \otimes (b*g^{-1})) = (\,^{g^{-1}}m) \cdot
(\,^{g^{-1}}a\otimes b)$.

To any dg $\Delta$-module $M$ is associated a dg $\Lambda^e$-module
$M*G$ with underlying complex of vector spaces $M\otimes {\mathbbm{k}} G$ and
with action of $\Lambda^e$ such that
$(m*g) \cdot (a*h\otimes b*k) = (\,^km\cdot (\,^{kg}a \otimes b))*kgh$
for all $m\in M$, $a,b\in A$, and $g,h,k\in G$.  Note that
$A\in \mathrm{Mod}(\Delta)$ by means of the following structure mapping
\[
\begin{array}{rcl}
  A \otimes \Delta & \to & A \\
  a \otimes (b*g \otimes c*g^{-1}) & \mapsto & (-1)^{\mathrm{deg}(c) \cdot(\mathrm{deg}(b) + \mathrm{deg}(a))} c\
                                               ^{g^{-1}}a\ ^{g^{-1}}b\,.
\end{array}
\]
Of course, in this case, the resulting dg $\Lambda^e$-module $A*G$ is
$\Lambda$.  The following basic properties are used throughout the
section.
\begin{lem}
  \label{sec:dg-algebra-delta}
  Assume setting~\ref{sec:main-results-1-bis}.
  \begin{enumerate}
  \item $\Delta\simeq (A^e)^{\mathrm{Card}(G)}$ as dg $A^e$-modules.
  \item $\Lambda^e\simeq \Delta^{\mathrm{Card}(G)}$ as dg $\Delta$-modules.
  \item Let $M\in\mathrm{Mod}(\Delta)$, then
    \begin{itemize}
    \item the following mapping is an isomorphism of dg $\Lambda^e$-modules,
      \begin{equation}
        \label{eq:7}
        \begin{array}{rcl}
          M \otimes_\Delta \Lambda^e & \to & M*G \\
          m \otimes (a*g \otimes b*h) & \mapsto & (m \cdot (a*h^{-1}
                                                  \otimes
                                                  b * h)) * hg\,.
        \end{array}
      \end{equation}
    \item $M$ is a direct summand of
      $M\underset{\Delta}{\otimes}\Lambda^e$ in $Z^0\mathrm{Mod}(\Delta)$,
    \item if $\mathrm{char}({\mathbbm{k}})\nmid \mathrm{Card}(G)$ then $M$ is a direct
      summand of $M\underset{A^e}{\otimes}\Delta$ in $Z^0\mathrm{Mod}(\Delta)$.
    \end{itemize}
  \item The restriction-of-scalars functor $\mathrm{Mod}(\Delta)\to \mathrm{Mod}(A^e)$
    maps cofibrant objects to cofibrant objects. Moreover, for all
    cofibrant resolutions $X\to A$ in $\mathrm{Mod}(\Delta)$, the composite
    morphism
    $X\otimes_\Delta \Lambda^e \to A \otimes_\Delta \Lambda^e
    \xrightarrow{~(\ref{eq:7})}\Lambda$
    is a cofibrant resolution in $\mathrm{Mod}(\Lambda^e)$.
  \end{enumerate}
\end{lem}
\begin{proof}
  (1) follows from the definition of $\Delta$ and the fact that
  $A*g\simeq A$ as dg $A$-modules.

  (2) follows from the direct sum decompositions
  \begin{equation}
    \label{eq:6}
    \Lambda^e = \bigoplus_{t\in G} (\bigoplus_{g\in G} (A*g) \otimes
    (A*g^{-1}t)) = \bigoplus_{t\in G} (\bigoplus_{g\in G} (A*g) \otimes
    (A*tg^{-1}))
  \end{equation}
  where each term $\oplus_{g\in G} (A*g) \otimes (A*g^{-1}t)$ (or,
  $\oplus_{g\in G} (A*g) \otimes (A*tg^{-1})$) is a dg
  $\Delta\otimes (A^e)^{\mathrm{op}}$-submodule (or, a dg
  $A^e\otimes \Delta^{\mathrm{op}}$-submodule) of $\Lambda^e$ freely generated
  by $(1*e)\otimes (1*t)$ as a dg $\Delta$-module (or, as a dg
  $\Delta^{\mathrm{op}}$-module), and hence is isomorphic to $\Delta$ as a dg
  $\Delta$-module (or, as a dg $\Delta^{\mathrm{op}}$-module, respectively).

  (3) The mapping (\ref{eq:7}) is a well-defined morphism of dg
  $\Lambda^e$-modules, its inverse is
  \[
  \begin{array}{rcl}
    M * G & \to & M \otimes_{\Delta} \Lambda^e \\
    m*g & \mapsto & m \otimes((1*g)\otimes (1*e))\,.
  \end{array}
  \]

  In order to prove that $M$ is a direct summand of
  $M\otimes_\Delta \Lambda^e$ in $Z^0\mathrm{Mod}(\Delta)$, denote by $X_t$ the dg
  $\Delta\otimes (A^e)^{\mathrm{op}}$-submodule
  $\oplus_{g\in G} (A*g) \otimes (A*g^{-1}t)$ of $\Lambda^e$, for all
  $t\in G\backslash\{e\}$. In the decomposition (\ref{eq:6}), both
  summands $X_e$ and $\oplus_{t\in G\backslash\{e\}} X_t$ are dg
  $\Delta \otimes \Delta^{\mathrm{op}}$-submodules. Accordingly, there is a
  direct sum decomposition in $Z^0\mathrm{Mod}(\Delta)$
    \[
    M = \left(M\otimes_\Delta X_e\right) \bigoplus
    \left(M\otimes_\Delta(\bigoplus_{t\neq e}X_t)\right)\,.
    \]
    Since $X_e = \Delta$, the dg module $M$ is a direct summand of
    $M\otimes_\Delta \Lambda^e$ in $Z^0\mathrm{Mod}(\Delta)$.

    Finally, note that Lemma~\ref{sec:adjunctions-1} may be applied
    here to the dg algebras $A^e$ and $\Delta\simeq A^e*G$ instead of
    to $A$ and $A*G$, respectively. Using part (3) of that lemma
    yields that $M$ is a direct summand of $M \otimes_{A^e} \Delta$ in
    $Z^0\mathrm{Mod}(\Delta)$.

    (4) Let $X\to A$ be a cofibrant resolution in $\mathrm{Mod}(\Delta)$.

    In view of (1), the functor
    $\mathrm{Hom}_{A^e}(\Delta,-) \colon \mathrm{Mod}(A^e)\to \mathrm{Mod}(\Delta)$ preserves
    quasi-isomorphisms and component-wise surjective
    morphisms. Accordingly, its left adjoint, the
    restriction-of-scalars functor $\mathrm{Mod}(\Delta) \to \mathrm{Mod}(A^e)$,
    preserves cofibrant objects. Thus, $X\to A$ is a cofibrant
    resolution in $\mathrm{Mod}(A^e)$.

    The restriction-of-scalars functor $\mathrm{Mod}(\Lambda^e)\to \mathrm{Mod}(\Delta)$
    preserves quasi-isomorphisms and component-wise surjective
    morphisms. Accordingly, its left adjoint, the
    extension-of-scalars functor
    $-\otimes_\Delta \Lambda^e \colon \mathrm{Mod}(\Delta) \to \mathrm{Mod}(\Lambda^e)$
    preserves cofibrant objects. In view of (2), it also preserves
    quasi-isomorphisms. Thus,
    $X\otimes_\Delta \Lambda^e \to A\otimes_\Delta \Lambda^e$ is a
    cofibrant resolution in $\mathrm{Mod}(\Lambda^e)$ and remains so after
    composing with the quasi-isomorphism (\ref{eq:7}).
\end{proof}

\subsection{Homological smoothness of \texorpdfstring{$A$}{A} and \texorpdfstring{$\Lambda$}{Λ}}
\label{sec:homol-smoothn-a}

\begin{prop}
  \label{sec:homol-smoothn-a-1}
  Assume setting~\ref{sec:main-results-1-bis}.
  \begin{enumerate}
  \item If $\Lambda$ is homologically smooth, then so is $A$.
  \item The converse holds true when $\mathrm{char}({\mathbbm{k}})$ does not divide $\mathrm{Card}(G)$.
  \end{enumerate}
\end{prop}
\begin{proof}
  (1) Assume that $\Lambda$ is homologically smooth. Consider the
  restriction-of-scalars functors,
  \[
  \mathcal{D}(\Lambda^e) \xrightarrow{\mathrm{res}_1} \mathcal{D}(\Delta)
  \xrightarrow{\mathrm{res}_2} \mathcal{D}(A^e)\,.
  \]
  Since $\Lambda^e\simeq \Delta^{\mathrm{Card}(G)}$ (see
  Lemma~\ref{sec:dg-algebra-delta}, part (2)), the triangle functor
  $\mathrm{res}_1$ maps $\mathrm{per}(\Lambda^e))$ into $\mathrm{per}(\Delta)$. Moreover,
  $\Lambda\in \mathrm{per}(\Lambda^e)$. Therefore,
  \begin{equation}
    \label{eq:8}
    \Lambda\in \mathrm{per}(\Delta)\,.
  \end{equation}
  Applying part (3) of Lemma~\ref{sec:dg-algebra-delta} to $M=A$ and
  using the restriction-of-scalars functor
  $\mathcal{D}(\Lambda^e)\to \mathcal{D}(\Delta)$ yields that
  $\Lambda \simeq A \otimes_\Delta \Lambda^e$ in $\mathcal{D}(\Delta)$ and that
  $A$ is a direct summand of $A\otimes_\Delta \Lambda^e$ in
  $\mathcal{D}(\Delta)$. Therefore, it follows from (\ref{eq:8}) that
  \begin{equation}
    \label{eq:9}
    A\in \mathrm{per}(\Delta)\,.
  \end{equation}
  Finally, $\Delta\in \mathrm{per}(A^e)$ because
  $\Delta \simeq (A^e)^{\mathrm{Card}(G)}$ (see
  Lemma~\ref{sec:dg-algebra-delta}, part (1)); accordingly,
  $\mathrm{res}_2(\mathrm{per}(\Delta)) \subseteq \mathrm{per}(A^e)$; therefore
  $A\in \mathrm{per}(A^e)$ (see \eqref{eq:9}), that is, $A$ is homologically
  smooth.

  (2) Assume that $\mathrm{char}({\mathbbm{k}})$ does not divide $\mathrm{Card}(G)$ and that $A$ is
  homologically smooth. Consider the extension-of-scalars functors,
  \[
  \mathcal{D}(A^e) \xrightarrow{-\otimes^{\mathrm L}_{A^e} \Delta} \mathcal{D}(\Delta)
  \xrightarrow{ - \otimes^{\mathrm L}_\Delta \Lambda^e} \mathcal{D}(\Lambda^e)\,.
  \]
  First, $A\otimes^{\mathrm L}_{A^e} \Delta \in \mathrm{per}(\Delta)$ because
  $A \in \mathrm{per}(A^e)$. Besides,
  $A \otimes^{\mathrm L}_{A^e} \Delta \simeq A\otimes_{A^e} \Delta$ in $\mathcal{D}(\Delta)$
  because $\Delta \simeq (A^e)^{\mathrm{Card}(G)}$ as dg $(A^e)^{\mathrm{op}}$-modules
  (see Lemma~\ref{sec:dg-algebra-delta}, part (1)). Now, $A$ is a
  direct summand of $A\otimes_{A^e} \Delta$ in $\mathcal{D}(\Delta)$ (see above
  or part (3) of Lemma~\ref{sec:dg-algebra-delta}). Thus,
  \begin{equation}
    \label{eq:10}
    A\in \mathrm{per}(\Delta)\,.
  \end{equation}
  Finally,
  $A\otimes_\Delta \Lambda^e \simeq A\otimes^{\mathrm L}_\Delta\Lambda^e$ because
  $\Lambda^e \simeq \Delta^{\mathrm{Card}(G)}$ as dg $\Delta^{\mathrm{op}}$-modules, and
  $\Lambda \simeq A\otimes_{\Delta}\Lambda^e$ in $\mathcal{D}(\Lambda^e)$ (see
  Lemma~\ref{sec:dg-algebra-delta}, parts (2) and (3)). Therefore
  $\Lambda \in \mathrm{per}(\Lambda^e)$ (see~\eqref{eq:10}), that is,
  $\Lambda$ is homologically smooth.
\end{proof}

\subsection{Interaction between duality and extension of scalars}
\label{sec:s2.3}

This subsection establishes a relationship between $\mathrm{RHom}_{A^e}(A,A^e)$ and
$\mathrm{RHom}_{\Lambda^e}(\Lambda,\Lambda^e)$. For this purpose, the
following lemma makes of $\mathrm{RHom}_{A^e}(-,A^e)$ a triangle functor
$\mathcal{D}(\Delta)\to \mathcal{D}(\Delta)$.
\begin{lem}
  \label{sec:inter-betw-dual}
  Assume setting~\ref{sec:main-results-1-bis}. Let $M\in\mathrm{Mod}(\Delta)$. For
  all $f\in \mathrm{Hom}_{A^e}(M,A^e)$ and $g\in G$, denote by $^gf$ the
  following mapping
  \begin{equation}
    \label{eq:5}
    \begin{array}{crcl}
      ^gf \colon & M & \to & A^e \\
                 & m & \mapsto & \,^g(f(m \cdot ((1*g)\otimes
                                 (1*g^{-1}))))\,.
    \end{array}
  \end{equation}
  This construction defines an action of $G$ on
  $\mathrm{Hom}_{A^e}(M,A^e)$ compatible with the natural structure of
  dg $A^e$-module, and hence an additive functor
  \begin{equation}
    \label{eq:18}
    \mathrm{Hom}_{A^e}(-,A^e) \colon \mathrm{Mod}(\Delta) \to \mathrm{Mod}(\Delta)
  \end{equation}
such that the two following diagrams are commutative (the vertical arrows
are the restriction-of-scalars functors).
  \begin{equation}
    \label{eq:16}
      \xymatrix{
        \mathrm{Mod}(\Delta) \ar[rrr]^{\mathrm{Hom}_{A^e}(-,A^e)} \ar[d] &&&
        \mathrm{Mod}(\Delta) \ar[d]\\
        \mathrm{Mod}(A^e) \ar[rrr]_{\mathrm{Hom}_{A^e}(-,A^e)} &&& \mathrm{Mod}(A^e)
  }
  \end{equation}
\begin{equation}
  \label{d1}
  \xymatrix{
    \mathcal{D}(\Delta) \ar[rrr]^{\mathrm{RHom}_{A^e}(-,A^e)} \ar[d] &&&
    \mathcal{D}(\Delta) \ar[d]\\
    \mathcal{D}(A^e) \ar[rrr]_{\mathrm{RHom}_{A^e}(-,A^e)} &&& \mathcal{D}(A^e)
  }
\end{equation}
\end{lem}
\begin{proof}
  The part of the statement before \eqref{eq:16} follows from a direct
  computation. In particular, \eqref{eq:16} is commutative.  Following
  Lemma~\ref{sec:dg-algebra-delta} (part (4)), the vertical arrows in
  \eqref{eq:16} preserve cofibrant objects. The diagram \eqref{d1}
  therefore arises from \eqref{eq:16}.
\end{proof}

Note that, given a quasi-isomorphism $X\to M$ in $\mathcal{D}(\Delta)$
such that $X$ is cofibrant in $\mathrm{Mod}(A^e)$, the image of $M$
under
$\mathrm{RHom}_{A^e}(-,A^e) \colon \mathcal{D}(\Delta) \to
\mathcal{D}(\Delta)$ is isomorphic to
$\mathrm{Hom}_{A^e}(X,A^e)$. Indeed, there exists a cofibrant
resolution $X'\to X$ in $\mathrm{Mod}(\Delta)$. Therefore, the image
of $M$ under
$\mathrm{RHom}_{A^e}(-,A^e)\colon \mathcal{D}(\Delta) \to
\mathcal{D}(\Delta)$ is isomorphic to
$\mathrm{Hom}_{A^e}(X',A^e)$. Now, $X'$ is cofibrant in
$\mathrm{Mod}(A^e)$ because so it is in $\mathrm{Mod}(\Delta)$;
actually, it follows from \cite[Proposition 2.5]{MR2809906} and from
the isomorphism \eqref{eq:82} that a dg $\Delta$-module is
homotopically projective if and only if it is so as a dg
$A^e$-module. The quasi-isomorphism $X'\to X$ between cofibrant
objects in $\mathrm{Mod}(A^e)$ therefore induces a quasi-isomorphism
$\mathrm{Hom}_{A^e}(X,A^e)\to \mathrm{Hom}_{A^e}(X',A^e)$ in
$\mathrm{Mod}(A^e)$, and hence also in $\mathrm{Mod}(\Delta)$. In
particular, there is no risk of misunderstanding in using the same
notation, $\mathrm{RHom}_{A^e}(-,A^e)$, for both functors
$\mathcal{D}(A^e)\to \mathcal{D}(A^e)$ and
$\mathcal{D}(\Delta) \to \mathcal{D}(\Delta)$.

It is now possible to establish the relationship between
$\mathrm{RHom}_{A^e}(A,A^e)$ and $\mathrm{RHom}_{\Lambda^e}(\Lambda,\Lambda^e)$
announced at the beginning of the subsection.
\begin{prop}
  \label{sec:inter-betw-dual-2}
  Assume setting~\ref{sec:main-results-1-bis}. The following diagram, where
  the top horizontal arrow is given by the top horizontal one of
  \eqref{d1}, is commutative up to an isomorphism of functors.
  \begin{equation}
    \label{d3}
    \xymatrix{
      \mathcal{D}(\Delta)\ar[rrr]^{\mathrm{RHom}_{A^e}(-,A^e)}
      \ar[d]_{-\otimes^{\mathrm L}_\Delta\Lambda^e} &&& 
      \mathcal{D}(\Delta) \ar[d]^{- \otimes^{\mathrm L}_\Delta \Lambda^e} \\
      \mathcal{D}(\Lambda^e) \ar[rrr]_{\mathrm{RHom}_{\Lambda^e}(-,\Lambda^e)}
      &&& \mathcal{D}(\Lambda^e)&
    }
  \end{equation}
  Consequently, there are isomorphisms in $\mathcal{D}(\Lambda^e)$,
  \begin{equation}
    \label{eq:31}
  \mathrm{RHom}_{A^e}(A,A^e)*G \simeq 
  \mathrm{RHom}_{A^e}(A,A^e)\otimes^{\mathrm L}_\Delta \Lambda^e \simeq 
  \mathrm{RHom}_{\Lambda^e}(\Lambda,\Lambda^e)\,.
  \end{equation}
\end{prop}
\begin{proof}
  For the ease of reading, $A$ and ${\mathbbm{k}} G$ are identified
  with their respective canonical images in $\Lambda$ and so are $A^e$
  and $({\mathbbm{k}} G)^e$ in $\Lambda^e$. Unless otherwise
  specified, ``$\cdot$'' denotes the product in $\Lambda^e$.  For all
  $M\in \mathrm{Mod}(\Delta)$, there is an isomorphism of dg
  $\Lambda^e$-modules,
  \begin{equation}
    \label{eq:38}
    \Phi_M \colon  \mathrm{Hom}_{A^e}(M,A^e) \otimes_\Delta \Lambda^e  \to 
    \mathrm{Hom}_\Delta(M,\Lambda^e)
  \end{equation}
  such that $\Phi_M(f \otimes((a*g)\otimes (b*h)))$ is given as follows,
  for all $f\in \mathrm{Hom}_{A^e}(M,A^e)$, $a,b\in A$ and $g,h\in G$,
  \[
  \begin{array}{rcl}
    M & \to & \Lambda^e \\
    m & \mapsto & \pm \sum_{k\in G} ((b*h) \otimes (a*g)) \cdot f(m\cdot
                  (k \otimes k ^{-1})) \cdot (k^{-1} \otimes k)\,,
  \end{array}
  \]
  where the sign is
  $(-1)^{\mathrm{deg}(b) \cdot(\mathrm{deg}(a) + \mathrm{deg}(f)) +\mathrm{deg}(a) \cdot \mathrm{deg}(f)}$.  Here
  is a description of $\Phi_M^{-1}$. Given
  $\varphi \in \mathrm{Hom}_\Delta(M,\Lambda^e)$ and $m\in M$, there is a
  unique decomposition
  $\varphi (m) = \sum_{s,t\in G} \varphi_{s,t}(m) \cdot (s \otimes t)$
  in $\Lambda^e$, where $\varphi_{s,t}(m)\in A^e$;
  in particular, for all $s,t\in G$,
  \begin{itemize}
  \item $\varphi_{s,t}(m \cdot (a \otimes b)) = \varphi_{s,t}(m) \cdot
    (\,^sa \otimes \,^{t^{-1}}b)$, for all $a,b\in A$, and
  \item $\varphi_{s,t}(m \cdot (k \otimes k^{-1})) = \varphi_{sk^{-1},
      kt}(m)$, for all $k\in G$;
  \end{itemize}
  hence, the following well-defined mapping lies in
  $\mathrm{Hom}_{A^e}(M,A^e)$, for all $g\in G$,
  \[
  \begin{array}{crcl}
    \varphi^{(g)} \colon & M & \to & A^e \\
    & m & \mapsto & (g^{-1} \otimes 1) \cdot \varphi_{g,e}(m) \cdot
                    (g\otimes 1)\,.
  \end{array}
  \]
  With this notation, $\Phi_M^{-1}$ is given by
  $\varphi \mapsto \sum_{g\in G} \varphi^{(g)} \otimes (g\otimes 1)$.

  This construction is natural. Since, moreover,
  $\mathrm{Hom}_\Delta(-,\Lambda^e) \cong \mathrm{Hom}_{\Lambda^e}(-\otimes_\Delta
  \Lambda^e, \Lambda^e)$, there is a functorial isomorphism of
  functors from $\mathrm{Mod}(A^e)$ to $\mathrm{Mod}(\Lambda^e)$,
  \begin{equation}
    \label{eq:28}
    \mathrm{Hom}_{A^e}(-,A^e) \otimes_\Delta \Lambda^e \xrightarrow{\cong}
    \mathrm{Hom}_{\Lambda^e}(- \otimes_\Delta \Lambda^e, \Lambda^e)\,.
  \end{equation}
  Now, recall that $\Lambda^e\simeq \Delta^{\mathrm{Card}(G)}$ both as dg
  $\Delta$-modules and as dg $\Delta^{\mathrm{op}}$-modules. Hence, deriving
  \eqref{eq:28} yields that \eqref{d3} commutes up to an isomorphism
  of functors, which entails \eqref{eq:31}, see
  subsection~\ref{sec:s2.1} for the first isomorphism of
  \eqref{eq:31}, and see the isomorphism \eqref{eq:7} with $M=A$ for
  the second one.
\end{proof}

\subsection{Main theorem}
\label{sec:s2.4}

Now it is possible to prove Theorem~\ref{sec:main-results-3}. Here is
a restatement of it.
\begin{thm}
  \label{sec:main-theorem-1}
  Assume setting~\ref{sec:main-results-1-bis}. Assume that
  $\mathrm{char}({\mathbbm{k}})\nmid\mathrm{Card}(G)$, that $A$ is homologically smooth and that
  $\Sigma^d\mathrm{RHom}_{A^e}(A,A^e)\simeq A$ in $\mathcal{D}(\Delta)$, then $A$ and
  $\Lambda$ are $d$-Calabi-Yau.
\end{thm}
\begin{proof}
  By restriction of scalars, $\Sigma^d\mathrm{RHom}_{A^e}(A,A^e)$ is
  isomorphic to $A$ in $\mathcal{D}(A^e)$ because this is the case in
  $\mathcal{D}(\Delta)$. Hence, $A$ is $d$-Calabi-Yau. By
  Proposition~\ref{sec:homol-smoothn-a-1}, the dg algebra $\Lambda$ is
  homologically smooth. Finally, since
  $\Sigma^d\mathrm{RHom}_{A^e}(A,A^e)\simeq A$ in $\mathcal{D}(\Delta)$, it follows from
  the isomorphism \eqref{eq:7} with $M=A$ and from
  Proposition~\ref{sec:inter-betw-dual-2} that $\Lambda$ is $d$-Calabi-Yau.
\end{proof}

Here is a necessary condition for $\Lambda$ to be Calabi-Yau. No
assumption on $\mathrm{char}({\mathbbm{k}})$ is made.
\begin{prop}
  \label{sec:main-theorem}
  Assume setting~\ref{sec:main-results-1-bis}. Assume that $\mathrm{HH}^0(A)$ is a
  local ${\mathbbm{k}}$-algebra. If $\Lambda$ is $d$-Calabi-Yau, then $A$ is
  homologically smooth and there exists $g\in G$ such that, in
  $\mathcal{D}(A^e)$,
  \[
  \Sigma^d\mathrm{RHom}_{A^e}(A,A^e) \simeq A_g\,,
  \]
  where $\bullet_g$ denotes the restriction-of-scalars functor along
  the dg automorphism $A^e \to A^e,\,a\otimes b \mapsto \,^ga \otimes b$.
\end{prop}
\begin{proof}
  The dg algebra $A$ is homologically smooth by
  Proposition~\ref{sec:homol-smoothn-a-1}. Since $\Lambda$ is
  $d$-Calabi-Yau, there is an isomorphism
  $\Sigma^d\mathrm{RHom}_{\Lambda^e}(\Lambda,\Lambda^e)\xrightarrow{\sim}
  \Lambda$
  in $\mathcal{D}(\Lambda^e)$. Using Proposition~\ref{sec:inter-betw-dual-2}
  and the isomorphism \eqref{eq:7} with $M=A$ yields that
  $\Sigma^d\mathrm{RHom}_{A^e}(A,A^e)\otimes^{\mathrm L}_\Delta \Lambda^e \simeq A
  \otimes^{\mathrm L}_\Delta \Lambda^e$
  in $\mathcal{D}(\Lambda^e)$. Using \eqref{eq:7} yields that there is an
  isomorphism in $\mathcal{D}(A^e)$,
  \begin{equation}
    \label{eq:32}
    \bigoplus\limits_{g\in G} \Sigma^d \mathrm{RHom}_{A^e}(A,A^e)_g \simeq
    \bigoplus\limits_{g\in G} A_g \,.
  \end{equation}
  Note that all the summands appearing in \eqref{eq:32} lie in
  $\mathrm{per}(A^e)$ because they are direct summands of $\Lambda$ in
  $\mathcal{D}(A^e)$.

  Now, since $\mathrm{RHom}_{A^e}(-,A^e)$ is a duality on $\mathrm{per}(A^e)$, the
  endomorphism algebras of $A_g$ and $\Sigma\mathrm{RHom}_{A^e}(A,A^e)_h$, for
  $g,h\in G$, are isomorphic to $\mathcal{D}(A^e)(A,A)$ ($\simeq
  \mathrm{HH}^0(A)$).
  Since $\mathrm{HH}^0(A)$ is local, it follows from \eqref{eq:32} that
  $\Sigma^d\mathrm{RHom}_{A^e}(A,A^e) \simeq A_g$ in $\mathcal{D}(A^e)$ for some
  $g\in G$.
\end{proof}

\subsection{Skew group algebras of deformed Calabi-Yau completions}
\label{sec:smash-prod-deform}

Assume setting~\ref{sec:main-results-1-bis} and assume that $\mathrm{char}({\mathbbm{k}})$ does
not divide $\mathrm{Card}(G)$ and that $A$ is homologically smooth. Hence
$\Lambda$ is homologically smooth. This subsection investigates when
$G$ acts on the deformed Calabi-Yau completions of $A$, as introduced
in \cite{MR2795754}, and when the associated skew group algebras are
deformed Calabi-Yau completions of $\Lambda$. For this subsection,
consider a cofibrant replacement in $\mathrm{Mod}(\Delta)$
\[
r_A \colon X \to A\,,
\]
and consider a quasi-isomorphism in $\mathrm{Mod}(\Delta)$ which is a cofibrant
replacement in $\mathrm{Mod}(A^e)$,
\[
p_A \colon  \Theta_A \to \mathrm{Hom}_{A^e}(X,A^e)\,.
\]
Recall that $r_A$ is a cofibrant resolution in $\mathrm{Mod}(A^e)$ (see
subsection~\ref{sec:s2.1}) and that dg $\Delta$-modules are also viewed
as dg $A^e$-modules with a compatible action of $G$.

The quasi-isomorphism $p_A$ induces a quasi-isomorphism in $\mathrm{Mod}(\mathbbm{k})$
because $A$ is homologically smooth,
\begin{equation}
  \label{eq:33}
  \begin{array}{rcl}
    \Sigma^{d-2}X\otimes_{A^e} A & \to & \mathrm{Hom}_{A^e}(\Sigma^{d-1}\Theta_{A},\Sigma
                             A) \\
    x \otimes a & \mapsto & (\theta \mapsto (p_A(\theta))(x)\cdot a)\,,
  \end{array}
\end{equation}
where the ``$\cdot$'' denotes the structure action of $A^e$ on the
left dg $A^e$-module $A$. The (deformed) Calabi-Yau completions of $A$ are
defined as follows.
\begin{definition}[\cite{MR2795754}]
  \label{sec:skew-group-algebras}
  Let $A$ be a homologically smooth dg algebra. Let $\Theta_A$ be a
  cofibrant replacement of $\mathrm{RHom}_{A^e}(A,A^e)$ in $\mathrm{Mod}(A^e)$.
  \begin{enumerate}
  \item The \emph{$d$-Calabi-Yau completion} $\mathbf{\Pi}_d(A)$ of $A$ is the
    following tensor dg algebra,
\begin{equation}
  \label{eq:34}
  \mathbf{\Pi}_d(A) = T_A(\Sigma^{d-1}\Theta_A)\,.
\end{equation}
\item Let $D\in Z^0\mathrm{Hom}_{A^e}(\Sigma^{d-1}\Theta_A,\Sigma A)$. Denote
  by $\delta$ the differential of $\Sigma^{d-1}\Theta_A$. Define
  $\delta_D$ to be the unique square-zero skew-derivation
  of degree $1$ on the graded algebra $T_A(\Sigma^{d-1}\Theta_A)$
  whose restriction to $A$ is equal to the differential of $A$ and
  whose restriction to $\Sigma^{d-1}\Theta_A$ is given by
\begin{equation}
  \label{eq:37}
  [D,\delta]\colon \Sigma^{d-1}\Theta_A \to \Sigma A \bigoplus \Sigma^d\Theta_A\,.
\end{equation}
\item Let $c\in \mathrm{HH}_{d-2}(A)$. Let
  $D\in Z^0\mathrm{Hom}_{A^e}(\Sigma^{d-1}\Theta_A,\Sigma A)$ represent $c$
  \emph{via} \eqref{eq:33}. The \emph{deformed $d$-Calabi-Yau completion}
  $\mathbf{\Pi}_d(A,c)$ is defined by
\begin{equation}
  \label{eq:35}
  \mathbf{\Pi}_d(A,c) = (T_A(\Sigma^{d-1}\Theta_A),\delta_D)\,.
\end{equation}
  \end{enumerate}
\end{definition}
Up to quasi-isomorphisms of dg algebras, neither $\mathbf{\Pi}_d(A)$ nor
$\mathbf{\Pi}_d(A,c)$ depend on the cofibrant replacement $r_A\colon X\to A$ or
on the cocycle $D$.

\begin{thm}[\cite{MR2795754,yeung}]
  Let $A$ be a homologically smooth dg algebra. Let $d$ be an
  integer. The Calabi-Yau completion $\mathbf{\Pi}_d(A)$ is $d$-Calabi-Yau (\cite[Theorem
  4.8]{MR2795754}). For all $c\in \mathrm{HH}_{d-2}(A)$ lying in the image of
  the canonical morphism $\mathrm{HN}_{d-2}(A) \to \mathrm{HH}_{d-2}(A)$ the
  deformed Calabi-Yau completion $\mathbf{\Pi}_d(A,c)$ is $d$-Calabi-Yau (\cite[Theorem
  3.17]{yeung}).
\end{thm}

In order to relate the (deformed) Calabi-Yau completions of $A$ to those of
$\Lambda$, it is necessary to endow the former ones with suitable
actions of $G$. For this purpose, note that the mapping \eqref{eq:33}
is $G$-equivariant provided that $\Sigma^{d-2}X\otimes_{A^e}A$ and
$\mathrm{Hom}_{A^e}(\Sigma^{d-1}\Theta_A,\Sigma A)$ are endowed with the
following actions of $G$,
\begin{itemize}
\item $^g(x \otimes a) = \,^gx\otimes\,^ga$, for all $g\in G$,
  $x\in X$ and $a\in A$,
\item $(\,^gD)(\theta) = \,^g(D(\,^{g^{-1}}\theta))$, for all
  $g\in G$, $D\in \mathrm{Hom}_{A^e}(\Sigma^{d-1}\Theta_A,\Sigma A)$ and
  $\theta\in \Sigma^{d-1}\Theta_A$.
\end{itemize}
In particular, the former action induces the natural action on the
Hochschild homology $\mathrm{HH}_{d-2}(A)$ which is computed by
$X \otimes_{A^e}A$.

\begin{lem}
  \label{sec:smash-prod-deform-1}
  Assume setting~\ref{sec:main-results-1-bis}. Assume that $\mathrm{char}({\mathbbm{k}})$ does
  not divide $\mathrm{Card}(G)$ and that $A$ is homologically smooth. Let
  $\Theta_A$ be a dg $\Delta$-module which is a cofibrant replacement
  of $\mathrm{RHom}_{A^e}(A,A^e)$ in $\mathrm{Mod}(A^e)$.
  \begin{enumerate}
  \item The actions of $G$ on $A$ and on $\Theta_A$ extend uniquely to
    an action by dg automorphisms on
    $T_A(\Sigma^{d-1}\Theta_A)=\mathbf{\Pi}_d(A)$.
  \item Let $c\in \mathrm{HH}_{d-2}(A)^G$. Then, there exists
    $D\in Z^0\mathrm{Hom}_{A^e}(\Sigma^{d-1}\Theta_A,\Sigma A)$ which
    represents $c$ \emph{via} \eqref{eq:33} and such that
    $D(\,^g\theta) = \,^g(D(\theta))$ for all $g\in G$ and
    $\theta\in \Theta_A$. For any such $D$, the action of $G$ on the
    graded algebra $T_A(\Sigma^{d-1}\Theta_A)$ is an action by dg
    automorphisms on $(T_A(\Sigma^{d-1}\Theta_A),\delta_D)=\mathbf{\Pi}_d(A,c)$
    (see part (2) of Definition~\ref{sec:skew-group-algebras} for the
    definition of $\delta_D$).
  \end{enumerate}
\end{lem}
\begin{proof}
  (1) In view of \eqref{eq:36}, there is an action of $G$ on
  $(\Sigma^{d-1}\Theta_A)^{\otimes_An}$, for all
  $n\in\mathbb N\backslash\{0\}$, such that
  $^g(\theta_1\otimes\cdots\otimes \theta_n) =
  \,^g\theta_1\otimes\cdots\otimes\,^g\theta_n$
  for all $\theta_1,\ldots,\theta_n\in \Sigma^{d-1}\Theta_A$ and
  $g\in G$. These actions together with the one on $A$ do form an
  action of $G$ on $T_A(\Sigma^{d-1}\Theta_A)$ by dg
  automorphisms. The uniqueness of this action follows from
  $T_A(\Sigma^{d-1}\Theta_A)$ being generated by $\Theta_A$ over $A$
  as an algebra.

  (2) Note that
  $(H^0\mathrm{Hom}_{A^e}(\Sigma^{d-1}\Theta_A , \Sigma A))^G \cong
  H^0(\mathrm{Hom}_{A^e}(\Sigma^{d-1}\Theta_A,\Sigma A)^G)$
  because the algebra ${\mathbbm{k}} G$ is semisimple.  This explains the
  existence of $D$. In particular, both the differential of $A$ and
  the mapping \eqref{eq:37} are $G$-equivariant. Hence, so is
  $\delta_D$. Accordingly, the action of $G$ on the graded algebra
  $T_A(\Sigma^{d-1}\Theta_A)$ is an action by dg automorphisms on
  $\mathbf{\Pi}_d(A,c)$ when the latter is taken equal to
  $(T_A(\Sigma^{d-1}(\Theta_A),\delta_D)$.
\end{proof}

Lemma~\ref{sec:smash-prod-deform-1} sets the framework in which the dg
algebras $\mathbf{\Pi}_d(A)*G$ and $\mathbf{\Pi}_d(A,c)*G$ are defined. Using this, it is
possible to prove Theorem~\ref{sec:main-results-2}. This theorem is
restated below for convenience.
\begin{thm}
  \label{sec:skew-group-algebras-2}
  Assume setting~\ref{sec:main-results-1-bis}. Assume that $A$ is
homologically smooth and that $\mathrm{char}({\mathbbm{k}})\nmid\mathrm{Card}(G)$.
  \begin{enumerate}
  \item $A*G$ is homologically smooth and, up to appropriate choices
    of inverse dualising complexes of $A$ and $A*G$, the group $G$
    acts on $\mathbf{\Pi}_d(A)$ by dg automorphisms and
    $\mathbf{\Pi}_d(A)*G\simeq \mathbf{\Pi}_d(A*G)$ as dg algebras. In particular,
    $\mathbf{\Pi}_d(A)*G$ is $d$-Calabi-Yau.
  \item Let $c\in \mathrm{HH}_{d-2}(A)^G$. Let $\overline c$ be the image of
    $c$ under the mapping $\mathrm{HH}_{d-2}(A)\to \mathrm{HH}_{d-2}(A*G)$ induced by
    $A \to A*G,\,a \mapsto a*e$. Up to appropriate choices as in (1)
    and up to appropriate choices of representatives of $c$ and
    $\overline c$, the group $G$ acts on $\mathbf{\Pi}_d(A,c)$ by dg
    automorphisms and $\mathbf{\Pi}_d(A,c)*G \simeq \mathbf{\Pi}_d(A*G,\overline c)$ as dg
    algebras. In particular, if $c$ lifts to $\mathrm{HN}_{d-2}(A)$
    (so that $\mathbf{\Pi}_d(A,c)$ is $d$-Calabi-Yau), then $\mathbf{\Pi}_d(A,c)*G$ is
    $d$-Calabi-Yau.
  \end{enumerate}  
\end{thm}
\begin{proof}
  It follows from Proposition~\ref{sec:homol-smoothn-a-1} that
  $\Lambda$ is homologically smooth.
  
  (1) Here, $\mathbf{\Pi}_d(A)$ is taken equal to the dg algebra
  $T_A(\Sigma^{d-1}\Theta_A)$ endowed with the action of $G$
  introduced in Lemma~\ref{sec:smash-prod-deform-1}.  Note that
  $\Theta_A*G$ is a cofibrant replacement of
  $\mathrm{RHom}_{\Lambda^e}(\Lambda,\Lambda^e)$ in $\mathrm{Mod}(\Lambda^e)$ (see
  subsection~\ref{sec:s2.1} and
  Proposition~\ref{sec:inter-betw-dual-2}). Hence, $\mathbf{\Pi}_d(A*G)$ may be
  taken equal to the tensor dg algebra
  $T_\Lambda(\Sigma^{d-1}\Theta_A*G)$. To avoid confusion, the
  piece of notation ``$*$'' is kept for elementary tensors in $A*G$ and
  $\Theta_A*G$ but not for those in $T_A(\Sigma^{d-1}\Theta_A)*G$. The
  following map is therefore a morphism of dg algebras,
  \begin{equation}
    \label{eq:41}
    \begin{array}{rcll}
      \mathbf{\Pi}_d(A)* G & \longrightarrow & \mathbf{\Pi}_d(A*G)\\
      a\otimes g & \longmapsto & a*g\in A*G\\
      (\phi_1\otimes \cdots\otimes \phi_n)\otimes g & \longmapsto &
                                                                    (\phi_1*e) \otimes \cdots \otimes (\phi_{n-1}*e)\otimes (\phi_n* g) &,
    \end{array}
  \end{equation}
  where $a\in A$, $\phi_1,\ldots,\phi_n\in\Sigma^{d-1}\Theta_A$ and $g\in
  G$. It is an isomorphism, its inverse maps $a*g\in A*G$ onto
  $a\otimes g$, for all $a\in A$ and $g\in G$, and it maps $(\phi_1 *
  g_1) \otimes \cdots \otimes (\phi_n *g_n)$ onto
  \[
  \left(\phi_1\otimes
    \,^{g_1}\phi_2\otimes\cdots \otimes \,^{g_1\cdots
      g_l}\phi_{l+1}\otimes\cdots\otimes \,^{g_1\cdots
      g_{n-2}}\phi_{n-1}\otimes \,^{g_1\cdots
      g_{n-1}}\phi_{n}\right)\otimes g_1\cdots
  g_n\,
  \]
  for all $(\phi_1*g_1),\ldots,(\phi_n*g_n)\in \Sigma^{d-1}\Theta_A*G$.

  (2) In view of the functoriality of the Hochschild homology and of
  the negative cyclic homology, if $c$ lifts to $\mathrm{HN}_{d-2}(A)$
  then $\overline c$ lifts to $\mathrm{HN}_{d-2}(A*G)$. It is hence
  sufficient to prove the first statement of (2). Thanks to the remark
  following Lemma~\ref{sec:inter-betw-dual}, it is possible to choose
  $p_A$ such that it is, in addition, a cofibrant resolution in
  $\mathrm{Mod}(\Delta)$. Here $\mathbf{\Pi}_d(A,c)$ is taken equal to
  $(T_A(\Sigma^{d-1}\Theta_A),\delta_D)$ where $D$ is as in
  Lemma~\ref{sec:smash-prod-deform-1}. Thus, $\mathbf{\Pi}_d(A,c)$ is endowed
  with an action of $G$ by dg automorphisms.

  It is necessary to first describe $\mathbf{\Pi}_d(\Lambda,\overline c)$ in
  the same way as $\mathbf{\Pi}_d(A,c)$. The following composite mapping is a
  cofibrant resolution in $\mathrm{Mod}(\Lambda^e)$,
  \[
  r_\Lambda \colon X\otimes_\Delta \Lambda^e\xrightarrow{ r_A \otimes
    \Lambda^e} A \otimes_\Delta \Lambda^e \xrightarrow{~\eqref{eq:7}} \Lambda\,.
  \]
  The following composite mapping is a cofibrant resolution in
  $\mathrm{Mod}(\Lambda^e)$,
  \[
  p_\Lambda \colon \Theta_A \otimes_\Delta \Lambda^e
  \xrightarrow{p_A \otimes \Lambda^e}
  \mathrm{Hom}_{A^e}(X,A^e) \otimes_\Delta \Lambda^e
  \underset{~\eqref{eq:38}}{\xrightarrow{ \frac{1}{\mathrm{Card}(G)} \cdot
      \Phi_X}}
  \mathrm{Hom}_\Delta(X,\Lambda^e)
  \xrightarrow{\cong}
  \mathrm{Hom}_{\Lambda^e}(X\otimes_\Delta \Lambda^e,\Lambda^e)\,.
  \]
  Hence, $\mathbf{\Pi}_d(\Lambda, \overline c)$ may be taken equal to the
  graded algebra $T_\Lambda(\Sigma^{d-1}\Theta_A \otimes_\Delta
  \Lambda^e)$ with the differential $(\delta \otimes_\Delta
  \Lambda^e)_{D'}$, where $D'$ is any $0$-cocycle in
  $\mathrm{Hom}_{\Lambda^e}(\Sigma^{d-1}\Theta_A \otimes_\Delta
  \Lambda^e,\Sigma \Lambda)$ which represents $\overline c$ \emph{via}
  the following quasi-isomorphism
  \begin{equation}
    \label{eq:39}
    \begin{array}{rcl}
      \Sigma^{d-2}(X\otimes_\Delta \Lambda^e)\otimes_{\Lambda^e} \Lambda
      & \to &
              \mathrm{Hom}_{\Lambda^e}(\Sigma^{d-1}\Theta_A\otimes_\Delta
              \Lambda^e,\Sigma\Lambda) \\
      \alpha \otimes \lambda & \mapsto &
                                         (\theta \mapsto
                                         (p_\Lambda(\theta))(\alpha)
                                         \cdot \lambda)
    \end{array}
  \end{equation}
  where the ``$\cdot$'' stands for the structure action of $\Lambda^e$
  on $\Lambda$.  Recall that $\delta$ is the differential of
  $\Sigma^{d-1}\Theta_A$.

  In order to prove (2), it is hence sufficient to prove the following
  claim.

  \emph{Claim 1 - $D'$ may be chosen such that, in addition,
    \eqref{eq:41} is an isomorphism of dg algebras from
    $(T_A(\Sigma^{d-1}\Theta_A),\delta_D)*G$ to
    $(T_\Lambda(\Sigma^{d-1}\Theta_A \otimes_\Delta
    \Lambda^e),(\delta\otimes_\Delta \Lambda^e)_{D'})$.}

  Note that, in view of the definition of $p_\lambda$, and up to the
  two canonical identifications
  \begin{equation}
    \label{eq:43}
    (X\otimes_\Delta \Lambda^e)
    \otimes_{\Lambda^e}\Lambda \cong X\otimes_\Delta \Lambda
  \end{equation}
  and
  \begin{equation}
    \label{eq:42}
  \mathrm{Hom}_{\Lambda^e}(\Sigma^{d-1}\Theta_A\otimes_\Delta
  \Lambda^e,\Sigma \Lambda) \cong
  \mathrm{Hom}_\Delta(\Sigma^{d-1}\Theta_A,\Sigma \Lambda)\,,
  \end{equation}
  the quasi-isomorphism \eqref{eq:39} identifies with the following
  one
  \begin{equation}
    \label{eq:40}
    \begin{array}{rcl}
      \Sigma^{d-2}X\otimes_\Delta \Lambda
      & \to &
              \mathrm{Hom}_\Delta(\Sigma^{d-1}\Theta_A,\Sigma\Lambda) \\
      x \otimes \lambda & \mapsto &
                                    (\theta \mapsto
                                    \frac{1}{\mathrm{Card}(G)}\cdot
                                    \sum_{k \in G}
                                    (p_A(\theta))(\,^{k^{-1}}x)\cdot
                                    (k^{-1}\otimes k)
                                    \cdot \lambda)\,,
    \end{array}
  \end{equation}
  where, for simplicity, $A^e$ and $({\mathbbm{k}} G)^e$ are viewed as subspaces
  of $\Lambda^e$ in a natural way.

  Denote by $\overline D$ the following cocycle
  \[
  \begin{array}{rcl}
    \overline D \colon \Sigma^{d-1}\Theta_A & \to & \Sigma \Lambda
    \\
    \theta & \mapsto & D(\theta) *e\,.
  \end{array}
  \]
  Then, the following diagram where the rightmost vertical arrow is
  the pre-image of $\overline D$ under \eqref{eq:42} is commutative,
  \[
    \xymatrix{ \Sigma^{d-1}\Theta_A * G \ar[r]^{~\eqref{eq:41}}
      \ar[d]_{\overline D}& \Sigma^{d-1}
      \Theta_A\otimes_\Delta \Lambda^e \ar[d] \\
      \Sigma A*G \ar[r]_{~\eqref{eq:41}} & \Sigma \Lambda }
  \]
  Hence, in order to prove claim 1, it is sufficient to prove the
  following claim.

  \emph{Claim 2 - $D'$ may be chosen such that its image under
    \eqref{eq:42} is equal to $\overline D$.}

  In order to prove this claim, consider the following diagram
  \begin{equation}
    \label{eq:44}
    \xymatrix{
      (\Sigma^{d-2}X \otimes_{A^e} A)^G \ar[r] \ar[d]_{~\eqref{eq:33}}
      & \Sigma^{d-2} X \otimes_\Delta \Lambda \ar[d]^{~\eqref{eq:40}}
      \ar[r]^{~\eqref{eq:43}} & (\Sigma^{d-2} X\otimes_\Delta
      \Lambda^e) \otimes_{\Lambda^e} \Lambda \ar[d]^{p_\Lambda} \\
      \mathrm{Hom}_\Delta(\Sigma^{d-1}\Theta_A, \Sigma A) \ar[r] &
      \mathrm{Hom}_\Delta(\Sigma^{d-1}\Theta_A,\Sigma\Lambda)
      \ar[r]_{~\eqref{eq:42}} &
      \mathrm{Hom}_{\Lambda^e}(\Sigma^{d-1}\Theta_A\otimes_\Delta \Lambda^e,
      \Sigma\Lambda) 
    }
  \end{equation}
  where
  \begin{itemize}
  \item the codomain of the $G$-invariant part of \eqref{eq:33} is
    indeed equal to $\mathrm{Hom}_\Delta(\Sigma^{d-1}\Theta_A,\Sigma A)$ as a
    subcomplex of $\mathrm{Hom}_{A^e}(\Sigma^{d-1}\Theta_A,\Sigma A)$,
  \item the unlabelled top horizontal arrow is induced by the mapping
    $x \otimes a \mapsto x \otimes(a*e)$, and
  \item the unlabelled bottom horizontal arrow is given by the
    composition with the embedding $\Sigma A \to \Sigma
    \Lambda,\,a\mapsto a*e$.
  \end{itemize}
  Note that,
  \begin{itemize}
  \item the mapping $\mathrm{HH}_{d-2}(A)^G \to \mathrm{HH}_{d-2}(A*G)$ which assigns
    $\overline c$ to $c$ is induced by the composition of the
    top horizontal arrows of \eqref{eq:44},
  \item $\overline D$ is the image of $D$ under the unlabelled bottom
    horizontal arrow of \eqref{eq:44},
  \item $c$ (or $\overline c$) is represented by any cocycle in
    $(\Sigma^{d-2}X \otimes_{A^e}A)^G$ (or, in $(\Sigma^{d-2}X
    \otimes_{\Delta} \Lambda^e) \otimes_{\Lambda^e} \Lambda$) whose
    image under the leftmost (or, rightmost) vertical arrow of
    \eqref{eq:44} is cohomologous to $D$ (or, to $D'$, respectively),
    when viewing $\mathrm{Hom}_\Delta(\Sigma^{d-1}\Theta_A,\Sigma A)$ as a
    subcomplex of $\mathrm{Hom}_{A^e}(\Sigma^{d-1}\Theta_A,\Sigma A)$.
  \end{itemize}
  Hence, in order to prove claim 2, it is sufficient to prove the
  following claim.

  \emph{Claim 3 - The diagram \eqref{eq:44} is commutative.}

  Here is a proof of this claim. Given that \eqref{eq:40} is obtained
  upon composing \eqref{eq:43}, \eqref{eq:39}, and \eqref{eq:42}, the
  rightmost square is commutative. Let
  $\sum_i x_i\otimes a_i \in (\Sigma^{d-2}X\otimes_{A^e}A)^G$. Hence,
  for all $k\in G$,
  \begin{equation}
    \label{eq:45}
    \sum_i x_i\otimes a_i = \sum_i \,^kx_i \otimes\,^ka_i\,.
  \end{equation}
  On one hand, the image of $\sum_i x_i\otimes a_i$ under the leftmost
  vertical arrow of \eqref{eq:44} is as follows
  \[
  \begin{array}{rcl}
    \Sigma^{d-1}\Theta_A & \to & \Sigma A\\
    \theta & \mapsto & \sum_i(p_A(\theta)(x_i))\cdot a_i\,,
  \end{array}
  \]
  which has image under the unlabelled bottom horizontal arrow as
  follows
  \begin{equation}
    \label{eq:46}
    \begin{array}{rcl}
      \Sigma^{d-1} \Theta_A & \to & \Sigma \Lambda \\
      \theta & \mapsto & \sum_i
                         \underset{=(p_A(\theta)(x_i))\cdot (a_i*e)\,.}
                         {\underbrace{((p_A(\theta)(x_i)\cdot a_i))*e}}
    \end{array}
  \end{equation}
  On the other hand, the image of $\sum_i x_i\otimes a_i$ under the
  unlabelled top horizontal arrow of \eqref{eq:44} is equal to
  $\sum_i x_i \otimes (a_i *e)$, which has image under the middle
  vertical arrow as follows
  \begin{equation}
    \label{eq:47}
    \begin{array}{rcl}
      \Sigma^{d-1}\Theta_A & \to & \Sigma\Lambda \\
      \theta & \mapsto & \frac{1}{\mathrm{Card}(G)}
                         \cdot
                         \sum_{k\in G,\,i}
                         (p_A(\theta)(\,^{k^{-1}}x_i)) \cdot
                         \underset{=\,^{k^{-1}}a_i*e}{
                         \underbrace{(k^{-1}\otimes k) \cdot (a_i*e)}\,.
                         }
    \end{array}
  \end{equation}
  In view of \eqref{eq:45}, the mappings \eqref{eq:46} and
  \eqref{eq:47} are equal.  Thus \eqref{eq:44} is commutative. This
  finishes proving claim 3, claim 2, claim 1, and, finally, (2).
\end{proof}

\section{Application to Ginzburg dg algebras}
\label{sec:smash-prod-ginzb}

Assume setting~\ref{sec:appl-ginzb-dg}. It is restated below for
convenience.
\begin{setting}
  \label{sec:appl-ginzb-dg-bis}
  Let $Q$ be a finite graded ${\mathbbm{k}}$-quiver. Let $G$ be a finite group
  such that $\mathrm{char}({\mathbbm{k}})\nmid \mathrm{Card}(G)$, acting on ${\mathbbm{k}} Q$ by degree
  preserving automorphisms in such a way that both the set of
  (idempotents associated with) vertices and the vector space
  generated by the arrows of $Q$ are stabilised by the action. Let $W$
  be a homogeneous potential of degree $d-3$ on $Q$.
\end{setting}

The purpose of this section is to prove
Corollary~\ref{sec:main-results-4} by applying section~\ref{sec:s2} to
Ginzburg dg algebras. Subsection~\ref{sec:reminder-ginzburg-dg} gives
a reminder on Ginzburg dg
algebras. Subsection~\ref{sec:equiv-non-comm} proves that cyclic
derivations are $G$-equivariant. This is used in
subsection~\ref{sec:group-acti-mathc} to prove that $G$ acts on
$\mathcal A(Q,W)$ by dg automorphisms if $W$ is $G$-invariant up to
cyclic permutation. Subsection~\ref{sec:smash-prod-ginzb-1} proves
Corollary~\ref{sec:main-results-4}.

\subsection{Reminder on Ginzburg dg algebras}
\label{sec:reminder-ginzburg-dg}

These were introduced in \cite[4.2]{G06}. The more general definition
used here follows \cite[6.3]{MR2795754}.

Let $Q$ be a finite \emph{graded ${\mathbbm{k}}$-quiver}, that is, $Q$ is a
quiver (or, an oriented graph) with (finite) set of vertices denoted
by $Q_0$ and (finite) set of arrows denoted by $Q_1$, and the
${\mathbbm{k}} Q_0$-bimodule ${\mathbbm{k}} Q_1$ spanned by $Q_1$ has a
$\mathbbm{Z}$-grading. Without loss of generality, assume that every
arrow of $Q$ is homogeneous. The path algebra ${\mathbbm{k}} Q$ is hence
$\mathbb Z$-graded, idempotents $e_x$ of vertices $x$ have degree
$0$. With this grading, ${\mathbbm{k}} Q$ is also viewed as a dg algebra with
zero differential.

The space of \emph{potentials} on $Q$ is the $\mathbb Z$-graded vector
space,
\begin{equation}
  \label{eq:59}
\frac{{\mathbbm{k}} Q\otimes_{{\mathbbm{k}}{Q_0}^e}{\mathbbm{k}} Q_0}
{\left\langle
  uv\otimes e_x-(-1)^{\mathrm{deg}(u).\mathrm{deg}(v)}vu\otimes e_y\ |\
  u\in e_x{\mathbbm{k}} Qe_y,\,v\in e_y{\mathbbm{k}} Qe_x\ \text{homogeneous, $x,y\in Q_0$}\right\rangle}\,.
\end{equation}
Hence, potentials are linear combinations of oriented cycles in $Q$,
each of which is considered up to cyclic permutation with Koszul-type
signs. Note that ${\mathbbm{k}} Q$ has a cofibrant resolution in $\mathrm{Mod}({\mathbbm{k}} Q^e)$
by the cone of
\[
\begin{array}{ccc}
{\mathbbm{k}} Q\otimes_{{\mathbbm{k}} Q_0} {\mathbbm{k}} Q_1 \otimes_{{\mathbbm{k}} Q_0}{\mathbbm{k}} Q & \to & {\mathbbm{k}} Q
                                                         \otimes_{{\mathbbm{k}} Q_0} {\mathbbm{k}} Q  \\
  u \otimes a \otimes v & \mapsto & ua\otimes v -u \otimes av
\end{array}
\]
so that the Hochschild homology of ${\mathbbm{k}} Q$ is computed by the cone of
\[
\begin{array}{ccc}
  ({\mathbbm{k}} Q_1\otimes_{{\mathbbm{k}} Q_0} {\mathbbm{k}} Q)  \otimes_{{\mathbbm{k}}{Q_0}^e}{\mathbbm{k}} Q_0 & \to & {\mathbbm{k}} Q
                                                                  \otimes_{{\mathbbm{k}}{Q_0}^e} {\mathbbm{k}} Q_0  \\
  a \otimes u \otimes 1 & \mapsto & (-1)^{\mathrm{deg}(a)\cdot \mathrm{deg}(u)} ua
                                    \otimes 1\,.
\end{array}
\]
Hence, a homogeneous potential of degree $n$ on $Q$ may be viewed as
an element of $\mathrm{HH}_n({\mathbbm{k}} Q)$.

Let $d\geqslant 3$ be an integer. Let $W$ be a homogeneous potential
on $Q$ of degree $d-3$. The Ginzburg dg algebra $\mathcal A(Q,W)$ is
defined as follows. As a graded algebra, $\mathcal A(Q,W)$ equals
${\mathbbm{k}} \widetilde Q$ where $\widetilde Q$ is the following graded
${\mathbbm{k}}$-quiver,
\begin{itemize}
\item $\widetilde Q$ and $Q$ have the same vertices,
\item as a graded ${\mathbbm{k}} Q_0$-bimodule, ${\mathbbm{k}} \widetilde Q_1$ is equal to
  \begin{equation}
    \label{eq:48}
  {\mathbbm{k}} \widetilde Q_1 = {\mathbbm{k}} Q_1 \bigoplus \Sigma^{d-2}\mathrm{Hom}_{{\mathbbm{k}}{Q_0}^e}({\mathbbm{k}}
  Q_1,{\mathbbm{k}}{Q_0}^e) \bigoplus \Sigma^{d-1} {\mathbbm{k}} Q_0\,,
\end{equation}
given an arrow $a\colon x\to y$ in $Q$ denote by $a^*$ the element of
$\mathrm{Hom}_{{\mathbbm{k}}{Q_0}^e}({\mathbbm{k}} Q_1,{\mathbbm{k}}{Q_0}^e)$ defined by its behaviour on arrows
as follows
\begin{equation}
  \label{eq:85}
  a^* \colon b\in Q_1 \mapsto \left\{
    \begin{array}{ll}
      e_x \otimes e_y & \text{if $b=a$} \\
      0 & \text{otherwise,}
    \end{array}\right.
\end{equation}
and given $i\in Q_0$ denote by $c_i$ the corresponding element of
$\Sigma^{d-1}{\mathbbm{k}} Q_0$; this defines a bijective ${\mathbbm{k}}$-linear mapping
\begin{equation}
  \label{eq:88}
  \begin{array}{rcl}
    {\mathbbm{k}} Q_1 & \longrightarrow & \mathrm{Hom}_{{\mathbbm{k}}{Q_0}^e}({\mathbbm{k}} Q_1,{\mathbbm{k}}{Q_0}^e) \\
    a \in Q_1 & \longmapsto & a^*\,.
  \end{array}
\end{equation}
\end{itemize}
Note that the grading on $\mathcal A(Q,W)$ is determined by the
grading on ${\mathbbm{k}} Q$.

The differential on $\mathcal A(Q,W)$ (of degree $+1$) is uniquely
determined by the following rules:
\begin{itemize}
\item it vanishes on ${\mathbbm{k}} Q$,
\item it maps $a^*\colon y\to x$ to the cyclic derivative
  $\partial_aW$, for all arrows $a\colon x \to y$ in $Q_1$,
\item it maps $c_i\colon i\to i$ to $\sum\limits_{a\colon i\to
    \cdot}aa^*-\sum\limits_{a\colon \cdot\to i}a^*a$, where the first
  sum runs over all the arrows of $Q$ starting in $i$ and the second
  sum runs over all the arrows of $Q$ arriving in $i$.
\end{itemize}
Recall that, when $a\in Q_1$, the cyclic derivation $\partial_a$ is
the linear mapping which takes a path $p$ to
$\partial_ap=\sum\limits_{p=p_1ap_2}(-1)^{\mathrm{deg}(p_1)\cdot\mathrm{deg}(ap_2)}p_2p_1$,
where the sum runs over all decompositions $p=p_1ap_2$ with paths
$p_1$ and $p_2$.

By \cite[Subsections 6.1 and 6.2, and Theorem A.12]{MR2795754},
$\mathcal A(Q,W)$ is $d$-Calabi-Yau and there is a quasi-isomorphism of dg
algebras
\begin{equation}
  \label{eq:49}
  \mathbf{\Pi}_d({\mathbbm{k}} Q,c) \xrightarrow\sim \mathcal A(Q,W)\,,
\end{equation}
where $c\in\mathrm{HH}_{d-2}({\mathbbm{k}} Q)$ is the image of $W$ under Connes' boundary
in Hochschild homology
$\mathrm{HH}_{d-3}({\mathbbm{k}} Q)\xrightarrow{B}\mathrm{HH}_{d-2}({\mathbbm{k}} Q)$.

\subsection{Equivariant non commutative differential calculus on
  quiver algebras}
\label{sec:equiv-non-comm}

Assume setting~\ref{sec:appl-ginzb-dg-bis}. In this subsection, tensor
products, quotients and $\mathrm{Hom}$-spaces involving $G$-modules are
considered as $G$-modules in the classical sense of group
representations.  Denote by $\mathrm{DR}({\mathbbm{k}} Q)$ the $G$-module of potentials
on ${\mathbbm{k}} Q$. The following lemma is essential for
subsection~\ref{sec:group-acti-mathc}.
\begin{lem}
  \label{sec:equiv-non-comm-1}
  Assume setting~\ref{sec:appl-ginzb-dg-bis}. The following ${\mathbbm{k}}$-linear
  mapping is $G$-equivariant
    \begin{equation}
  \label{eq:86}
  \begin{array}{ccc} \mathrm{Hom}_{{\mathbbm{k}}{Q_0}^e}({\mathbbm{k}} Q_1,{\mathbbm{k}}{Q_0}^e) \otimes \mathrm{DR}({\mathbbm{k}}
Q)& \longrightarrow & {\mathbbm{k}} Q \\ \varphi \otimes w & \longmapsto &
d_{\mathcal A(Q,w)}(\varphi)\,,
  \end{array}
\end{equation} where $d_{\mathcal A(Q,w)}$ denotes the differential of
the dg algebra $\mathcal A(Q,w)$, for all potentials $w$ on $Q$.
\end{lem}
\begin{proof}
  The proof uses the interpretation (\cite[Subsections 3.5 and
  3.6]{G06}) of cyclic derivatives as reduced contractions
  (\cite[Section 2]{CEG07}).  Recall that the tensor product of two
  homogeneous ${\mathbbm{k}}$-linear mappings $f$ and $f'$ is defined using the
  Koszul-sign rule, that is
  $(f\otimes f')(u\otimes v) = (-1)^{\mathrm{deg}(f') \cdot \mathrm{deg}(u)}f(u)
  \otimes f'(v)$
  for all homogeneous $u,v$ such that $f(u)$ and $f'(v)$ are defined.
  Denote by $\Delta_{{\mathbbm{k}} Q}$ the natural coassociative and counital
  comultiplication of ${\mathbbm{k}} Q$, viewed as $T_{{\mathbbm{k}} Q_0}({\mathbbm{k}} Q_1)$. Denote
  by $p$ the natural projection from ${\mathbbm{k}} Q$ to ${\mathbbm{k}} Q_1$. Denote by
  $-^\diamond$ the ${\mathbbm{k}}$-linear mapping from ${\mathbbm{k}} Q \otimes {\mathbbm{k}} Q$ to
  ${\mathbbm{k}} Q$ given by
  $(u\otimes v)^\diamond = (-1)^{\mathrm{deg}(u) \cdot \mathrm{deg}(v)}vu$ for all
  homogeneous $u,v\in {\mathbbm{k}} Q$.

  By considering tensors of the shape $a^*\otimes w$ where $a\in Q_1$
  and $w\in \mathrm{DR}({\mathbbm{k}} Q)$, it is elementary to check that \eqref{eq:86}
  is equal to the following mapping (for the ease of reading, no
  distinction is made between a linear combination of oriented cycles
  and its associated potential)
  \[
  \begin{array}{ccc}
    \mathrm{Hom}_{{\mathbbm{k}}{Q_0}^e}({\mathbbm{k}} Q_1, {\mathbbm{k}}{Q_0}^e) \otimes \mathrm{DR}({\mathbbm{k}} Q) & \longrightarrow & {\mathbbm{k}} Q \\
    \varphi \otimes w & \longmapsto & \left(
                                      \left(
                                      (
                                      \mathrm{Id} \otimes (\varphi
                                      \circ p) \otimes \mathrm{Id}) \circ
                                      (\Delta_{{\mathbbm{k}} Q} \otimes \mathrm{Id}) \circ
                                      \Delta_{{\mathbbm{k}} Q}
                                      \right)
                                      (w)
                                      \right)^\diamond\,.
  \end{array}
  \]
  This is $G$-equivariant, hence so is \eqref{eq:86}
\end{proof}

\begin{rem}
  The proof and the statement of Lemma~\ref{sec:equiv-non-comm-1}
  still hold true when $G$ is infinite.
\end{rem}

\begin{rem}
  \label{sec:equiv-non-comm-2}
  Let $g\in G$. It is elementary to check that \eqref{eq:86} is
  $\langle g\rangle$-equivariant if
  \begin{equation}
    \label{eq:87}
    (\forall a\in Q_1)\ \ ^ga \in \{t\cdot b\ |\ t\in
    {\mathbbm{k}}^\times,\,b\in  Q_1\}\,.
  \end{equation}
  Since $G$ is finite, one may try to reduce the proof of
  Lemma~\ref{sec:equiv-non-comm-1} to that particular case using a
  ${\mathbbm{k}}$-algebra automorphism $f \colon {\mathbbm{k}} Q \to {\mathbbm{k}} Q$ such that
  $f^{-1}(\,^gf(-))$ satisfies \eqref{eq:87} in place of $g$. The
  reduction then consists in proving that the conjugation of
  \eqref{eq:86} under $f^{-1}$ is $\langle g\rangle$-equivariant.
  However, this is not immediate because the conjugate of \eqref{eq:86}
  under $f^{-1}$ is not \eqref{eq:86}. Rather, its expression is
  governed by the cyclic chain rule (\cite[Lemma 3.9]{MR2480710}). The
  details of the reduction using this rule are left to the reader.
\end{rem}

The following example illustrates the fact that conjugating the cyclic
derivation by an automorphism does not yield the cyclic derivation.

\begin{ex}
  \label{sec:equiv-non-comm-3}
  Assume that ${\mathbbm{k}}=\mathbb C$, that $G=\{e,g\}$, that $Q$ is the quiver
  $\xymatrix{\bullet \ar@(rd,ru)_b \ar@(ld,lu)^a}$, where
  $\mathrm{deg}(a)=\mathrm{deg}(b)=0$, that $^ga= \frac{a+b}{\sqrt{2}}$ and
  $^gb = \frac{a-b}{\sqrt{2}}$, and that $W$ is the
  $\langle g\rangle$-invariant potential
  $a^2 + \,^g(a^2) = \frac{1}{2}(3a^2+ab+ba+b^2)$. In
  Remark~\ref{sec:equiv-non-comm-2}, one may assume that $f(a) = a-ib$
  and $f(b) = a+ib$. Indeed, $^g f(a)= \frac{1+i}{2\sqrt{2}}f(b)$ and
  $^g f(b) = \frac{1-i}{2\sqrt{2}}f(a)$. Denote \eqref{eq:86} by
  $\iota$ and denote its conjugate under $f^{-1}$ by $\iota'$. Note
  that $f$ acts on $\mathrm{Hom}_{{\mathbbm{k}}{Q_0}^e}({\mathbbm{k}} Q_1, {\mathbbm{k}}{Q_0}^e)$ by
  precomposition with $f^{-1}$. On one hand,
  $\iota(a^*\otimes W) = \partial_a(W) = 3a+b$. On the other hand,
  $\iota'(a^*\otimes W) = f^{-1}(\iota((a^* \circ f^{-1}) \otimes
  f(W)))$;
  note that $a^*\circ f^{-1} = \frac{1}{2}(a^*-ib^*)$ and
  $f(W) = 3a^2+iab+iba-b^2$; hence $\iota'(a^*\otimes W) = 4a+2ib$.
  Thus $\iota$ is not invariant under conjugation by $f^{-1}$.
\end{ex}

\subsection{Group actions on \texorpdfstring{$\mathcal A(Q,W)$}{A(Q,W)}}
\label{sec:group-acti-mathc}

Assume setting~\ref{sec:appl-ginzb-dg-bis} and assume that $W$ is
$G$-invariant up to cyclic permutation. Like in
subsection~\ref{sec:equiv-non-comm}, tensor products, quotients and
$\mathrm{Hom}$-spaces involving $G$-modules are considered as $G$-modules in
the classical sense of group representations. In particular,
$\mathrm{Hom}_{{{\mathbbm{k}}{Q_0}}^e}({\mathbbm{k}} Q_1, {{\mathbbm{k}}{Q_0}}^e)$ is a $G$-module. Denote by
$c$ the image of $W$ under
$B \colon \mathrm{HH}_{d-3}({\mathbbm{k}} Q) \to \mathrm{HH}_{d-2}({\mathbbm{k}} Q)$.  Then $c$ lies in
$\mathrm{HH}_{d-2}({\mathbbm{k}} Q)^G$. The purpose of this subsection is to show that
$G$ acts on $\mathcal A(Q,W)$ by algebra automorphisms, that this is
an action by dg automorphisms, and that, when $\mathbf{\Pi}_d({\mathbbm{k}} Q,c)$ is
endowed with the action of $G$ of Lemma~\ref{sec:smash-prod-deform-1},
the following quasi-isomorphism of \cite[Theorem 6.3]{MR2795754} is
$G$-equivariant
\begin{equation}
  \label{eq:62}
  \mathbf{\Pi}_d({\mathbbm{k}} Q,c) \to \mathcal A(Q,W)\,.
\end{equation}

Here is a reminder on this quasi-isomorphism. In $\mathrm{Mod}({\mathbbm{k}} Q^e)$, there
is a cofibrant replacement of ${\mathbbm{k}} Q$ equal to the cone of the
following morphism
\begin{equation}
  \label{eq:63}
  \begin{array}{ccc}
    {\mathbbm{k}} Q \otimes_{{\mathbbm{k}}{Q_0}} {{\mathbbm{k}} Q_1} \otimes_{{\mathbbm{k}}{Q_0}} {\mathbbm{k}} Q & \xrightarrow{\alpha} & {\mathbbm{k}}
                                                               Q\otimes_{{\mathbbm{k}}{Q_0}}
                                                               {\mathbbm{k}} Q \\
    1 \otimes a \otimes 1 & \mapsto & a \otimes 1 - 1 \otimes a\,.
  \end{array}
\end{equation}
Denote by $\Theta_{{\mathbbm{k}} Q}$ the cylinder of the following morphism. 
\begin{equation}
  \label{eq:64}
  \mathrm{Hom}_{{\mathbbm{k}} Q^e}({\mathbbm{k}} Q\otimes_{{\mathbbm{k}}{Q_0}} {\mathbbm{k}} Q , {\mathbbm{k}} Q^e) \xrightarrow{\alpha^t}
  \mathrm{Hom}_{{\mathbbm{k}} Q^e}( {\mathbbm{k}} Q\otimes_{{\mathbbm{k}}{Q_0}} {{\mathbbm{k}} Q_1} \otimes_{{\mathbbm{k}}{Q_0}} {\mathbbm{k}} Q, {\mathbbm{k}} Q^e)\,.
\end{equation}
Hence, $\Theta_{{\mathbbm{k}} Q}$ is a cofibrant replacement of
$\mathrm{RHom}_{{\mathbbm{k}} Q^e}({\mathbbm{k}} Q,{\mathbbm{k}} Q\otimes {\mathbbm{k}} Q)$ in $\mathrm{Mod}({\mathbbm{k}} Q^e)$ and it
determines the dg algebra $\mathbf{\Pi}_d({\mathbbm{k}} Q,c)$ (see part (2) of
Lemma~\ref{sec:smash-prod-deform-1}). Note that $\Theta_{{\mathbbm{k}} Q}$ is
canonically isomorphic to the cone of the following morphism
\begin{equation}
  \label{eq:65}
  {{\mathbbm{k}}{Q_0}} \otimes_{{{\mathbbm{k}}{Q_0}}^e} {\mathbbm{k}} Q^e \to \mathrm{Hom}_{{{\mathbbm{k}}{Q_0}}^e}({{\mathbbm{k}} Q_1},{{\mathbbm{k}}{Q_0}}^e) \otimes_{{{\mathbbm{k}}{Q_0}}^e} {\mathbbm{k}} Q^e\,.
\end{equation}
Since $\mathbf{\Pi}_d({\mathbbm{k}} Q,c)=T_{{\mathbbm{k}} Q}(\Sigma^{d-1}\Theta_{{\mathbbm{k}} Q})$ as graded
algebras, \eqref{eq:62} is an isomorphism of graded algebras, not just
a quasi-isomorphism, and it is induced by the natural inclusions
\begin{equation}
  \label{eq:66}
  {\mathbbm{k}} Q \hookrightarrow \mathcal A(Q,W),\
  \Sigma^{d-2} \mathrm{Hom}_{{{\mathbbm{k}}{Q_0}}^e}({{\mathbbm{k}} Q_1},{{\mathbbm{k}}{Q_0}}^e) \hookrightarrow \mathcal A(Q,W)\
  \text{and}\
  \left\{
  \begin{array}{ccc}
    \Sigma^{d-1}{{\mathbbm{k}}{Q_0}} & \hookrightarrow & \mathcal A(Q,W) \\
    e_i & \mapsto & c_i\,.
  \end{array}\right.
\end{equation}

Note that each one of the ${\mathbbm{k}} Q^e$-modules appearing in \eqref{eq:64}
and \eqref{eq:65} is endowed with an action of $G$. These actions are
compatible so that the corresponding objects are dg
$\Delta$-modules. In this sense, \eqref{eq:64} and \eqref{eq:65} are
morphisms of dg $\Delta$-modules and their cones lie in $\mathrm{Mod}(\Delta)$.

On one hand, $\Theta_{{\mathbbm{k}} Q}$ is a dg $\Delta$-module on its own by
means of Lemma~\ref{sec:inter-betw-dual} and this structure is used to
endow $\mathbf{\Pi}_d({\mathbbm{k}} Q,c)$ with action of $G$ (see
Lemma~\ref{sec:smash-prod-deform-1}). By construction, the structure
of dg $\Delta$-module of $\Sigma\Theta_{{\mathbbm{k}} Q}$ coincides with the one
of the cone of \eqref{eq:64}.

On the other hand, the actions of $G$ on ${\mathbbm{k}} Q$ and
$\mathrm{Hom}_{{{\mathbbm{k}}{Q_0}}^e}({{\mathbbm{k}} Q_1},{{\mathbbm{k}}{Q_0}}^e)$ endow ${\mathbbm{k}} \widetilde Q$ with an action by
automorphisms of graded algebra (recall that $\tilde Q$ is the graded
quiver whose path algebra equals the underlying graded algebra of
$\mathcal A(Q,W)$, see~\eqref{eq:48}). And it follows from
Lemma~\ref{sec:equiv-non-comm-1} that this is actually an action by dg
automorphisms. The mappings of \eqref{eq:66} are $G$-equivariant, and
hence so is \eqref{eq:62}.

\begin{lem}
  \label{sec:group-acti-texorpdfs}
  Assume setting~\ref{sec:appl-ginzb-dg-bis}. Assume that $W$ is
  $G$-invariant up to cyclic permutation.
  \begin{itemize}
  \item Take $\Theta_{{\mathbbm{k}} Q}$ to be the cofibrant replacement of
    $\mathrm{RHom}_{{\mathbbm{k}} Q^e}({\mathbbm{k}} Q,{\mathbbm{k}} Q^e)$ whose suspension is equal to the cone of
    \eqref{eq:64} and endow $\mathbf{\Pi}_d({\mathbbm{k}} Q,c)$ with its action of $G$ by
    dg automorphisms  as in Lemma~\ref{sec:smash-prod-deform-1}.
  \item Endow $\mathrm{Hom}_{{\mathbbm{k}} Q_0^e}({\mathbbm{k}} Q_1,{\mathbbm{k}}{Q_0}^e)$ with its natural
    action of $G$.
  \end{itemize}
  Then,
  \begin{enumerate}
  \item   the latter action together with the one on ${\mathbbm{k}} Q$ endow $\mathcal
    A(Q,W)$ with an action of $G$ by graded algebra automorphisms. For
    this action, \eqref{eq:62} is $G$-equivariant.
  \item the action of $G$ on $\mathcal A(Q,W)$ is an action by dg
    automorphisms and, taking $\mathbf{\Pi}_d({\mathbbm{k}} Q,c)$ like in part (2) of
    Lemma~\ref{sec:smash-prod-deform-1}, the mapping \eqref{eq:62} is
    a $G$-equivariant isomorphism of dg algebras.
  \end{enumerate}
\end{lem}
\begin{proof}
  This follows from the discussion made from the beginning of the
  subsection.
\end{proof}

\subsection{Skew group algebras of Ginzburg dg algebras}
\label{sec:smash-prod-ginzb-1}

Assume setting~\ref{sec:appl-ginzb-dg-bis}. Assume that $W$ is $G$-invariant
up to cyclic permutation. Following
subsection~\ref{sec:group-acti-mathc}, endow $\mathcal A(Q,W)$ with the
resulting action of $G$ by dg automorphisms. The purpose of this
subsection is to prove Corollary~\ref{sec:main-results-4}.

Following \cite{RR85} and \cite{MR2578593}, there exists a graded
${\mathbbm{k}}$-quiver $Q'$ and a (non unital) injective homomorphism
of graded algebras ${\mathbbm{k}} Q'\to {\mathbbm{k}} Q*G$ whose
restriction-of-scalars functor is an equivalence
$\mathrm{Mod}({\mathbbm{k}} Q*G) \to \mathrm{Mod}({\mathbbm{k}}
Q')$. The quiver $Q'$ is constructed in \cite{RR85} when $G$ is cyclic
and in \cite{MR2578593} in full generality. Here is a reminder of its
definition.
\begin{definition}[\cite{MR2578593}]
  \label{sec:skew-group-algebras-1}
Let $[G\backslash Q_0]$ be a complete set of representatives of the
$G$-orbits of vertices of $Q$. For each $i\in Q_0$, denote by $G_i$
the stabiliser of $i$ and let $[G/G_i]$ be a complete set of
representatives of the cosets of $G$ modulo $G_i$. Finally, for all
$i\in Q_0$, let $\mathrm{irr}(G_i)$ be a complete set of representatives of
the isomorphism classes of the irreducible representations of $G_i$;
it is convenient to assume that $\rho = {\mathbbm{k}} G_i e_\rho$ for some
primitive idempotent $e_\rho$ of ${\mathbbm{k}} G_i$, for all
$\rho\in \mathrm{irr}(G_i)$.
\begin{enumerate}
\item Let $\varepsilon$ be the following idempotent of ${\mathbbm{k}} Q*G$,
\begin{equation}
  \label{eq:84}
\varepsilon = \sum_{
  \begin{subarray}{l}
    i\in [G\backslash Q_0]\\
    \rho \in \mathrm{irr}(G_i)
  \end{subarray}} e_i * e_\rho\,,
\end{equation}
\item Let $Q'$ be any quiver as follows.
\begin{itemize}
\item Its vertices are the pairs $(i,\rho)$ where $i\in [G\backslash
  Q_0]$ and $\rho \in \mathrm{irr}(G_i)$.
\item For all vertices $(i,\rho)$ and $(j,\tau)$, denote by
  $M(i,j;\tau)$ the following vector subspace of ${\mathbbm{k}} Q*G$
  \begin{equation}
    \label{eq:83}
  M(i,j;\tau) = \bigoplus_{y\in [G/G_j]} (e_i{\mathbbm{k}} Q_1 e_{y \cdot j}) * (y {\mathbbm{k}} G_j e_\tau)\,;
\end{equation}
by means of the multiplication, ${\mathbbm{k}} Q*G$ is a representation of $G$
and, by restriction, it is a representation of $G_i$; then
$M(i,j;\tau)$ is a subrepresentation of this $G_i$-module; the vector
space with basis being the family of arrows of $Q'$ from $(i,\rho)$
to $(j,\tau)$ is $\mathrm{Hom}_{G_i}(\rho, M(i,j;\tau))$.
\end{itemize}
\end{enumerate}
\end{definition}
There exists an isomorphism of algebras
\begin{equation}
  \label{eq:54}
  {\mathbbm{k}} Q' \xrightarrow{\sim} \varepsilon\cdot ({\mathbbm{k}} Q*G)\cdot \varepsilon
\end{equation}
given as follows,
\begin{itemize}
\item for all vertices $(i,\rho)$ of $Q'$, the corresponding
  idempotent of ${\mathbbm{k}} Q'$ is mapped onto $e_i*e_\rho$,
\item for all vertices $(i,\rho)$ and $(j,\tau)$ of $Q'$, and for all
  $f\in \mathrm{Hom}_{G_i}(\rho,M(i,j;\tau))$, then $f$ is mapped onto
  $f(e_\rho)$. Hence, the arrow $f\colon (i,\rho)\to (j,\tau)$ of $Q'$
  may be identified with the corresponding element $f(e_\rho)$ of
  $(e_i*e_\rho)\cdot ({\mathbbm{k}} Q*G)\cdot (e_j *e_\tau)$.
\end{itemize}


The following lemma introduces a potential $W'$ on $Q'$. The purpose
of this subsection is to prove that $\mathcal A(Q,W)*G$ and $\mathcal
A(Q',W')$ have equivalent derived categories.
\begin{lem}
  \label{sec:smash-prod-ginzb-2}
  Assume setting~\ref{sec:appl-ginzb-dg-bis}. Assume that $W$ is
  $G$-invariant up to cyclic permutation.  There exists a homogeneous
  potential $W'$ of degree $d-3$ on $Q'$ such that the image of $W$
  under the mapping $\mathrm{HH}_{d-3}({\mathbbm{k}} Q) \to \mathrm{HH}_{d-3}({\mathbbm{k}} Q*G)$ induced by
  the natural embedding ${\mathbbm{k}} Q\to {\mathbbm{k}} Q*G$ is equal to the image of $W'$
  under the isomorphism $\mathrm{HH}_{d-3}({\mathbbm{k}} Q') \to \mathrm{HH}_{d-3}({\mathbbm{k}} Q*G)$
  induced by ${\mathbbm{k}} Q' \xrightarrow{~\eqref{eq:54}} {\mathbbm{k}} Q*G\hookrightarrow
    {\mathbbm{k}} Q*G$.
\end{lem}
\begin{proof}
  This follows from standard facts on Hochschild homology, see
  \cite[Lemma 1.1.7, p.~10]{L98}. Here are some details since they are
  needed for the examples. Recall that, for a graded algebra $R$, the
  piece of notation
  $\mathrm{DR}(R)$ stands for the quotient $R$ by the vector subspace generated by the
  commutators $uv-(-1)^{\mathrm{deg}(u) \cdot \mathrm{deg}(v)}vu$, for all homogeneous
  $u,v\in R$.

Following \cite{MR2578593}, there exists a complete family
$(\varepsilon_j)_{j\in J}$ of primitive pairwise orthogonal idempotents of
${\mathbbm{k}} Q*G$ and there exists a subset $I\subseteq J$ such that,
\begin{itemize}
\item the image of ${\mathbbm{k}} Q'\to {\mathbbm{k}} Q*G$ is equal to
  $\varepsilon({\mathbbm{k}} Q*G)\varepsilon$, where $\varepsilon$ denotes $\sum_{i\in I}\varepsilon_i$ and
\item for all $j\in J$ there exists a unique $\alpha(j)\in I$ such
  that $\varepsilon_j\cdot ({\mathbbm{k}} Q*G)\simeq \varepsilon_{\alpha(j)}\cdot ({\mathbbm{k}} Q*G)$ as graded ${\mathbbm{k}} Q*G$-modules.
\end{itemize}
For all $j\in J$, there hence exist homogeneous $a_j,b_j\in {\mathbbm{k}} Q*G$
such that $\varepsilon_j = a_jb_j$ and $\varepsilon_{\alpha(j)} = b_ja_j$.

Now, since $1_{{\mathbbm{k}} Q*G} =\sum_{j\in J}\varepsilon_j$, the following
equalities hold true in $\mathrm{DR}({\mathbbm{k}} Q*G)$, where elements of ${\mathbbm{k}} Q$ are
identified with their natural images in ${\mathbbm{k}} Q*G$,
  \[
  \begin{array}{rcl}
    W & = & \sum_{j\in J}\varepsilon_j W\varepsilon_j \\
        & = & \sum_{j\in J} a_jb_jWa_jb_ja_jb_j \\
    & = & \sum_{i\in I}\varepsilon_i\left(\sum\limits_{j\,s.t.\,\alpha(j)=i} \pm
          b_jWa_j\right) \varepsilon_i\,;
  \end{array}
  \]
  the sign is $(-1)^{\mathrm{deg}(b_j)\cdot (\mathrm{deg}(a_j)+\mathrm{deg}(W))}$. The
  last term above defines an element of
  $\mathrm{DR}(\varepsilon({\mathbbm{k}} Q*G)\varepsilon)$ whose image
  under the mapping
  $\mathrm{DR}(\varepsilon({\mathbbm{k}} Q*G)\varepsilon)\to
  \mathrm{DR}({\mathbbm{k}} Q*G)$ induced by the inclusion mapping
  $\varepsilon\cdot ({\mathbbm{k}} Q*G)\cdot\varepsilon\to
  {\mathbbm{k}} Q*G$ is equal to the image of $W$ under the mapping
  $\mathrm{DR}({\mathbbm{k}} Q) \to \mathrm{DR}({\mathbbm{k}} Q*G)$
  induced by ${\mathbbm{k}} Q \to {\mathbbm{k}} Q*G$. Since
  ${\mathbbm{k}} Q' \simeq \varepsilon({\mathbbm{k}} Q*G)\varepsilon$
  as graded algebras, this proves the lemma.
\end{proof}

The potential provided by Lemma~\ref{sec:smash-prod-ginzb-2} fits the
requirements of Corollary~\ref{sec:main-results-4} which can now be
proved. This corollary is restated below for convenience.
\begin{cor}
  \label{sec:skew-group-algebras-3}
  Assume setting~\ref{sec:appl-ginzb-dg-bis}. Assume that $W$ is
  $G$-invariant up to cyclic permutation.
  \begin{enumerate}
  \item The action of $G$ on ${\mathbbm{k}} Q$ extends to an action of $G$ on
    $\mathcal A(Q,W)$ by dg automorphisms and $\mathcal A(Q,W)*G$ is
    $d$-Calabi-Yau.
  \item For all graded ${\mathbbm{k}}$-quivers $Q'$ and for all (non unital)
    graded algebra homomorphisms ${\mathbbm{k}} Q'\to {\mathbbm{k}} Q*G$ whose
    restriction-of-scalars functor is an equivalence from $\mathrm{Mod}({\mathbbm{k}}
      Q*G)$ to $\mathrm{Mod}({\mathbbm{k}} Q')$ (see \cite{MR2578593}), there exists a
    homogeneous of degree $d-3$ potential $W'$ on $Q'$ such that ${\mathbbm{k}}
    Q'\to {\mathbbm{k}} Q*G$ extends to a (non unital) dg algebra homomorphism
    \[
    \mathcal A(Q',W')\to \mathcal A(Q,W)*G
    \]
    whose restriction-of-scalars functor induces an equivalence
    \[
    \mathcal{D}(\mathcal A(Q,W)*G) \xrightarrow{\sim} \mathcal{D}(\mathcal A(Q',W'))\,.
    \]
  \end{enumerate}
\end{cor}
\begin{proof}
  Let $W'$ be as in Lemma~\ref{sec:smash-prod-ginzb-2}.  The following
  diagram is commutative. In this diagram, $B$ stands for the Connes'
  boundary and the horizontal arrows are induced by ${\mathbbm{k}} Q\to {\mathbbm{k}} Q*G$
  and ${\mathbbm{k}} Q'\to {\mathbbm{k}} Q*G$, and the rightmost horizontal arrows are
  isomorphisms because ${\mathbbm{k}} Q'$ and ${\mathbbm{k}} Q*G$ are Morita equivalent.
  \begin{equation}
    \label{eq:67}
    \xymatrix{
      \mathrm{HH}_{d-3}({\mathbbm{k}} Q)^G \ar@{^(->}[r] \ar[d] & \mathrm{HH}_{d-3}({\mathbbm{k}} Q)
      \ar[r] \ar[d]_B & \mathrm{HH}_{d-3}({\mathbbm{k}} Q* G)  \ar[d]_B &
      \mathrm{HH}_{d-3}({\mathbbm{k}} Q') \ar[l]_{\sim} \ar[d]^B\\
      \mathrm{HH}_{d-2}({\mathbbm{k}} Q)^G \ar@{^(->}[r] & \mathrm{HH}_{d-2}({\mathbbm{k}} Q) \ar[r] &
      \mathrm{HH}_{d-2}({\mathbbm{k}} Q* G) & \mathrm{HH}_{d-2}({\mathbbm{k}} Q') \ar[l]^{\sim}
    }
  \end{equation}
  This defines homology classes $c,\overline c,c'$ as follows,
  \begin{enumerate}[(a)]
  \item $c\in \mathrm{HH}_{d-2}({\mathbbm{k}} Q)^G$ is the image of $W$ under 
    $B\colon \mathrm{HH}_{d-3}({\mathbbm{k}} Q) \to \mathrm{HH}_{d-2}({\mathbbm{k}} Q)$,
  \item $\overline c\in \mathrm{HH}_{d-2}({\mathbbm{k}} Q*G)$ is the image of $c$ under
     $\mathrm{HH}_{d-2}({\mathbbm{k}} Q) \to \mathrm{HH}_{d-2}({\mathbbm{k}} Q*G)$ and
   \item $c' \in \mathrm{HH}_{d-2}({\mathbbm{k}} Q')$ is the pre-image of $\overline c$
     under $\mathrm{HH}_{d-2}({\mathbbm{k}} Q') \xrightarrow{\sim}\mathrm{HH}_{d-2}({\mathbbm{k}} Q*G)$.
  \end{enumerate}
  In particular,
  \begin{enumerate}[(a)]
    \setcounter{enumi}{3}
  \item $c'$ is the image of $W'$ under
    $B\colon \mathrm{HH}_{d-3}({\mathbbm{k}} Q') \to \mathrm{HH}_{d-2}({\mathbbm{k}} Q')$.
  \end{enumerate}
   The following
  realisations of the graded algebras $\mathbf{\Pi}_d({\mathbbm{k}} Q,c)$ and $\mathbf{\Pi}_d({\mathbbm{k}}
  Q*G,\overline c)$ are used in this proof,
  \begin{itemize}
  \item as for $\mathbf{\Pi}_d({\mathbbm{k}} Q,c)$, take the de-suspension of the cone of
    \eqref{eq:64} as a cofibrant replacement  $\Theta_{{\mathbbm{k}} Q}$ of
    $\mathrm{RHom}_{{\mathbbm{k}} Q^e}({\mathbbm{k}} Q,{\mathbbm{k}} Q^e)$ and
  \item as for $\mathbf{\Pi}_d({\mathbbm{k}} Q,\overline c)$, take $\Theta_{{\mathbbm{k}} Q}*G$ as
    cofibrant replacement of $\mathrm{RHom}_{({\mathbbm{k}} Q*G)^e}({\mathbbm{k}} Q*G,({\mathbbm{k}} Q*G)^e)$
    (see the proof of Theorem~\ref{sec:skew-group-algebras-2}).
  \end{itemize}
  
  The proof of the proposition is based on a sequence of
  quasi-isomorphic dg algebras. Actually, because of the choices just
  made, this sequence comes from a sequence of composable
  isomorphisms of dg algebras.

  As proved in Lemma~\ref{sec:group-acti-texorpdfs}, there is
  $G$-equivariant isomorphism
  $\mathbf{\Pi}_d({\mathbbm{k}} Q,c) \xrightarrow{\sim} \mathcal A(Q,W)$ of dg algebras.
  Accordingly, there is an isomorphism of dg algebras
  \begin{equation}
    \label{eq:50}
    \mathbf{\Pi}_d({\mathbbm{k}} Q,c)*G \xrightarrow{\sim} \mathcal A(Q,W)*G\,.    
  \end{equation}
  Following Theorem~\ref{sec:main-results-2}, there is an
  isomorphism of dg algebras
  \begin{equation}
    \label{eq:51}
    \mathbf{\Pi}_d({\mathbbm{k}} Q,c) * G \xrightarrow{\sim} \mathbf{\Pi}_d({\mathbbm{k}}
    Q*G,\overline c)\,.
  \end{equation}

  Now, by \cite[Theorem 5.8]{MR2795754}, the (non unital) algebra
  homomorphism ${\mathbbm{k}} Q'\to {\mathbbm{k}} Q*G$ extends to a
  (non unital) dg algebra homomorphism
  $\mathbf{\Pi}_d({\mathbbm{k}} Q',c') \to
  \mathbf{\Pi}_d({\mathbbm{k}} Q*G,\overline c)$ whose
  restriction-of-scalars functor is a triangle equivalence
  \begin{equation}
    \label{eq:53}
    \mathcal{D}(\mathbf{\Pi}_d({\mathbbm{k}} Q',c'))
    \simeq \mathcal{D}(\mathbf{\Pi}_d({\mathbbm{k}} Q*G,\overline c))\,.
  \end{equation}
  Finally, in view of point (d) above relating $c'$ to $W'$, it
  follows from \cite[Theorem 6.3]{MR2795754} that the identity mapping
  ${\mathbbm{k}} Q'\to {\mathbbm{k}} Q'$ extends to a quasi-isomorphism of dg algebras
  \begin{equation}
    \label{eq:68}
    \mathbf{\Pi}_d({\mathbbm{k}} Q',c') \xrightarrow{\mathrm{qis}} \mathcal A(Q',W')\,.
  \end{equation}
  The assertion (2) follows from \eqref{eq:50}, \eqref{eq:51},
  \eqref{eq:53}, and \eqref{eq:68}. The assertion (1) follows from (2)
  and from the fact that $\mathcal A(Q',W')$ is $d$-Calabi-Yau.
\end{proof}

\subsection{Examples}
\label{sec:examples}

Assume setting~\ref{sec:appl-ginzb-dg-bis}. Assume that $W$ is
$G$-invariant up to cyclic permutation. This subsection illustrates
the computation of a quiver with potential $(Q',W')$ such as in
Corollary~\ref{sec:skew-group-algebras-3}.  In the examples of this
subsection, graded ${\mathbbm{k}}$-quivers are concentrated in degree $0$ and the
resulting Ginzburg dg algebras are $3$-Calabi-Yau. For convenience, for all
$g\in G$ and $i\in Q_0$, the vertex of $Q$ associated with the
idempotent $^ge_i$ is denoted by $g\cdot i$. In some places it is also
convenient to identify the elements of ${\mathbbm{k}} Q$ and $G$ with their
respective images in ${\mathbbm{k}} Q*G$.

The computation of a potential $W'$ such as above is therefore made in
two steps.
\begin{enumerate}
\item Express $W'$ as an element of $\varepsilon\cdot({\mathbbm{k}} Q*G)\cdot \varepsilon$, see the
  proof of Lemma~\ref{sec:smash-prod-ginzb-2}.
\item Express the result of the first step as a linear combination of
  paths in $Q'$.
\end{enumerate}

\begin{rem}
  \label{sec:examples-3}
For the second step, it is worth noting that any oriented cycle in $Q$
with source $i_0\in [G\backslash Q_0]$ is equal in ${\mathbbm{k}} Q*G$ to a
linear combination of products of the shape
\begin{equation}
  \label{eq:75}
  (\alpha_1*\sigma_1)\cdots(\alpha_j*\sigma_j)\cdots(\alpha_t*\sigma_t)\cdot(1*g)\,,
\end{equation}
where there exists a sequence $i_1,\ldots,i_t\in [G\backslash Q_0]$
such that $\alpha_j$ is an arrow of $Q$ from $i_{j-1}$ to
$\sigma_j\cdot i_j$ and $\sigma_j\in [G/G_j]$,
for all $j$, and where $g\in G_{i_0}$; besides, when a given index
$j\in \{1,\ldots,t-1\}$ is such that $G_{i_j}$ is abelian, then
$\sum_{\rho\in \mathrm{irr}(G_{i_j})} e_\rho$ is the unity of ${\mathbbm{k}} G_i$, and
hence
\begin{equation}
  \label{eq:76}
  (\alpha_j *\sigma_j) \cdot (\alpha_{j+1}*\sigma_{j+1}) =
  \sum_{\rho\in \mathrm{irr}(G_{i_j})} (\alpha_j *\sigma_j) \cdot
  (e_{i_j}*e_\rho)\cdot (\alpha_{j+1}*\sigma_{j+1})\,.
\end{equation}
\end{rem}

\begin{ex}
  \label{sec:examples-1}
  Assume that $(Q,W)$ is the following quiver with potential
  
  \begin{equation}
    \label{eq:74}
    \xymatrix{
    & 1 \ar@<2pt>[rd] \ar@<2pt>[ld]
    \ar@(ul,ur)[]
    \\
    2 \ar@<2pt>[rr] \ar@<2pt>[ru]
    & & 3
    \ar@<2pt>[lu] \ar@<2pt>[ll]
    &
    W = x_{12}x_{23}x_{32}x_{21} - x_{13}x_{32}x_{23}x_{31} +
    x_{13}x_{31}x_{11} - x_{12}x_{21}x_{11}
  }
  \end{equation}
  where $x_{ij}$ denotes the unique arrow $i\to j$ of $Q$ if any, for
  all vertices $i,j$. Recall that $\mathcal A(Q,W)$ is
  quasi-isomorphic to $H^0\mathcal A(Q,W)$ and that the latter is a
  non commutative crepant resolution of the suspended pinch-point, see
  \cite[Section 8]{MR2592501}. Let $G$ be the subgroup
  $\langle (23) \rangle$ of order $2$ of the symmetric group
  $\mathscr{S}_3$. Then, $G$ acts on ${\mathbbm{k}} Q$ as follows
  \begin{itemize}
  \item $\sigma \cdot i = \sigma(i)$ for all $\sigma\in G$ and $i\in
    \{1,2,3\}$,
  \item $^{(23)}x_{12} = -x_{13}$ and $^{(23)}x_{13}= - x_{12}$, and
  \item $^\sigma x_{ij}= x_{\sigma(i)\sigma(j)}$ for all $\sigma\in G$
    and for all arrows $x_{ij}$ of $Q$ distinct from $x_{12}$ and
    $x_{13}$.
  \end{itemize}
  Note that, for this action, $W$ is invariant up to cyclic
  permutation. One may take the following data to define $Q'$.
\begin{itemize}
\item $[G\backslash Q_0] = \{1,2\}$, $[G/G_1] = \{\mathrm{Id}\}$ and $[G/G_2]
  = G$,
\item $\mathrm{irr}(G_1)$ consists of the trivial representation
  $\rho_+={\mathbbm{k}} G_1\cdot(\mathrm{Id}+(23))$ of $G_1$ and the non trivial one
  $\rho_- = {\mathbbm{k}} G_1\cdot (\mathrm{Id}-(23))$,
\item $\mathrm{irr}(G_2)$ consists of the trivial representation ${\mathbbm{k}}$ of $G_2$.
\end{itemize}
Then (see~\eqref{eq:84}), $\varepsilon = \varepsilon_++\varepsilon_-+\varepsilon_2$, where
$\varepsilon_+ = \frac{1}{2}(e_1*\mathrm{Id} + e_1*(23))$,
$\varepsilon_- = \frac{1}{2}(e_1*\mathrm{Id} - e_1*(23))$ and $\varepsilon_2 = e_2*e$, recall
that $e$ denotes the neutral element of $G$. Hence
(see~\eqref{eq:83}),
\[
\begin{array}{rclccrcl}
  M(1,1;\rho_+) & = & \mathrm{span}(x_{11}\varepsilon_+) \simeq \rho_+
  &&  M(2,1;\rho_+) & = & \mathrm{span}(x_{21}\varepsilon_+)\simeq {\mathbbm{k}} \\
  M(1,1;\rho_-) & = &  \mathrm{span}(x_{11}\varepsilon_-)\simeq \rho_-
  && M(2,1;\rho_-) & = & \mathrm{span}(x_{21}\varepsilon_-) \simeq {\mathbbm{k}} \\
  M(1,2;{\mathbbm{k}}) & = & \underset{\simeq \rho_+}{\underbrace{\mathrm{span}(\varepsilon_+x_{12}e_2)}} \bigoplus
                                                               \underset{\simeq
                                                               \rho_-}{\underbrace{\mathrm{span}(\varepsilon_-x_{12}\varepsilon_2)}}
  && M(2,2;{\mathbbm{k}}) & = & \mathrm{span}(x_{23}(23)\varepsilon_2)\simeq{\mathbbm{k}}  \,,\\
\end{array}
\]
note that, because of the equalities $^{(23)}x_{12}=-x_{13}$ and
$^{(23)}x_{13}=-x_{12}$,
\begin{equation}
  \label{eq:78}
  \begin{array}{rclccrcl}
    \varepsilon_-x_{11}\varepsilon_+ & = & 0, &&     \varepsilon_+x_{12}\varepsilon_2 & = &
                                                                                            -\varepsilon_+x_{13}(23)\varepsilon_2,\\
    \varepsilon_+x_{11}\varepsilon_-& = &0, && \varepsilon_- x_{12}\varepsilon_2 & = & \varepsilon_-x_{13}(23)\varepsilon_2\,.
  \end{array}
\end{equation}
Hence $Q'$ is the following quiver
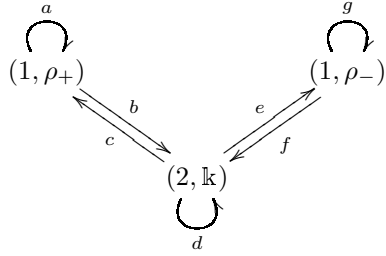
\[
\xymatrix{
  (1,\rho_+) \ar@(lu,ru)[]^a \ar@<2pt>[rd]^b &&
  (1,\rho_-) \ar@(lu,ru)[]^g \ar@<2pt>[ld]^f
  \\
  & (2,{\mathbbm{k}}) \ar@(dl,dr)[]_d \ar@<2pt>[lu]^c
  \ar@<2pt>[ru]^e
}
\]
where
\[
\begin{array}{rclcrclcrcl}
  a & = & \varepsilon_+x_{11}\varepsilon_+  && b & = & \varepsilon_+x_{12}\varepsilon_2 && c & = & \varepsilon_2x_{21}\varepsilon_+\\
  d & = & \varepsilon_2 x_{23}(23)\varepsilon_2 && e & = & \varepsilon_2x_{21}\varepsilon_- && f & = &
                                                                   \varepsilon_-x_{12}\varepsilon_2
  \\
  g & = & \varepsilon_-x_{11}\varepsilon_-\,.
\end{array}
\]
Now, $W=e_1We_1 = \varepsilon_+ W \varepsilon_+ + \varepsilon_-W\varepsilon_-$ because $e_1 = \varepsilon_++\varepsilon_-$.
Note that $\varepsilon_-W\varepsilon_+ = \varepsilon_+W\varepsilon_-=0$ because $W$ is $G$-invariant up to
cyclic permutation. Therefore, by Lemma~\ref{sec:smash-prod-ginzb-2},
it is possible to assume that
\begin{equation}
  \label{eq:77}
  W' = \varepsilon_+W\varepsilon_+ + \varepsilon_-W\varepsilon_-\,.
\end{equation}
There only remains to decompose $W'$ as a linear combination of
oriented cycles in $Q'$. For this purpose, it is worth noting that
$\sum_{\rho\in \mathrm{irr}(G_1)} e_1*e_\rho = \varepsilon_++\varepsilon_-$ and
$\sum_{\rho\in \mathrm{irr}(G_2)} e_2*e_\rho = \varepsilon_2$. Using
Remark~\ref{sec:examples-3} and \eqref{eq:78}, the computations yield,
\[
\begin{array}{rcl}
  \varepsilon_+ x_{12}x_{23}x_{32}x_{21}\varepsilon_+
  & = &
        (\varepsilon_+x_{12}\varepsilon_2)\cdot(\varepsilon_2x_{23}(23)\varepsilon_2)\cdot (\varepsilon_2
        x_{23}(23)\varepsilon_2) \cdot (\varepsilon_2x_{21}\varepsilon_+) \\
  \varepsilon_- x_{12}x_{23}x_{32}x_{21}\varepsilon_-
  & = &
        (\varepsilon_-x_{12}\varepsilon_2)\cdot(\varepsilon_2x_{23}(23)\varepsilon_2)\cdot (\varepsilon_2
        x_{23}(23)\varepsilon_2) \cdot (\varepsilon_2x_{21}\varepsilon_-) \\
  \varepsilon_+ x_{13}x_{32}x_{23}x_{31}\varepsilon_+
  & = &
        -(\varepsilon_+x_{12}\varepsilon_2) \cdot (\varepsilon_2 x_{23}(23)\varepsilon_2) \cdot (\varepsilon_2
        x_{23}(23) \varepsilon_2) \cdot (\varepsilon_2x_{21}\varepsilon_+) \\
  \varepsilon_- x_{13}x_{32}x_{23}x_{31}\varepsilon_-
  & = &
        -(\varepsilon_-x_{12}\varepsilon_2) \cdot (\varepsilon_2 x_{23}(23)\varepsilon_2) \cdot (\varepsilon_2
        x_{23}(23) \varepsilon_2) \cdot (\varepsilon_2x_{21}\varepsilon_-) \\
  \varepsilon_+x_{13}x_{31}x_{11}\varepsilon_+
  & = &
        -(\varepsilon_+x_{12}\varepsilon_2) \cdot (\varepsilon_2x_{21}\varepsilon_+) \cdot
        (\varepsilon_+x_{11}\varepsilon_+) \\
  \varepsilon_-x_{13}x_{31}x_{11}\varepsilon_-
  & = &
        -(\varepsilon_-x_{12}\varepsilon_2) \cdot (\varepsilon_2x_{21}\varepsilon_-) \cdot
        (\varepsilon_-x_{11}\varepsilon_-) \\
  \varepsilon_+x_{12}x_{21}x_{11}\varepsilon_+
  & = &
        (\varepsilon_+x_{12}\varepsilon_2) \cdot (\varepsilon_2x_{21}\varepsilon_+) \cdot (\varepsilon_+
        x_{11}\varepsilon_+) \\
  \varepsilon_-x_{12}x_{21}x_{11}\varepsilon_-
  & = &
        (\varepsilon_-x_{12}\varepsilon_2) \cdot (\varepsilon_2x_{21}\varepsilon_-) \cdot (\varepsilon_-
        x_{11}\varepsilon_-)\,.
\end{array}
\]
It follows from these computations and from \eqref{eq:74} and
\eqref{eq:77} that $W'=2(bd^2c+fd^2e-bca-feg)$.
\end{ex}

\begin{ex}
  \label{sec:examples-2}
  Assume that $(Q,W)$ is the following quiver with potential
  \begin{equation}
    \label{eq:79}
  \xymatrix{
    & 1 \ar@<2pt>[rd] \ar@<2pt>[ld]
    &&
        W = x_{12}x_{23}x_{31}- x_{13}x_{32}x_{21}
    \\
    2 \ar@<2pt>[rr] \ar@<2pt>[ru]
    & & 3
        \ar@<2pt>[lu] \ar@<2pt>[ll]}
  \end{equation}
  where the naming of the arrows follows the same convention as in
  Example~\ref{sec:examples-1}. Let $G$ be the symmetric group
  $\mathscr{S}_3$ on $\{1,2,3\}$.  It acts on ${\mathbbm{k}} Q$ as follows.
  \begin{itemize}
  \item The action on the set of vertices is the natural one.
  \item $^\sigma x_{i,j} = (-1)^\sigma x_{\sigma(i)\sigma(j)}$ for all
    $\sigma\in G$ and $i,j\in \{1,2,3\}$, where $(-1)^\sigma$ is the
    sign of the permutation $\sigma$.
  \end{itemize}
  For this action, $W$ is invariant up to cyclic
  permutation. Actually,
  $W=\frac{1}{3}\sum_{\sigma\in G}\,^\sigma(x_{12}x_{23}x_{31})$.

Take the following sets of representatives to define $Q'$,
\begin{itemize}
\item $[G\backslash Q_0] = \{1\}$, $[G/G_1] = \{\mathrm{Id}, c,c^2\}=
  \langle c\rangle$, where $c=(123)$,
\item $\mathrm{irr}(G_1)$ consists of the trivial representation
  $\rho_+={\mathbbm{k}} G_1\cdot(\mathrm{Id}+(23))$ of $G_1$ ($=G$) and the non trivial
  one $\rho_- = {\mathbbm{k}} G_1\cdot (\mathrm{Id}-(23))$.
\end{itemize}
Then, $\varepsilon$ is the sum of two primitive pairwise orthogonal
idempotents of ${\mathbbm{k}} Q*G$,
\[
\varepsilon = \varepsilon_++\varepsilon_-\,,
\]
where $\varepsilon_+ = \frac{1}{2}(e_1*\mathrm{Id} + e_1*(23))$ and
$\varepsilon_- = \frac{1}{2}(e_1*\mathrm{Id} - e_1*(23))$. Hence,
\begin{itemize}
\item $M(1,1;\rho_+) = \mathrm{span}(x_{12}c\varepsilon_+,x_{13}c^2\varepsilon_+) =
  \underset{\simeq \rho_+}{\underbrace{\mathrm{span}(\varepsilon_+x_{12}c\varepsilon_+)}} \oplus
  \underset{\simeq \rho_-}{\underbrace{\mathrm{span}(\varepsilon_-x_{12}c\varepsilon_+)}}$,
\item $M(1,1;\rho_-) = \mathrm{span}(x_{12}c\varepsilon_-,x_{13}c^2\varepsilon_-) =
    \underset{\simeq \rho_+}{\underbrace{\mathrm{span}(\varepsilon_+x_{12}c\varepsilon_-)}} \oplus
  \underset{\simeq \rho_-}{\underbrace{\mathrm{span}(\varepsilon_-x_{12}c\varepsilon_-)}}$,
\end{itemize}
note that $(23) x_{13}c^2 = -x_{12}c(23)$ and $(23)x_{12}c =
-x_{13}c^2(23)$ so that
\begin{equation}
  \label{eq:80}
  \begin{array}{rclccrcl}
    \varepsilon_+x_{12}c\varepsilon_+ &=& -\varepsilon_+x_{13}c^2\varepsilon_+\,, && \varepsilon_-x_{12}c\varepsilon_+ &=&
                                                                  \varepsilon_-x_{13}c^2\varepsilon_+\,,
    \\
    \varepsilon_+x_{12}c\varepsilon_- &=& -\varepsilon_+x_{13}c^2\varepsilon_-\,, && \varepsilon_-x_{12}c\varepsilon_- &=& -\varepsilon_-x_{13}c^2\varepsilon_-\,.
  \end{array}
\end{equation}
Hence $Q'$ is the following quiver
\[
\xymatrix{
  (1,\rho_+) \ar@(ul,ur)^\alpha \ar@<2pt>[rr]^\beta &&
  (1,\rho_-) \ar@(ul,ur)^\gamma \ar@<2pt>[ll]^\delta
}
\]
where
\[
\begin{array}{rclccrclccrclccrcl}
  \alpha & = & \varepsilon_+x_{12}c\varepsilon_+\,,
  && \beta & = & \varepsilon_+x_{12}c\varepsilon_-\,,
  && \gamma & = & \varepsilon_-x_{12}c\varepsilon_-\,,
  && \delta & = & \varepsilon_- x_{12}c\varepsilon_+\,.
\end{array}
\]
Since $e_1 = \varepsilon_++\varepsilon_-$,
\[
W = e_1We_1 = \varepsilon_+ W \varepsilon_+ + \varepsilon_-W\varepsilon_-\,,
\]
here $\varepsilon_-W\varepsilon_+ = \varepsilon_+W\varepsilon_-=0$ because $W$ is $G$-invariant up to
cyclic permutation. Therefore, by Lemma~\ref{sec:smash-prod-ginzb-2},
it is possible to assume that
\begin{equation}
\label{eq:81}
  W' = \varepsilon_+W\varepsilon_+ + \varepsilon_-W\varepsilon_-\,.
\end{equation}
There only remains to decompose $W'$ as a linear combination of
oriented cycles in $Q'$. For this purpose, note that
$x_{12}x_{23}x_{31} = x_{12}cx_{12}cx_{12}c$, that
$x_{13}x_{32}x_{21} = x_{13}c^2x_{13}c^2x_{13}c^2$ and that
$\sum_{\rho\in \mathrm{irr}(G_1)} e_1*e_\rho = \varepsilon_++\varepsilon_-$. Hence, using
Remark~\ref{sec:examples-3} and \eqref{eq:80}, the computations yield,
\[
\begin{array}{rcl}
  \varepsilon_+ x_{12}x_{23}x_{31}\varepsilon_+ & = &
\alpha\alpha\alpha+\alpha\beta\delta+\beta\gamma\delta+\beta\delta\alpha\,,
\\ \varepsilon_- x_{12}x_{23}x_{31}\varepsilon_- & = &
\gamma\gamma\gamma+\gamma\delta\beta+\delta\beta\gamma+\delta\alpha\beta\,,
\\ \varepsilon_+ x_{13}x_{32}x_{21}\varepsilon_+ & = &
-(\alpha\alpha\alpha+\alpha\beta\delta+\beta\gamma\delta+\beta\delta\alpha)\,,
\\ \varepsilon_- x_{13}x_{32}x_{21}\varepsilon_- & = &
-(\gamma\gamma\gamma+\gamma\delta\beta+\delta\beta\gamma+\delta\alpha\beta)\,. \\
\end{array}
\]
It follows from these computations and from \eqref{eq:79} and
\eqref{eq:81} that
\[
W'=2\sum_{\text{$w$ cycle of length $3$}}w\,.
\]
\end{ex}

\section{Application to cluster tilting objects}
\label{sec:appl-gener-clust}

Assume setting~\ref{sec:main-results-5}. It is restated below for
convenience.
\begin{setting}
  \label{sec:main-results-5-bis}
  Let ${\mathbbm{k}}$ be a field. Let $G$ be a finite group such that
  $\mathrm{char}({\mathbbm{k}})$ does not divide
  $\mathrm{Card}(G)$. Let $A$ be dg ${\mathbbm{k}}$-algebra acted upon
  by $G$ by dg automorphisms. Assume that $A$ is concentrated in non
  positive degrees and $H^0(A)$ is finite dimensional. Assume that
  both $A$ and $A*G$ are $3$-Calabi-Yau.
\end{setting}
Note that $G$ acts on
$H^0(A)$ by algebra automorphisms and that $H^0(A*G) \simeq H^0(A)*G$
as algebras. In particular, $H^0(A*G)$ is finite dimensional. These
conditions ensure the existence of the associated generalised cluster
categories ${\mathcal{C}_A}$ and ${\mathcal{C}_{A*G}}$. The purpose of this section is to relate
these two categories by means of biadjoint functors
\[
\xymatrix{
  F_\lambda \colon {\mathcal{C}_A}
  \ar@<+0.5ex>[rr] &&
  {\mathcal{C}_{A*G}} \colon F_\cdot
  \ar@<+0.5ex>[ll]
}
\]
in order to compare the cluster tilting objects of ${\mathcal{C}_A}$ and those of
${\mathcal{C}_{A*G}}$.

Such a relationship appeared first in  \cite{lemeur}:
given a Galois covering $R'\to R$ with group $G$, which was not assumed
to be finite, where $R$ is a finite dimensional piecewise hereditary
algebra of quiver type, \cite{lemeur} shows that $R'$ is piecewise
hereditary, it introduces an action of $G$ on the cluster category
$C_{R'}$ of $R'$ and a \emph{push-down} functor from $C_{R'}$ to the
cluster category $C_R$ of $R$. This provides a description of the
cluster tilting objects of $C_R$ in terms of $G$-invariant cluster
tilting objects of $C_{R'}$. In particular, it is proved in
\cite{lemeur} that every indecomposable direct summand of any cluster
tilting object of $C_R$ is isomorphic to the image under the push-down
functor of an indecomposable direct summand of a $G$-invariant cluster
tilting object of $C_{R'}$. Recall that in that setting, $R$ is Morita
equivalent to $R'*G$ \cite[Proposition 2.4]{MR2170541}.

Following \cite{A08} the generalised cluster categories of $A$ and
$A*G$ are, respectively,
\[
\begin{array}{l}
  {\mathcal{C}_A} = \mathrm{per}(A)/\mathcal{D}_{\mathrm{fd}}(A)\,,\\
  {\mathcal{C}_{A*G}} = \mathrm{per}(A*G)/\mathcal{D}_{\mathrm{fd}}(A*G)\,.
\end{array}
\]
They are $\mathrm{Hom}$-finite, Krull-Schmidt and $2$-Calabi-Yau, that is, the square
of their suspensions are Serre functors. Denote by $\pi_A$ and
$\pi_{A*G}$ the canonical functors
\[
\begin{array}{l}
  \pi_A \colon \mathrm{per}(A) \to {\mathcal{C}_A}\,,\\
  \pi_{A*G} \colon \mathrm{per}(A*G) \to {\mathcal{C}_{A*G}}\,.
\end{array}
\]
For simplicity, given an object $T$ in ${\mathcal{C}_A}$ or ${\mathcal{C}_{A*G}}$, its
endomorphism algebra is denoted by $\mathrm{End}(T)$. The objects
$\overline A := \pi_A(A)$ and $\overline{A*G} := \pi_{A*G}(A*G)$ are
cluster tilting in ${\mathcal{C}_A}$ and ${\mathcal{C}_{A*G}}$, respectively.

The section is organised as follows. Subsection~\ref{sec:action-g-ca}
introduces an action of $G$ on ${\mathcal{C}_A}$ and defines the adjoint functors
$F_\lambda \colon {\mathcal{C}_A}\to {\mathcal{C}_{A*G}}$ and $F_\cdot\colon {\mathcal{C}_{A*G}} \to {\mathcal{C}_A}$.
Subsection~\ref{sec:cy-categ-relat} proves
Theorem~\ref{sec:main-results-6} which establishes a first general
comparison between the cluster tilting subcategories of ${\mathcal{C}_A}$ and
those of ${\mathcal{C}_{A*G}}$. This comparison is strengthened in
subsection~\ref{sec:clust-tilt-objects} when $G$ acts freely on the
isomorphism classes of indecomposable projective
$H^0(A)$-modules. Finally, these results are applied in
subsection~\ref{sec:appl-finite-type} to prove
Corollary~\ref{sec:main-results-7} which relates ${\mathcal{C}_A}$ being acyclic
to ${\mathcal{C}_{A*G}}$ being acyclic.


\subsection{The action of \texorpdfstring{$G$}{G} on
  \texorpdfstring{${\mathcal{C}_A}$}{C(A)} and the adjunction between
  \texorpdfstring{${\mathcal{C}_A}$}{C(A)} and \texorpdfstring{${\mathcal{C}_{A*G}}$}{C(A*G)}}
\label{sec:action-g-ca}

The action of $G$ on $\mathcal{D}(A)$ stabilises $\mathrm{per}(A)$ and
$\mathcal{D}_{\mathrm{fd}}(A)$.  Accordingly, it induces a strict action of $G$ on ${\mathcal{C}_A}$ by
strict automorphisms of triangulated categories. For all $M\in {\mathcal{C}_A}$,
define the \emph{stabiliser} of $M$ as the subgroup
$\{g\in G\ |\ M^g\simeq M\}$ of $G$, it is denoted by $G_M$.

Consider the adjoint pair of triangle functors of
subsection~\ref{sec:adjunctions},
\[
\xymatrix{ - \otimes^{\mathrm L}_A (A*G) \colon \mathcal{D}(A) \ar@<+0.5ex>[rr]
  && \mathcal{D}(A*G) \colon \mathrm{Rres}  \ar@<+0.5ex>[ll] }
\]
\begin{itemize}
\item $-\otimes^{\mathrm L}_A(A*G)$ maps $\mathrm{per}(A)$ into $\mathrm{per}(A*G)$, and
  $\mathrm{Rres}$ maps $\mathrm{per}(A*G)$ into $\mathrm{per}(A)$ because $A*G \simeq
  \oplus_{g\in G} A$ in $\mathrm{Mod}(A)$, and
\item $-\otimes^{\mathrm L}_A(A*G)$ maps $\mathcal{D}_{\mathrm{fd}}(A)$ into $\mathcal{D}_{\mathrm{fd}}(A*G)$
  because $A*G \simeq \oplus_{g\in G} A$ in $\mathrm{Mod}(A^{\mathrm{op}})$, and
  $\mathrm{Rres}$ maps $\mathcal{D}_{\mathrm{fd}}(A*G)$ into $\mathcal{D}_{\mathrm{fd}}(A)$.
\end{itemize}
\begin{definition}
  \label{sec:acti-texorpdfstr-tex}
  The functors $F_\lambda \colon {\mathcal{C}_A} \to {\mathcal{C}_{A*G}}$ and
  $F_\cdot \colon {\mathcal{C}_{A*G}} \to {\mathcal{C}_A}$ are the triangle functors induced by
  $-\otimes^{\mathrm L}_A(A*G) \colon \mathcal{D}(A)\to \mathcal{D}(A*G)$ and
  $\mathrm{Rres} \colon \mathcal{D}(A*G) \to \mathcal{D}(A)$, respectively.
\end{definition}

The functor $F_\lambda$ between (generalised) cluster categories has
already been considered in other situations, for instance,
\begin{itemize}
\item in \cite[Section 4]{lemeur} between cluster categories
  $C_{R'}\to C_R$ where $R$ is a piecewise hereditary algebra and $R'$
  is a Galois covering of $R$ with group $G$,
\item in \cite[Section 3]{MR3899030}, where $A$ is the Ginzburg dg algebra of
  any quiver with potential $(Q,P)$ and the action of $G$ is given by a
  free action on $Q$ by quiver automorphisms,
\item in \cite[Corollary 3.6]{AP}, where $A$ is the Ginzburg dg
  algebra of any quiver with potential $(Q,P)$, and the action of
  $G=\mathbb Z/2\mathbb Z$ is given by an action on $Q$ by quiver
  automorphisms, the last condition is actually superfluous in that
  result once one knows that the skew group algebra $\mathcal A(Q,W)*\mathbb Z/2\mathbb Z$ is
  Morita equivalent to a Ginzburg dg algebra.
  \end{itemize}
\begin{prop}
  \label{sec:adjunct-betw-ca-2}
  Assume setting~\ref{sec:main-results-5-bis}.
  \begin{enumerate}
  \item The pairs $(F_\lambda,F_\cdot)$ and $(F_\cdot,F_\lambda)$ are
    adjoint.
  \item The units of these adjunctions split (functorially).
  \item There exists a family
    $(\lambda_g \colon F_\lambda(\bullet^g) \xrightarrow{\sim}
    F_\lambda)_{g\in G}$
    of isomorphisms of triangle functors such that the following
    equalities hold true, for all $M\in {\mathcal{C}_A}$ and
    $g,h\in G$,
    \[
    \left\{
      \begin{array}{rcl}
        (\lambda_e)_M & = & \mathbbm 1_M \\
        (\lambda_g)_M \circ (\lambda_h)_{M^g} & = & (\lambda_{gh})_M\,.
      \end{array}\right.
    \]
  \item For all $M\in {\mathcal{C}_A}$ there is a functorial isomorphism
    $F_\cdot F_\lambda M \simeq \oplus_{g\in G}M^g$ in ${\mathcal{C}_A}$.
  \item For all $M,M'\in {\mathcal{C}_A}$ there is an isomorphism
    $\oplus_{g\in G}{\mathcal{C}_A}(M,M'^g) \xrightarrow{\sim} {\mathcal{C}_{A*G}}(F_\lambda
    M,F_\lambda M')$ functorial in $M$ and $M'$.
  \item $(F_\cdot N)^g  \simeq F_\cdot N$ for all $N\in {\mathcal{C}_{A*G}}$.
  \item For all indecomposable $M\in {\mathcal{C}_A}$, if $G_M = \{e\}$, then
    $F_\lambda M$ is indecomposable.
  \end{enumerate}
\end{prop}
\begin{proof}
  (1) Consider the following subcategories of $\mathrm{per}(A)$ and
  $\mathrm{per}(A*G)$, respectively,
  \[
  \left\{
    \begin{array}{l}
      \mathcal F_A = \mathcal D^{\leqslant 0}(A) \cap\,^\perp
      \mathcal D^{\leqslant -2}(A) \cap \mathrm{per}(A)\,, \\
      \mathcal F_{A*G} = \mathcal D^{\leqslant 0}({A*G}) \cap\,^\perp
      \mathcal D^{\leqslant -2}({A*G}) \cap \mathrm{per}(A*G) \,,\\
    \end{array}\right.
  \]
  where $\mathcal D^{\leqslant i}$ stands for the full subcategory of
  objects $Z$ such that $H^pZ$ vanishes for $p>i$ and $^\perp\bullet$
  denotes the full subcategory of objects with no non zero morphism to
  any object in a given subcategory $\bullet$. Following
  \cite[Proposition 2.9]{A08}, the following restrictions of $\pi_A$
  and $\pi_{A*G}$ are equivalences
  \[
  \left\{
    \begin{array}{l}
      \mathcal F_A \to {\mathcal{C}_A}\,, \\
      \mathcal F_{A*G} \to {\mathcal{C}_{A*G}}\,.
    \end{array}\right.
  \]
  Since the pairs $(-\otimes^{\mathrm L}_A(A*G),\mathrm{Rres})$ and
  $(\mathrm{Rres},-\otimes^{\mathrm L}_A(A*G))$ are adjoint, it is sufficient to prove
  the following assertions in order to prove that
  $(F_\lambda,F_\cdot)$ and $(F_\cdot,F_\lambda)$ are adjoint,
  \begin{enumerate}[(a)]
  \item $X\otimes^{\mathrm L}_A (A*G) \in \mathcal F_{A*G}$ for all
    $X\in \mathcal F_A$, and
  \item $\mathrm{Rres} Y\in \mathcal F_A$ for all $Y\in \mathcal F_{A*G}$.
  \end{enumerate}
  Here is a proof of these assertions. First, recall that
  $-\otimes^{\mathrm L}_A(A*G)$ and $\mathrm{Rres}$ map $\mathrm{per}(A)$ and
  $\mathrm{per}(A*G)$ into $\mathrm{per}(A*G)$ and $\mathrm{per}(A)$, respectively. Next, for
  all $i\in \mathbb Z$,
  \begin{equation}
    \label{eq:1}
    \mathrm{Rres}(\mathcal D^{\leqslant i}(A*G)) \subseteq \mathcal
    D^{\leqslant i}(A)\,.
  \end{equation}
  Next, since $A$ is concentrated in non positive degrees, the
  following holds true for all $i\in \mathbb Z$,
  \begin{equation}
    \label{eq:2}
    (\forall M\in \mathcal D^{\leqslant i}(A))\ \ X\otimes^{\mathrm L}_A (A*G)\in \mathcal D^{\leqslant i}(A)\,.
  \end{equation}
  Because of the adjunction $(-\otimes^{\mathrm L}_A(A*G),\mathrm{Rres})$, it follows from (\ref{eq:1}) that
  \begin{equation}
    \label{eq:3}
    (\forall X\in \,^\perp \mathcal D^{\leqslant -2}(A))\ \
    X\otimes^{\mathrm L}_A (A*G) \in \,^\perp \mathcal D^{\leqslant -2}(A*G)\,.
  \end{equation}
  Because of the adjunction
  $(\mathrm{Rres},-\otimes^{\mathrm L}_A(A*G))$, it follows from
  (\ref{eq:2}) that
  \begin{equation}
    \label{eq:4}
    (\forall Y\in \,^\perp\mathcal D^{\leqslant -2}(A*G))\ \ \mathrm{Rres} Y \in \,^\perp \mathcal D^{\leqslant -2}(A)\,.
  \end{equation}
  The assertions (a) and (b) above are consequences of (\ref{eq:1}),
  (\ref{eq:2}), (\ref{eq:3}) and (\ref{eq:4}). Thus,
  $(F_\lambda,F_\cdot)$ and $(F_\cdot,F_\lambda)$ are adjoint.

  (2), (3), (4), (5) and (6) follow from (1) and from
  Lemma~\ref{sec:adjunctions-1}.

  (7) The arguments in the proof of \cite[3.5]{MR654725} apply
  literally here: for a direct sum decomposition
  $F_\lambda M = N\oplus N'$ such that $N$ is non zero, there exists a
  non empty subset $H$ of $G$ such that
  $F_\cdot N \simeq \oplus_{g\in H} M^g$ and
  $F_\cdot N'\simeq \oplus_{g\in G\backslash H}M^g$ because $M$ is
  indecomposable and because of (4); hence $gH\subseteq H$ for all
  $g\in G$ (because $(F_\cdot N)^g\simeq F_\cdot N$); accordingly,
  $H=G$, and hence $N'=0$.
\end{proof}

\subsection{Comparison of cluster tilting subcategories of \texorpdfstring{${\mathcal{C}_A}$}{C(A)} and
  \texorpdfstring{${\mathcal{C}_{A*G}}$}{C(A*G)}}
\label{sec:cy-categ-relat}

This subsection compares the cluster tilting subcategories of ${\mathcal{C}_A}$ to
the ones of ${\mathcal{C}_{A*G}}$. Using completely analogous, and hence omitted,
arguments, the same conclusions are obtained when replacing cluster
tilting subcategories by maximal rigid ones. The comparison is based
on the following more general property where cluster tilting
subcategories are used instead of cluster tilting objects for the ease
of exposition.
\begin{prop}
  \label{sec:cy-categ-relat-1}
  Let $\mathcal C$ and $\mathcal C'$ be triangulated ${\mathbbm{k}}$-linear
  categories which are $\mathrm{Hom}$-finite, Krull-Schmidt and
  $2$-Calabi-Yau. Assume that there exist triangle functors
  $\xymatrix{ F\colon \mathcal C \ar@<+0.5ex>[rr] && \mathcal C'
    \colon F' \ar@<+0.5ex>[ll]}$ such that
  \begin{itemize}
  \item $(F,F')$ and $(F',F)$ are adjoint and
  \item for all $M\in \mathcal C$ and $N\in \mathcal C'$ the unit
    morphisms $M \to F'FM$ and $N\to FF'N$ split in $\mathcal C$
    and $\mathcal C'$, respectively.
  \end{itemize}
  Then, the two assignments $\mathcal T \mapsto \mathrm{add}(F \mathcal T)$
  and $\mathcal T' \mapsto \mathrm{add}(F' \mathcal T')$ induce mutually
  inverse bijections between
  \begin{itemize}
  \item cluster tilting subcategories $\mathcal T$ of $\mathcal C$ which are
    stable under $F' F$, and
  \item cluster tilting subcategories $\mathcal T'$ of $\mathcal C'$ which are
    stable under $F F'$.
  \end{itemize}
  The same result holds true when replacing cluster tilting
  subcategories by maximal rigid ones.
\end{prop}
\begin{proof}
  By adjunction, $F$ and $F'$ preserve functorially finite
  subcategories.

  Let $\mathcal T$ be a cluster tilting subcategory of $\mathcal C$
  stable under $F'F$. Then,
  $\mathcal C'(F M, \Sigma F M') \simeq \mathcal C(M,\Sigma F'FM') =
  0$
  for all $M,M'\in \mathcal T$ (note that, by assumption,
  $F'FM'\in \mathcal T$). Hence, $\mathcal C'(N_1,\Sigma N_2)=0$, for
  all $N_1,N_2\in \mathrm{add}( F \mathcal T)$. Moreover, for all
  $N\in \mathcal C'$, if $\mathcal C'(N',\Sigma N) = 0$ for all
  $N'\in \mathrm{add}(F \mathcal T)$, then
  $\mathcal C(M,\Sigma F' N) \simeq \mathcal C'(F M,\Sigma N)= 0$ for
  all $M\in \mathcal T$; therefore, $F' N \in \mathcal T$, and hence
  $F F' N \in F\mathcal T$; since $N\in \mathrm{add}(F F' N)$, the object $N$
  lies in $\mathrm{add}(F \mathcal T)$. Thus, $\mathrm{add}(F \mathcal T)$ is a
  cluster tilting subcategory of $\mathcal C'$.

  \emph{Mutatis mutandis,} the same arguments show that, for all
  cluster tilting subcategories $\mathcal T'$ of $\mathcal C'$ stable
  under $F F'$, the subcategory $\mathrm{add}(F' \mathcal T')$ of $\mathcal C$
  is cluster tilting.

  Finally, given a subcategory $\mathcal T$ of $\mathcal C$ (or
  $\mathcal T'$ of $\mathcal C'$) such that $\mathcal T$ (or
  $\mathcal T'$) is closed under direct sums, direct summands and
  $F' F$ (or $F F'$), then $\mathrm{add}(F' \mathrm{add}(F \mathcal T)) = \mathcal T$
  (or $\mathrm{add}(F \mathrm{add}(F' \mathcal T')) = \mathcal T'$) because, for all
  $M\in \mathcal C$ (or $N\in \mathcal C'$), the unit morphism
  $M\to F' FM$ (or $N \to F F'N$, respectively) splits.
\end{proof}

In view of Proposition~\ref{sec:adjunct-betw-ca-2},
Proposition~\ref{sec:cy-categ-relat-1} applies to the current
situation.
\begin{cor}
  \label{sec:clust-tilt-objects-3}
  Assume setting~\ref{sec:main-results-5-bis}. The assignments
  $\mathcal T \mapsto \mathrm{add}(F_\lambda \mathcal T)$ and
  $\mathcal T' \mapsto \mathrm{add}(F_\cdot \mathcal T')$ induce mutually
  inverse bijections between
  \begin{itemize}
  \item cluster tilting subcategories $\mathcal T$ of ${\mathcal{C}_A}$ which are
    stable under the action of $G$, and
  \item cluster tilting subcategories $\mathcal T'$ of ${\mathcal{C}_{A*G}}$ which are
    stable under $F_\lambda F_\cdot$.
  \end{itemize}
  The same result holds true if cluster tilting subcategories are
  replaced by maximal rigid ones.
\end{cor}
\begin{proof}
  Note that, since $F_\cdot F_\lambda M \simeq \oplus_{g\in G}M^g$,
  for all $M\in {\mathcal{C}_A}$, a full subcategory of
  ${\mathcal{C}_A}$ which is closed under direct sums and direct
  summands is stable under the action of $G$ if and only if it is
  stable under $F_\cdot F_\lambda$. The conclusion therefore follows
  from Proposition~\ref{sec:cy-categ-relat-1} applied to
  $\mathcal C = {\mathcal{C}_A}$, $\mathcal C' = {\mathcal{C}_{A*G}}$,
  $F = F_\lambda$ and $F'= F_\cdot$.
\end{proof}

Given cluster tilting (or maximal rigid) objects $T,T'$ of ${\mathcal{C}_A}$ and
${\mathcal{C}_{A*G}}$, respectively, such that $\mathrm{add}(T') = \mathrm{add}(F_\lambda T)$ as in
Corollary~\ref{sec:clust-tilt-objects-3}, it is natural to ask whether
there exists a relation between $\mathrm{End}(T)$ and $\mathrm{End}(T')$. It is
provided by the following result.

For a better understanding, recall (\cite[3.1]{RR85}) that when $C$ is
a ${\mathbbm{k}}$-linear category endowed with an action of $G$ (on the right) by
automorphisms, the \emph{skew group} category $C*G$ is the
idempotent completion of the orbit category $C[G]$ (see also
\cite[Definition 2.3]{MR2844757} for an approach avoiding idempotent
completion). Recall that the \emph{orbit category} $C[G]$ has the same
objects as $C$, that, for all $M,M'\in C$, its space of morphisms
$M\to M'$ is defined by $C[G](M,M') = \oplus_{g\in G}C(M,M'^g)$, and
that the composition of morphisms is induced naturally by the one in
$C$.  Recall also that, for a ${\mathbbm{k}}$-linear category $C$, the piece of
notation $\mathrm{mod}(C)$ stands for the category of finitely
presented contravariant functors from $C$ to finite dimensional vector
spaces.
\begin{prop}
  \label{sec:clust-tilt-objects-4}
  Assume setting~\ref{sec:main-results-5-bis}. Let $\mathcal T$ be a
  cluster tilting subcategory of $\mathcal C_A$ stable under the
  action of $G$. Let $\mathcal T' = \mathrm{add}(F_\lambda \mathcal T)$. Then,
  \begin{enumerate}
  \item the skew group category $\mathcal T * G$ is Morita
    equivalent to $\mathcal T'$ and
  \item there exists a cluster tilting object $T\in {\mathcal{C}_A}$ such that
    $T^g = T$ for all $g\in G$ and $\mathcal T = \mathrm{add}(T)$. For any
    such $T$, the algebras $\mathrm{End}(T)*G$ and $\mathrm{End}(F_\lambda T)$ are
    Morita equivalent.
  \end{enumerate}
  The same result holds true when cluster tilting subcategories are
  replaced by maximal rigid ones.
\end{prop}
\begin{proof}
  (1) It follows from Proposition~\ref{sec:adjunct-betw-ca-2} that
  the categories $\mathcal T[G]$ and $F_\lambda\mathcal T$ are
  equivalent. Taking idempotent completions yields that $\mathcal T *
  G$ and $\mathrm{add}(F_\lambda\mathcal T)$ are equivalent. Hence,
  \begin{equation}
    \label{eq:11}
    \mathrm{mod}(\mathcal T*G) \simeq \mathrm{mod}(\mathcal T')\,.
  \end{equation}

  (2) Let $\{T_1,\ldots,T_n\}$ be a complete set of representatives of
  the $G$-orbits of the indecomposable objects of $\mathcal T$. Let
  $T$ be the object $\bigoplus\limits_{
    \begin{subarray}{c}
      1\leqslant i \leqslant n \\
      g\in G
    \end{subarray}} T_i^g$
  of ${\mathcal{C}_A}$.  Then $T^g = T$ for all $g\in G$ and
  $\mathcal T = \mathrm{add}(T)$.  The assignment $(g,f) \mapsto f^{g^{-1}}$
  defines an action of $G$ on $\mathrm{End}(T)$ by algebra automorphisms,
  which defines $\mathrm{End}(T)*G$.

  Since $\mathcal T$ is idempotent complete, the following holds true
  (\cite[Proposition 3.4]{RR85}),
  \begin{equation}
    \label{eq:12}
    \mathrm{mod}(\mathcal T * G) \simeq \mathrm{mod}(\mathcal T)*G\,.
  \end{equation}
  By construction of $T$, the categories
  $\mathrm{mod}(\mathcal T)$ and $\mathrm{mod}(\mathrm{End}(T))$ are
  equivalent and
  \begin{equation}
    \label{eq:13}
    \mathrm{mod}(\mathrm{End}(T))*G \simeq \mathrm{mod}(\mathcal T)*G\,.
  \end{equation}
  Note also that (\cite[3.1]{RR85})
  \begin{equation}
    \label{eq:14}
    \mathrm{mod}(\mathrm{End}(T))*G \simeq \mathrm{mod}(\mathrm{End}(T)*G)\,.
  \end{equation}
  Finally, since $\mathcal T' = \mathrm{add}(F_\lambda T)$,
  \begin{equation}
    \label{eq:15}
    \mathrm{mod}(\mathcal T') \simeq \mathrm{mod}(\mathrm{End}(F_\lambda T))\,.
  \end{equation}
  The conclusion therefore follows from (\ref{eq:11}), (\ref{eq:12}),
  (\ref{eq:13}), (\ref{eq:14}) and (\ref{eq:15}).
\end{proof}

\subsection{Cluster tilting objects of the first kind}
\label{sec:clust-tilt-objects}

The considerations of subsection~\ref{sec:cy-categ-relat} raise the
question whether a given cluster tilting or maximal rigid object of
${\mathcal{C}_{A*G}}$ lies in the essential image of $F_\lambda \colon {\mathcal{C}_A} \to
{\mathcal{C}_{A*G}}$. For instance, $F_\lambda \overline A = \overline{A*G}$. An
object $N\in {\mathcal{C}_{A*G}}$ is called here of the \emph{first kind} with
respect to $F_\lambda$ if, for all indecomposable direct summands $N'$
of $N$, there exists $M'\in {\mathcal{C}_A}$ such that $F_\lambda M'\simeq
N'$. This terminology comes from the theory of Galois coverings of
algebras (see \cite{MR896100}).

This subsection investigates cluster tilting  or maximal rigid
objects of ${\mathcal{C}_{A*G}}$ which are of the first kind with respect to
$F_\lambda$. Clearly, further assumptions are needed for this
investigation to be fruitful because, although
$F_\lambda\overline A = \overline{A*G}$, the indecomposable direct
summands of $\overline{A*G}$ need not lie in the essential image of
$F_\lambda$. This is because $F_\lambda$ need not take an
indecomposable direct summand of $\overline A$ to an indecomposable
object of ${\mathcal{C}_{A*G}}$. In view of
Proposition~\ref{sec:adjunct-betw-ca-2}, having some information on
the stabilisers of the indecomposable direct summands of $\overline A$
would be useful. This is provided by the following lemma.

For all $g\in G$, the mapping $A\to A$ given by $a
\mapsto \,^ga$ induces an isomorphism in ${\mathcal{C}_A}$ denoted by $\alpha_g$,
\begin{equation}
  \label{eq:23}
  \alpha_g \colon \overline A \to \overline A^g\,.
\end{equation}
Recall \cite[Proposition 2.9]{A08} that the mapping $A\to \mathrm{Hom}_A(A,A)$
given by $x \mapsto (a \mapsto xa)$ induces an isomorphism of algebras
\begin{equation}
  \label{eq:25}
  H^0(A) \xrightarrow{\sim} \mathrm{End}(\overline A)\,.
\end{equation}
In particular, for all direct summands $X$ of $\overline A$, the
$\mathrm{End}(\overline A)$-module ${\mathcal{C}_A}(\overline A,X)$ is projective and its
restriction of scalars along \eqref{eq:25} is isomorphic to
$\varepsilon H^0(A)$ for some idempotent $\varepsilon\in H^0(A)$. Any idempotent of
$H^0(A)$ may be obtained in this way.
\begin{lem}
  \label{sec:action-g-ca-1}
  Assume setting~\ref{sec:main-results-5-bis}. Let $\varepsilon\in H^0(A)$ be an
  idempotent. Let $X$ be a direct summand of $\overline A$ such that
  $\varepsilon H^0(A)$ is isomorphic, as an $H^0(A)$-module, to the restriction
  of scalars of ${\mathcal{C}_A}(\overline A, X)$ along (\ref{eq:25}). Let
  $g\in G$. Let $Y$ be a direct summand of $\overline A$ such that
  $X^g \simeq Y$ (see (\ref{eq:23})). Then, $^g\varepsilon H^0(A)$ is
  isomorphic to the restriction of scalars of ${\mathcal{C}_A}(\overline A,Y)$
  along (\ref{eq:25}).
\end{lem}
\begin{proof}
  The proof is based on the following diagram where all the arrows
  are isomorphisms. Note that the top
  horizontal arrow is given by the action of $G$ on $A$ and the
  leftmost bottom horizontal arrow is given by the action of $G$ on ${\mathcal{C}_A}$.
  \begin{equation}
    \label{eq:26}
  \xymatrix{
    H^0(A) \ar[rrr]^{x \mapsto \,^{g^{-1}}x} \ar[d]_{~(\ref{eq:25})} &&& H^0(A) \ar[d]^{~(\ref{eq:25})} \\
    \mathrm{End}(\overline A) \ar[r]_{-^g} & \mathrm{End}(\overline A^g)
    \ar[rr]_{\alpha_g^{-1}\circ - \circ \alpha_g} && \mathrm{End}(\overline
    A)}
  \end{equation}

  This diagram is commutative as shown by a direct computation using
  that
  \begin{itemize}
  \item $\mathrm{End}(\overline A) \xrightarrow{-^g} \mathrm{End}(\overline A^g)$ is
    induced by the identity mapping $\mathrm{Hom}_A(A,A) \to \mathrm{Hom}_A(A,A)$,
  \item
    $\mathrm{End}(\overline A^g) \xrightarrow{\alpha_g^{-1} \circ - \circ
      \alpha_g} \mathrm{End}(\overline A)$
    is induced by the mapping $\mathrm{Hom}_A(A,A) \to \mathrm{Hom}_A(A,A)$ which
    assigns to any $f\in \mathrm{Hom}_A(A,A)$  the mapping
    $a \mapsto \,^{g^{-1}}f(\,^ga)$.
  \end{itemize}
  
  On one hand,
  \begin{itemize}
  \item by assumption, the restriction of scalars of
    ${\mathcal{C}_A}(\overline A,X)$ along the rightmost vertical arrow of
    (\ref{eq:26}) is isomorphic to $\varepsilon H^0(A)$ as an $H^0(A)$-module and
  \item the following mapping is an isomorphism of $H^0(A)$-modules
    from the restriction of scalars of $\varepsilon H^0(A)$ along the top
    horizontal arrow of (\ref{eq:26}) to $^g\varepsilon H^0(A)$,
    \[
    \begin{array}{rcl}
      \varepsilon H^0(A) & \to & \,^g\varepsilon H^0(A) \\
      x & \mapsto & \,^gx\,.
  \end{array}
  \]
  \end{itemize}
  On the other hand, denoting by $\lambda \colon X \to \,^gY$ any
  isomorphism in ${\mathcal{C}_A}$, 
  \begin{itemize}
  \item the following mapping is an isomorphism of $\mathrm{End}(\overline
    A^g)$-modules from the restriction of scalars of ${\mathcal{C}_A}(\overline A,
    X)$ along the arrow $\mathrm{End}(\overline A^g) \to \mathrm{End}(\overline A)$ of
    (\ref{eq:26}) to ${\mathcal{C}_A}(\overline A^g,Y^g)$,
    \[
    \begin{array}{rcl}
      {\mathcal{C}_A}(\overline A,X) & \to & {\mathcal{C}_A}(\overline A^g, Y^g) \\
      f & \mapsto & \lambda \circ f \circ \alpha_g^{-1}\,,
    \end{array}
    \]
  \item and the following mapping is an isomorphism of $\mathrm{End}(\overline
    A)$-modules from the restriction of scalars of ${\mathcal{C}_A}(\overline
    A^g,Y^g)$ along the arrow $\mathrm{End}(\overline A)\to \mathrm{End}(\overline
    A^g)$  of (\ref{eq:26}) to ${\mathcal{C}_A}(\overline A,Y)$,
    \[
    \begin{array}{rcl}
      {\mathcal{C}_A}(\overline A^g,Y^g) & \to & {\mathcal{C}_A}(\overline A, Y) \\
      f & \mapsto & f^{g^{-1}}\,.
    \end{array}
    \]
  \end{itemize}
  Since the diagram (\ref{eq:26}) commutes, it follows that the
  restriction of scalars of ${\mathcal{C}_A}(\overline A, Y)$ along
  (\ref{eq:25}) is isomorphic to $^g\varepsilon H^0(A)$.    
\end{proof}

The main results of this subsection prove that being of the first kind
with respect to $F_\lambda$ is an invariant property within a mutation
class of cluster tilting or of maximal rigid objects. This fact is
based on the following lemma. It is already proved in other contexts
for Galois coverings instead of skew group algebras, namely, in
\cite[Lemma 3.2]{MR2423814} in the context of module categories, in
\cite[Lemma 4.2]{MR2746678} in the context of derived categories of
finite dimensional algebras and in \cite[Lemma 6.2]{lemeur} in the
context of cluster categories of piecewise hereditary algebras.
\begin{lem}
  \label{sec:clust-tilt-objects-1}
  Assume setting~\ref{sec:main-results-5-bis}. Let
  $\Delta \colon X\to M \to Y \to \Sigma X$ be a triangle in ${\mathcal{C}_{A*G}}$,
  such that
  \begin{itemize}
  \item there exists $\widetilde X\in {\mathcal{C}_A}$ such that $F_\lambda
    \widetilde X \simeq X$,
  \item there exists a direct sum decomposition into indecomposable
    objects $M = M_1 \oplus \cdots \oplus M_t$ in ${\mathcal{C}_{A*G}}$, and there
    exist $\widetilde M_1,\ldots,\widetilde M_t\in {\mathcal{C}_A}$ such that
    $F_\lambda \widetilde M_i\simeq M_i$ for all $i$,
  \item ${\mathcal{C}_{A*G}}(Y,\Sigma M) =0$.
  \end{itemize}
  Then, there exist $\widetilde Y\in {\mathcal{C}_A}$ and $g_1,\ldots,g_t\in G$,
  and there exists a triangle
  $\widetilde X \to \bigoplus_{i=1}^t\widetilde M_i^{g_i} \to
  \widetilde Y \to \Sigma \widetilde X$
  in ${\mathcal{C}_A}$ whose image under $F_\lambda$ is isormorphic to
  $\Delta$. In particular, $F_\lambda \widetilde Y \simeq Y$.
\end{lem}
\begin{proof}
  The proof of \cite[Lemma 4.2]{MR2746678} applies here provided that
  the occurrences of $\mathcal D^b(\mathrm{mod}\,\mathcal C)$ and
  $\mathcal D^b(\mathrm{mod}\,A)$ there are replaced by ${\mathcal{C}_A}$ and
  ${\mathcal{C}_{A*G}}$ here. Note that, by Proposition~\ref{sec:adjunct-betw-ca-2},
  the functor $F_\lambda\colon {\mathcal{C}_A} \to {\mathcal{C}_{A*G}}$ has all the needed
  properties of the functor that is also denoted by $F_\lambda$ in the
  proof of \cite[Lemma 4.2]{MR2746678}.
\end{proof}

As explained at the beginning of
subsection~\ref{sec:clust-tilt-objects}, it is preferable to have some
control on the stabilisers of the indecomposable direct summands of
the cluster tilting or maximal rigid objects of ${\mathcal{C}_A}$. For this
purpose, the following lemma is useful. Like
Lemma~\ref{sec:clust-tilt-objects-1}, the following result is already
known in the framework of Galois coverings in \cite[Proposition
4.6]{MR2746678} for tilting complexes over piecewise hereditary
algebras and in \cite[Proposition 6.7]{lemeur} for cluster tilting
objects of cluster categories of piecewise hereditary algebras.
\begin{lem}
  \label{sec:clust-tilt-objects-8}
  Assume the setting of Lemma~\ref{sec:clust-tilt-objects-1}. Assume,
  moreover, that
  \begin{itemize}
  \item $X$ is indecomposable and $X\not \simeq M_i$ for all
    $i\in \{1,\ldots,t\}$,
  \item ${\mathcal{C}_{A*G}}(M,\Sigma X) = 0$,
  \item the stabiliser of $\widetilde X$ is trivial.
  \end{itemize}
  Then, for all triangles
  $\widetilde X \to \oplus_{i=1}^t\widetilde M_i^{g_i} \to
  \widetilde Y\to \Sigma \widetilde Y$
  of ${\mathcal{C}_A}$ whose images under $F_\lambda$ are isomorphic to $\Delta$, 
  the stabiliser of $\widetilde Y$ is trivial.
\end{lem}
\begin{proof}
  The proof of \cite[Proposition 6.7]{lemeur} may be adapted to the
  current situation. Here are the details. Let
  $\widetilde \Delta\colon \widetilde X \to \oplus_{i=1}^t\widetilde
  M_i^{g_i} \to \widetilde Y\to \Sigma \widetilde X$
  be a triangle of ${\mathcal{C}_A}$ whose image under $F_\lambda$ is isomorphic to
  $\Delta$. Let $w\colon \widetilde Y^g \to \widetilde Y$ be an
  isomorphism in ${\mathcal{C}_{A*G}}$, where $g\in G$. For all $h\in G$, the space
  ${\mathcal{C}_A}(\widetilde M_i^h,\Sigma \widetilde X)$ is isomorphic to a
  subspace of ${\mathcal{C}_{A*G}}(M,\Sigma X)$ (see
  Proposition~\ref{sec:adjunct-betw-ca-2}), since the latter is zero,
  there exist morphisms $u,v,u',v'$ fitting into commutative diagrams
  as follows
  \[
  \xymatrix@C=2.5ex{ \widetilde X^g \ar[r] \ar[d]_u &
    \bigoplus_{i=1}^t\widetilde M_i^{g_ig} \ar[r] \ar[d]_v& \widetilde Y^g
    \ar[r] \ar[d]_w&
    \Sigma \widetilde X^g \ar[d]^{\Sigma u} &
    \widetilde X \ar[r] \ar[d]_{u'} & \bigoplus_{i=1}^t\widetilde
    M_i^{g_i} \ar[r] \ar[d]_{v'}
    & \widetilde Y \ar[r] \ar[d]_{w^{-1}} & \Sigma \widetilde X
    \ar[d]^{\Sigma u'}
    \\
    \widetilde X \ar[r] & \bigoplus_{i=1}^t\widetilde M_i^{g_i} \ar[r]
    & \widetilde Y \ar[r]  & \Sigma \widetilde X &
    \widetilde X^g \ar[r] &
    \bigoplus_{i=1}^t\widetilde M_i^{g_ig} \ar[r] & \widetilde Y^g
    \ar[r] &
    \Sigma \widetilde X^g
  }
  \]
  Therefore, the following diagram commutes
  \[
  \xymatrix{
    \widetilde X \ar[r] \ar[d]_{u\circ u'-\mathbbm 1} & \bigoplus_{i=1}^t\widetilde
    M_i^{g_i} \ar[r] \ar[d]_{v\circ v'- \mathbbm 1}
    & \widetilde Y \ar[r] \ar[d]_{0} & \Sigma \widetilde X \ar[d]^{\Sigma(u
      \circ u'- \mathbbm 1)}
    \\
    \widetilde X \ar[r] & \bigoplus_{i=1}^t\widetilde M_i^{g_i} \ar[r]
    & \widetilde Y \ar[r]  & \Sigma \widetilde X
  }
  \]
  and hence $\Sigma(u \circ u' - \mathbbm 1)$ factors through
  $\Sigma \widetilde X \to \oplus_{i=1}^t\Sigma \widetilde M_i^{g_i}$.
  Note that $\widetilde X \not\simeq \widetilde M_i^{g_i}$ because
  $X\not \simeq M_i$, $X \simeq F_\lambda \widetilde X$ and
  $F_\lambda \widetilde M_i \simeq M_i$, for all
  $i\in \{1,\ldots,t\}$. Since $\widetilde X$ is indecomposable, it
  follows that $\Sigma(u \circ u' - \mathbbm 1)$ lies in the radical
  of the algebra $\mathrm{End}(\Sigma\widetilde X)$, and hence $u \circ u'$ is
  an automorphism of $\widetilde X$. Thus,
  $u\colon \widetilde X^g \to \widetilde X$ is an isomorphism. Given
  that $\widetilde X$ has trivial stabiliser, it follows that
  $g=e$. This proves that $\widetilde Y$ has trivial stabiliser.
\end{proof}

Combining Lemmas~\ref{sec:clust-tilt-objects-1} and
\ref{sec:clust-tilt-objects-8} yields the following result. Recall
that, given cluster tilting objects $T,T'\in {\mathcal{C}_{A*G}}$, the
object $T'$ is called a \emph{mutation} of $T$ if there exist
indecomposable direct summands $X$ and $Y$ of $T$ and $T'$,
respectively, such that the set of isomorphism classes of
indecomposable direct summands of $T$ not isomorphic to $X$ is equal
to the set of isomorphism classes of indecomposable direct summands of
$T'$ not isomorphic to $Y$. In such a case, there exists a triangle
$X\to M \to Y \to \Sigma X$ in ${\mathcal{C}_{A*G}}$, where
$M\in \mathrm{add}(T) \cap \mathrm{add}(T')$. The \emph{mutation
  class} of $T$ is the transitive closure of $\{T\}$ under
mutation. The same concepts are defined for maximal rigid objects. See
\cite{MR2385669} for more details.
\begin{lem}
  \label{sec:clust-tilt-objects-2}
  Assume setting~\ref{sec:main-results-5-bis}. Let $T,T'\in {\mathcal{C}_{A*G}}$ be
  cluster tilting objects in the same mutation class. The following
  assertions are equivalent.
  \begin{enumerate}[(i)]
  \item All indecomposable direct summands of $T$ lie in the essential
    image of
    $F_\lambda$.
  \item All indecomposable direct summands of $T'$ lie in the
    essential image of $F_\lambda$.
  \end{enumerate}
  When these assertions are true, then the following ones are
  equivalent.
  \begin{enumerate}[(i)]
    \setcounter{enumi}{2}
  \item $G_X=\{e\}$ for all indecomposable objects $X\in {\mathcal{C}_A}$ such
    that $F_{\lambda}X$ is isomorphic to a direct summand of $T$.
  \item $G_X=\{e\}$ for all indecomposable objects $X\in {\mathcal{C}_A}$ such
    that $F_{\lambda}X$ is isomorphic to a direct summand of $T'$.
  \end{enumerate}
  The same result holds true when cluster tilting objects are replaced
  by maximal rigid ones.
\end{lem}
\begin{proof}
  This follows from
  Lemmas~\ref{sec:clust-tilt-objects-1} and \ref{sec:clust-tilt-objects-8}.
\end{proof}

Using Lemma~\ref{sec:clust-tilt-objects-2}, it is possible to prove a
first result answering the question stated at the beginning of
subsection~\ref{sec:clust-tilt-objects} and regarding cluster tilting
objects of ${\mathcal{C}_{A*G}}$ which are of the first kind with respect to
$F_\lambda$. Note that the action of $G$ on $H^0(A)$ defines an action
on the finite dimensional $H^0(A)$-modules like the action of $G$ on
$A$ defines an action on dg $A$-modules.

\begin{prop}
  \label{sec:clust-tilt-objects-6}
  Assume setting~\ref{sec:main-results-5-bis}. Assume that $G$ acts freely
  on the isomorphism classes of indecomposable projective
  $H^0(A)$-modules. Then, for all basic cluster tilting objects
  $T = T_1\oplus \cdots \oplus T_n$ of ${\mathcal{C}_{A*G}}$ in the mutation class of
  $\overline{A*G}$ and where $T_1,\ldots,T_n$ are indecomposable,
  \begin{enumerate}
  \item there exist $\widetilde T_1,\ldots,\widetilde T_n\in {\mathcal{C}_A}$ such
    that $F_\lambda \widetilde T_i \simeq T_i$ and $G_{\widetilde
      T_i}= \{e\}$, for all $i$.
  \end{enumerate}
  Moreover, for any such $\widetilde T_1, \ldots \widetilde T_n$,
  denoting
  $\bigoplus\limits_{ 1\leqslant i \leqslant n\,, g\in G} \widetilde
  T_i^g$ by $\widetilde T$,
  \begin{enumerate}
    \setcounter{enumi}{1}
  \item $\widetilde T$ is a cluster tilting object in ${\mathcal{C}_{A*G}}$ such that
    $\mathrm{add}(T) = F_\lambda\mathrm{add}(\widetilde T)$ and $\widetilde T^g =
    \widetilde T$ for all $g\in G$, and
  \item $F_\lambda$ induces a Galois covering with group $G$ from
    $\mathrm{End}(\widetilde T)$ to $\mathrm{End}(T)$.
  \end{enumerate}
  The same result holds true if cluster tilting objects are replaced
  by maximal rigid ones.
\end{prop}
\begin{proof}
  (1) First, the stated property holds true when
  $\mathrm{add}(T) = \mathrm{add}(\overline{A*G})$. Indeed,
  \begin{itemize}
  \item $F_\lambda\overline A = \overline{A*G}$,
  \item by assumption, every indecomposable direct summand of
    $\overline A$ has a trivial stabiliser (see
    Lemma~\ref{sec:action-g-ca-1}) and consequently,
  \item $F_\lambda M$ is isomorphic to an indecomposable direct
    summand of $\overline{A*G}$ for all indecomposable direct summands
    $M$ of $\overline A$ (see
    Proposition~\ref{sec:adjunct-betw-ca-2}, part (7)).
  \end{itemize}
    Next, it follows from Lemma~\ref{sec:clust-tilt-objects-2} that
    property (1) is preserved within a mutation class.

    (2) By construction,
    \begin{itemize}
    \item $\mathrm{add}(T) = F_\lambda\mathrm{add}(\widetilde T)$ and $\widetilde T^g =
      \widetilde T$ for all $g\in G$. Moreover, 
    \item $F_\lambda F_\cdot (F_\lambda \widetilde T) \simeq
      F_\lambda(\widetilde T^{\mathrm{Card}(G)})$ (see
      Proposition~\ref{sec:adjunct-betw-ca-2}, part (4)), and hence
      $\mathrm{add}(F_\lambda\widetilde T)$ is stable under $F_\lambda F_\cdot$.
    \end{itemize}
    Therefore, $\mathrm{add}(\widetilde T)$ is a cluster tilting subcategory
    (see Corollary~\ref{sec:clust-tilt-objects-3}), and hence
    $\widetilde T$ is a cluster tilting object of ${\mathcal{C}_A}$.

    (3) This was proved already in other contexts (\cite[Lemma
    2.2]{MR2423814} and \cite[Lemma 4.8]{MR2746678}) and may be proved similarly
    here. Here are some explanations for the ease of reading.

    Here, $\mathrm{End}(\widetilde T)$ and $\mathrm{End}(T)$ are identified with the
    categories whose objects sets are equal to the sets
    $\{\widetilde T_i^g\ |\ 1\leqslant i \leqslant n\,,g\in G\}$ and
    $\{T_1,\ldots,T_n\}$, respectively.  For all
    $i\in \{1,\ldots,n\}$, let $\mu_i \colon F_\lambda \widetilde T_i
    \to T_i$
    be an isomorphism in ${\mathcal{C}_{A*G}}$. Then, define a functor
    \begin{equation}
      \label{eq:27}
    \begin{array}{rcl}
      \mathrm{End}(\widetilde T) & \to & \mathrm{End}(T)
    \end{array}
  \end{equation}
  as follows, where $g,h\in G$, $i,j\in \{1,\ldots,n\}$,
  \begin{itemize}
  \item the object $\widetilde T_i^g$ is mapped onto $T_i$,
  \item any morphism
  $u \colon \widetilde T_i^g \to \widetilde T_j^h$ is mapped onto the following
  composite morphism
  \[
  T_i \xrightarrow{\mu_i^{-1}} F_\lambda \widetilde T_i
  \xrightarrow{(\lambda_g)_{\widetilde T_i}^{-1}} F_\lambda\widetilde
  T_i^g \xrightarrow{F_\lambda u} F_\lambda\widetilde T_j^h\xrightarrow{(\lambda_h)_{\widetilde T_j}} F_\lambda\widetilde
  T_j \xrightarrow{\mu_j} T_j\,,
  \]
  recall that  $\lambda_g,\lambda_h$ are
  introduced in Proposition~\ref{sec:adjunct-betw-ca-2} (part (3)).
\end{itemize}
In view of Proposition~\ref{sec:adjunct-betw-ca-2} (parts (3) and (5))
and since $\widetilde T_i^g$ has a trivial stabiliser, for all $g\in
G$ and $i\in \{1,\ldots,n\}$,  this construction does yield a
Galois covering with group $G$.
\end{proof}

When the cluster tilting objects of ${\mathcal{C}_{A*G}}$ form a mutation class,
Proposition~\ref{sec:clust-tilt-objects-6} provides the following
refinement of Corollary~\ref{sec:clust-tilt-objects-3}. Note that this
hypothesis is known to hold true when, for instance, ${\mathcal{C}_{A*G}}$ is acyclic
(\cite[Proposition 3.5]{BMRRT06}).
\begin{cor}
  \label{sec:clust-tilt-objects-9}
  Assume setting~\ref{sec:main-results-5-bis}. Assume that $G$ acts freely
  on the isomorphism classes of indecomposable projective
  $H^0(A)$-modules and that the cluster tilting objects of ${\mathcal{C}_{A*G}}$ form
  a mutation class. Then,
\begin{enumerate}
\item statements (1), (2) and (3) of
  Proposition~\ref{sec:clust-tilt-objects-6} hold true for all basic
  cluster tilting objects of ${\mathcal{C}_{A*G}}$,
\item the mapping $\mathcal T \mapsto F_\lambda\mathcal T$ induces a
  bijection from the set of cluster tilting
  subcategories of ${\mathcal{C}_A}$ stable under the action of $G$ to the set of
  cluster tilting  subcategories of
  ${\mathcal{C}_{A*G}}$.
\end{enumerate}
When cluster tilting objects are replaced by maximal rigid ones both
in the hypotheses and the conclusion, then (1) and (2) still hold true
and, in addition,
\begin{enumerate}
  \setcounter{enumi}{2}
\item given $N\in {\mathcal{C}_{A*G}}$ indecomposable, the following assertions are
  equivalent.
  \begin{enumerate}[(i)]
  \item $N$ is rigid.
  \item There exists a maximal rigid object $T\in {\mathcal{C}_A}$ such that
    $T^g\simeq T$ for all $g\in G$ and there exists an indecomposable
    direct summand $M$ of $T$ such that $N \simeq F_\lambda M$.
  \end{enumerate}
\end{enumerate}
\end{cor}
\begin{proof}
  (1) follows from
  Proposition~\ref{sec:clust-tilt-objects-6}. (2) follows from (1) and
  Corollary~\ref{sec:clust-tilt-objects-3}. (3) follows from (2).
\end{proof}

Parts of Corollary~\ref{sec:clust-tilt-objects-9} are proved in the
following more particular cases,
\begin{itemize}
\item in \cite[Section 6]{lemeur}, when replacing ${\mathcal{C}_{A*G}}$ by the
  cluster category of a piecewise hereditary algebra $R$ of quiver
  type and ${\mathcal{C}_A}$ by the cluster category of $R'$, and assuming that
  $R'$ is a Galois covering with group $G$ of $R$,
\item in \cite[Subsection 7.4]{MR3899030}, where $A$ is given by a quiver with
  potential $(Q,W)$, where the action of $G$ is given by a free action
  on $Q$ by quiver automorphisms such that $W$ is $G$-invariant up to
  cyclic permutation, and assuming that ${\mathcal{C}_{A*G}}$ is acyclic and with
  finitely many isomorphism classes of indecomposable objects.
\end{itemize}

\subsection{The acyclic case}
\label{sec:appl-finite-type}

This subsection proves Corollary~\ref{sec:main-results-7} relating ${\mathcal{C}_A}$
being acyclic to ${\mathcal{C}_{A*G}}$ being acyclic. Cluster categories of
Dynkin quivers are particular because these are exactly the
$\mathrm{Hom}$-finite, Krull-Schmidt and $2$-Calabi-Yau categories with only
finitely many isomorphism classes of indecomposable objects
(\cite[Corollary 6.6]{A07}). For these, there is a general phenomenon
as follows.
\begin{prop}
  \label{sec:appl-finite-type-1}
  Assume the setting of Proposition~\ref{sec:cy-categ-relat-1}.  The
  following assertions are equivalent.
  \begin{enumerate}[(i)]
  \item $\mathcal C$ has finitely many isomorphism classes of
    indecomposable objects.
  \item $\mathcal C'$ has finitely many isomorphism classes of
    indecomposable objects.
  \end{enumerate}
\end{prop}
\begin{proof}
  It suffices to prove ``$(i)\Longrightarrow (ii)$''. Assume (i). Let
  $M\in \mathcal C$ such that $\mathcal C = \mathrm{add}(M)$.  Then
  $\mathcal C' = \mathrm{add}(F M)$ by assumption on the unit of the adjoint
  pair $(F , F')$. Thus, (ii) holds true.
\end{proof}

Using subsections~\ref{sec:cy-categ-relat} and
\ref{sec:clust-tilt-objects} together with
Proposition~\ref{sec:appl-finite-type-1}, it is possible to prove
Corollary~\ref{sec:main-results-7}. This corollary is restated below
for convenience.
\begin{cor}
  \label{sec:acyclic-case}
    Assume setting~\ref{sec:main-results-5-bis}.
  \begin{enumerate}
  \item The following assertions are equivalent.
    \begin{enumerate}[(i)]
    \item ${\mathcal{C}_A}$ is equivalent to the cluster category of a Dynkin
      quiver.
    \item ${\mathcal{C}_{A*G}}$ is equivalent to the cluster category of a Dynkin quiver.
    \end{enumerate}
  \item If ${\mathcal{C}_A}$ is acyclic and has infinitely many isomorphism
    classes of indecomposable objects, then
    \begin{enumerate}[(a)]
    \item there exists a cluster tilting object $T\in {\mathcal{C}_A}$ such that
      $\mathrm{End}(T)$ is hereditary and $T^g = T$ for all $g\in G$ and
    \item for any such $T$, there exists a cluster tilting object
      $T'\in {\mathcal{C}_{A*G}}$ such that the algebras $\mathrm{End}(T)*G$ and $\mathrm{End}(T')$
      are Morita equivalent.
  \end{enumerate}
  Consequently, ${\mathcal{C}_{A*G}}$ is acyclic.
  \item Assume that $G$ acts freely on the set of isomorphism classes
    of indecomposable projective $H^0(A)$-modules. If ${\mathcal{C}_{A*G}}$ is
    acyclic, then so is ${\mathcal{C}_A}$.
  \end{enumerate}
\end{cor}
\begin{proof}
  (1) follows from Propositions \ref{sec:adjunct-betw-ca-2} and
  \ref{sec:appl-finite-type-1} and from \cite[Corollary 6.6]{A07}.

  (2) The following arguments are adapted from \cite{MR2467401} where
  it is proved that the skew group algebras of piecewise hereditary
  algebras are piecewise hereditary.

  Since ${\mathcal{C}_A}$ is acyclic, there exists a unique connected component
  $\Gamma$ of the Auslander-Reiten quiver of ${\mathcal{C}_A}$ having only
  finitely many $\tau$-orbits, where $\tau$ ($=\Sigma$) is the
  Auslander-Reiten translation. Note that $\Gamma$ is the repetitive
  quiver of a finite acyclic quiver because ${\mathcal{C}_A}$ has infinitely many
  isomorphism classes of indecomposable objects.

  The action of $G$ preserves the Auslander-Reiten structure of ${\mathcal{C}_A}$.
  In particular, it stabilises $\Gamma$. Following \cite[Definition
  4.1.1]{MR2467401}, let $M_0\in \Gamma$ and denote by $\mathcal X$ the
  full subquiver of $\Gamma$ consisting of the vertices $M\in \Gamma$
  such that
  \begin{itemize}
  \item there exists $g\in G$ and there exists a path from $M_0^g$ to
    $M$ in $\Gamma$ and
  \item any such path is sectional.
  \end{itemize}
  Note that $\mathcal X$ is stable under the action of $G$ up to
  isomorphisms. Moreover (\cite[Proposition 4.1.4]{MR2467401})
  $\mathcal X$ is a section of $\Gamma$, that is, $\mathcal X$ is a
  connected quiver, it is convex in $\Gamma$ and it intersects each
  $\tau$-orbit of $\Gamma$ exactly once.  Let
  $\{T_i\}_{i\in \{1,\ldots,n\}}$ be a complete family of
  representatives of $G$-orbits of indecomposable objects of
  $\mathrm{add}(\mathcal X)$ and denote by $T$ the object
  $\bigoplus\limits_{
    \begin{subarray}{c}
      1\leqslant i \leqslant n \\
      g\in G
    \end{subarray}} T_i^g$
  of ${\mathcal{C}_A}$. Then, $\mathrm{add}(\mathcal X) = \mathrm{add}(T)$ and $T^g = T$ for all
  $g\in G$. Since $\mathcal X$ is a section of $\Gamma$, the object
  $T$ is a cluster tilting object of ${\mathcal{C}_A}$ and $\mathrm{End}(T)$ is a
  hereditary algebra, isomorphic to ${\mathbbm{k}} \mathcal X$.  This proves
  (a). Note that (b) follows from (a) and from
  Proposition~\ref{sec:clust-tilt-objects-4} taking $T'=F_\lambda T$.

  Finally, it follows from (b) that $\mathrm{End}(T)$ and $\mathrm{End}(F_\lambda T)$
  have the same global dimension (\cite[Theorem 1.3]{RR85}). Hence,
  $\mathrm{End}(F_\lambda T)$ is hereditary. Thus, ${\mathcal{C}_{A*G}}$ is acyclic (see
  \cite{KR08}).

  (3) There exists a cluster tilting object $T'\in {\mathcal{C}_{A*G}}$ such that
  $\mathrm{End}(T')$ is a hereditary algebra. In particular, the cluster
  tilting objects of ${\mathcal{C}_{A*G}}$ form a mutation class
  (see \cite[Proposition 3.5]{BMRRT06}). By
  Corollary~\ref{sec:clust-tilt-objects-9}, there exists a cluster
  tilting object $T\in {\mathcal{C}_A}$ such that $T^g\simeq T$ for all $g\in G$
  and such that $T' \simeq F_\lambda T$. By
  Proposition~\ref{sec:clust-tilt-objects-4}, the algebras $\mathrm{End}(T)*G$
  and $\mathrm{End}(T')$ are Morita equivalent. Since $\mathrm{End}(T')$ is
  hereditary, so is $\mathrm{End}(T)$. Thus, ${\mathcal{C}_A}$ is acyclic (see
  \cite{KR08}).
\end{proof}

\section{Application to higher Auslander-Reiten theory}
\label{sec:appl-high-ausl}

This section applies the study of
subsection~\ref{sec:smash-prod-deform} on Calabi-Yau completions to the
interplay between skew group algebras and $d$-representation infinite
algebras in the sense of \cite{MR3144232}.  For the sake of
comprehensiveness, the corresponding interplay for $d$-representation
finite algebras is presented in
subsection~\ref{sec:n-repr-finite}. The case of $d$-representation
infinite algebras is presented in
subsection~\ref{sec:texorpdfstr-repr-inf}. The setting is the
following.
\begin{setting}
  \label{sec:appl-high-ausl-1}
  Let $A$ be a finite dimensional algebra. Let $d$ be a non negative
  integer. Let $G$ be a finite group acting on $A$ by algebra
  automorphisms. Assume that $\mathrm{char}({\mathbbm{k}})\nmid\mathrm{Card}(G)$.
\end{setting}
Until the end of the text, the category of finitely generated modules
over a finite dimensional algebra $R$ is denoted by $\mathrm{mod}(R)$.

\subsection{\texorpdfstring{$d$}{d}-representation finite algebras}
\label{sec:n-repr-finite}

Given a finite dimensional algebra $R$, recall (\cite{MR2735750}) that
a \emph{$d$-cluster tilting subcategory} of $\mathrm{mod}(R)$ is a
functorially finite subcategory $\mathcal C$ such that
\[
\left\{
  \begin{array}{rcl}
    \mathcal C& = & \{M\in \mathrm{mod}(R)\ |\ (\forall X\in \mathcal C)\
                    (\forall i=1,\ldots,d-1)\ \mathrm{Ext}_R^i(X,M)=0\} \\
    & = & \{M\in \mathrm{mod}(R)\ |\ (\forall X\in \mathcal C)\
                    (\forall i=1,\ldots,d-1)\ \mathrm{Ext}_R^i(M,X)=0\}\,,
  \end{array}\right.
  \]
  and a \emph{$d$-cluster tilting module} is a module $M\in \mathrm{mod}(R)$
  such that $\mathrm{add}(M)$ is a $d$-cluster tilting subcategory. The algebra
  $R$ is called \emph{$d$-representation finite} when such an $M$
  exists and when $R$ has global dimension at most $d$. When such an
  $M$ exists and when $d$ is equal to the global dimension of $R$,
  then $\mathrm{add}(M)$ is necessarily equal to the following subcategory
  denoted by $\mathcal I_R$ (see \cite[Proposition 0.2]{MR2820136}),
  \begin{equation}
    \label{eq:73}
    \mathcal I_R = \mathrm{add}(\{\tau_d^kI\ |\ k\in \mathbb N,\ I\in \mathrm{mod}(A)\
    \text{injective}\})\,,
\end{equation}
where $\tau_d := D\mathrm{Ext}_R^d(-,R) \colon \mathrm{mod}(R)\to \mathrm{mod}(R)$ is the
higher Auslander-Reiten translation ($D:=\mathrm{Hom}_{\mathbbm{k}}(-,{\mathbbm{k}})$).

The aim of this subsection is to prove that, when $A$ has global
dimension $d$, then $A$ is $d$-representation finite if and
only if so is $A*G$. This is based on the following lemma.
\begin{lem}
  \label{sec:n-repr-finite-1}
  Assume setting~\ref{sec:appl-high-ausl-1}. Denote by $F_\lambda$ and
  by $F_\cdot$ the extension-of-scalars functor
  $\mathrm{mod}(A) \to \mathrm{mod}(A*G)$ and the restriction-of-scalars functor
  $\mathrm{mod}(A*G)\to \mathrm{mod}(A)$, respectively.
  \begin{enumerate}
  \item $\tau_dF_\lambda M \simeq F_\lambda \tau_d M$ for all
    finite dimensional $A$-modules $M$.
  \item $\tau_d F_\cdot N \simeq F_\cdot \tau_dN$ for all finite
    dimensional $A$-modules $N$.
  \end{enumerate}
\end{lem}
\begin{proof}
  (1) Let $M\in \mathrm{mod}(A)$. By adjunction, there is an isomorphism
  $D(\mathrm{Hom}_{A*G}(F_\lambda M, A*G))\simeq D(\mathrm{Hom}_A(M,A*G))$ in
  $\mathrm{mod}(A*G)$. The following composite mapping is bijective,
  \begin{equation}
    \label{eq:69}
    D(\mathrm{Hom}_A(M,A*G)) \to \bigoplus_{g\in G} D(\mathrm{Hom}_A(M,A^g)) \to
    \bigoplus_{g\in G} D(\mathrm{Hom}_A(M,A))\,,
  \end{equation}
  where the leftmost arrow is induced by the mapping
  $\oplus_{g\in G} A^g \to A*G$ given by
  $(a_g)_{g\in G} \mapsto \sum_{g\in G}a_g*g$ and the rightmost arrow
  is induced by the mappings $A \to A^g$ given by $a \mapsto \,^ga$,
  for all $g\in G$.  Let $\varphi \in D(\mathrm{Hom}_A(M,A*G))$. If its image
  under \eqref{eq:69} is denoted by $(\lambda_g)_{g\in G}$, then
  \[
  (\forall g\in G)\ (\forall f\in \mathrm{Hom}_A(M,A))\ \ \lambda_g(f) =
  \varphi\left(
    \begin{array}{ccc}
      M & \to & A*G \\
      m & \mapsto & \,^gf(m)*g
    \end{array}
  \right)\,.
  \]
  Given $a\in A$, if $(\mu_g)_{g\in G}$ denotes the image of
  $\varphi\cdot (a*e)$ under \eqref{eq:69}, then
  $\mu_g(f) = \lambda_g(\,^{g-1}a\cdot f)$ for all $g\in G$ and
  $f\in \mathrm{Hom}_A(M,A)$; in other words,
  $\mu_g = \lambda_g\cdot (\,^{g^{-1}}a)$ for all $g\in G$. Besides,
  given $h\in H$, if $(\nu_g)_{g\in G}$ denotes the image of
  $\varphi \cdot (1*h)$ under \eqref{eq:69}, then
  $\nu_g(f) = \lambda_{hg}(f)$ for all $f\in \mathrm{Hom}_A(M,A)$; in other
  words, $\nu_g = \lambda_{hg}$ for all $g\in G$.
  
  It follows from these considerations that the composition of
  the following mapping with \eqref{eq:69} is an isomorphism in $\mathrm{mod}(A*G)$
  from $D(\mathrm{Hom}_A(M,A*G))$ to $D(\mathrm{Hom}_A(M,A))\otimes_AA*G$.
  \begin{equation}
    \label{eq:52}
    \begin{array}{ccc}
      \bigoplus_{g\in G} D(\mathrm{Hom}_A(M,A)) & \to & D(\mathrm{Hom}_A(M,A)) \otimes_A A*G
      \\
      (\lambda_g)_{g\in G} & \mapsto & \sum_{g\in G}\lambda_g \otimes
                                       (1*g^{-1})\,.
    \end{array}
  \end{equation}
   Therefore,
  there is a functorial isomorphism in $\mathrm{mod}(A*G)$,
  \begin{equation}
    \label{eq:70}
    D(\mathrm{Hom}_{A*G}(F_\lambda M, A*G))\simeq F_\lambda D(\mathrm{Hom}_A(M,A))\,.
  \end{equation}
  Now, the isomorphism $\tau_d F_\lambda M \simeq F_\lambda \tau_d M$
  follows from \eqref{eq:70} after deriving because $F_\lambda$ is
  exact and preserves projective modules.

  (2) Let $N\in \mathrm{mod}(A*G)$. By adjunction,
  there is an isomorphism in $\mathrm{mod}(A)$,
  \begin{equation}
    \label{eq:71}
    D(\mathrm{Hom}_A(F_\cdot N,A)) \simeq D(\mathrm{Hom}_{A*G}(N, A*G))\,.
  \end{equation}
  And the isomorphism $\tau_d F_\cdot N \simeq F_\cdot \tau_d N$
  follows from \eqref{eq:71} after deriving also because $F_\cdot$ is
  exact and preserves projective modules.
\end{proof}

Here is the first main result of this section.
\begin{prop}
  \label{sec:texorpdfstr-repr-fin}
  Assume setting~\ref{sec:appl-high-ausl-1}. The following assertions
  are equivalent.
  \begin{enumerate}[(i)]
  \item $A$ is $d$-representation finite.
  \item $A*G$ is $d$-representation finite.
  \end{enumerate}
\end{prop}
\begin{proof}
  Denote by $F_\lambda$ and by $F_\cdot$ the extension-of-scalars
  functor $\mathrm{mod}(A) \to \mathrm{mod}(A*G)$ and the restriction-of-scalars
  functor $\mathrm{mod}(A*G)\to \mathrm{mod}(A)$, respectively.  Note that $A*G$ has
  global dimension $d$.  Since $F_\lambda D(A) \simeq D(A*G)$ as
  $A*G$-modules and $F_\cdot D(A*G) \simeq D(A)^{\mathrm{Card}(G)}$ as
  $A$-modules, it follows from Lemma~\ref{sec:n-repr-finite-1} that
   \begin{equation}
     \label{eq:72}
     \left\{
     \begin{array}{rcl}
       \mathrm{add}(F_\lambda \mathcal I_A) & = & \mathcal I_{A*G} \\
       \mathrm{add}(F_\cdot \mathcal I_{A*G}) & = & \mathcal I_A\,.
     \end{array}\right.
 \end{equation}

 Assume (i) and let $M\in \mathrm{mod}(A)$ be a $d$-cluster tilting module
 such that $\mathrm{add}(M) = \mathcal I_A$. By construction, $\mathcal I_A$ is
 stable under the action of $G$. Hence, so is $\mathrm{add}(M)$. Therefore,
 $\mathrm{add}(M)$ is stable under $\mathrm{add}(F_\cdot F_\lambda)$ (see
 Lemma~\ref{sec:adjunctions-1}, part (4)). Hence, the proof of
 Proposition~\ref{sec:cy-categ-relat-1} can be adapted here to prove
 that $\mathrm{add}(F_\lambda M)$ is a $d$-cluster tilting subcategory of
 $\mathrm{mod}(A*G)$ provided that Lemma~\ref{sec:adjunctions-1} is used
 instead of Proposition~\ref{sec:adjunct-betw-ca-2}. Therefore,
 $(i)\Longrightarrow (ii)$.

 Conversely, assume (ii) and let $N\in \mathrm{mod}(A*G)$ be a $d$-cluster
 tilting module such that $\mathrm{add}(N) = \mathcal I_{A*G}$. Since
 $\mathcal I_{A*G}$ is stable under $\mathrm{add}(F_\lambda F_\cdot)$
 (\eqref{eq:72}), the proof of Proposition~\ref{sec:cy-categ-relat-1}
 can also be adapted here to prove that $\mathrm{add}(F_\cdot N)$ is a
 $d$-cluster tilting subcategory of $\mathrm{mod}(A)$. Therefore
 $(ii)\Longrightarrow (i)$.
\end{proof}

\subsection{\texorpdfstring{$d$}{d}-representation infinite algebras}
\label{sec:texorpdfstr-repr-inf}

Let $R$ be a finite dimensional ${\mathbbm{k}}$-algebra. Following
\cite{MR3144232}, it is called \emph{$d$-representation infinite} when
it has global dimension at most $d$ and $\mathbb S_d^{-i}R$ is an
$R$-module for all $i\in \mathbb N$, where $\mathbb S_d$ denotes
$(-\otimes^{\mathrm L}_RDR)\circ \Sigma^{-d}$. This entails that $R$ has global
dimension equal to $d$. When $R$ has global dimension at most $d$, it
is $d$-representation infinite if and only if the cohomology of
$\mathbf{\Pi}_{d+1}(R)$ is concentrated in degree $0$. Recall that the
Calabi-Yau completion $\mathbf{\Pi}_{d+1}(R)$ is also called the \emph{derived
  $(d+1)$-preprojective} algebra and its $0$-th cohomology is
isomorphic to the \emph{$(d+1)$-preprojective algebra} $\Pi_{d+1}(R)$
of $R$ which is defined by $\Pi_{d+1}(R)= T_R(\mathrm{Ext}_R^d(DR,R))$
(\cite{MR3077865}).

\begin{prop}
  \label{sec:texorpdfstr-repr-inf-1}
  Assume setting~\ref{sec:appl-high-ausl-1}.
  \begin{enumerate}
  \item $G$ acts  on the derived $d+1$-preprojective algebra
    $\mathbf{\Pi}_{d+1}(A)$ by dg automorphisms and $\mathbf{\Pi}_{d+1}(A)*G$ is
    quasi-isomorphic to the derived $d+1$-preprojective algebra
    $\mathbf{\Pi}_{d+1}(A*G)$.
  \item $G$ acts on the $d+1$-preprojective algebra $\Pi_{d+1}(A)$ by
    algebra automorphisms and the skew group algebra $\Pi_{d+1}(A)*G$
    is isomorphic to the $d+1$-preprojective algebra $\Pi_{d+1}(A*G)$.
  \item The following assertions are equivalent.
    \begin{enumerate}[(i)]
    \item $A$ is $d$-representation infinite.
    \item $A*G$ is $d$-representation infinite.
    \end{enumerate}
  \end{enumerate}
\end{prop}
\begin{proof}
  (1) is a rephrasing of Theorem~\ref{sec:main-results-2}. (2) and (3)
  follow from (1).
\end{proof}

Note that subsection~\ref{sec:smash-prod-deform} applies here. Since
$A$ has global dimension to $d$, there exists a bounded complex $P$ of
finitely generated projective $A$-bimodules endowed with a compatible
action of $G$ and with $G$-equivariant differential, and there exists
a $G$-equivariant quasi-isomorphism $P\to A$ in $\mathrm{Mod}(A^e)$. Then,
there is a compatible action of $G$ on $\mathrm{Hom}_{A^e}(P,A^e)$
such that $(\,^gf)(p) = \,^g(f(\,^{g^{-1}}p))$ for all $p\in P$,
$g\in G$ and $f\in \mathrm{Hom}_{A^e}(P,A^e)$. Hence, there is an
isomorphism
$\mathrm{Hom}_{A^e}(P,A^e) \simeq \mathrm{RHom}_{A^e}(A,A^e)$ in
$\mathcal{D}(\Delta)$ and, since $P$ is bounded,
$\mathrm{Hom}_{A^e}(P,A^e)$ is cofibrant in $\mathrm{Mod}(A^e)$. Therefore,
one may assume that
\begin{itemize}
\item $\mathbf{\Pi}_{d+1}(A) = T_A(\Sigma^d\mathrm{Hom}_{A^e}(P,A^e))$ as dg algebras
  respectively endowed with an action of $G$ by dg automorphisms, and 
  $\mathbf{\Pi}_{d+1}(A*G)\underset{\mathrm{qis}}{\simeq}
  T_A(\Sigma^d\mathrm{Hom}_{A^e}(P,A^e))*G$,
  and, applying $H^0(-)$, that
\item $\Pi_{d+1}(A) \simeq T_A(H^d\mathrm{Hom}_{A^e}(P,A^e))$ as
  algebras respectively endowed
  with an action of $G$ by algebra automorphisms, and $\Pi_{d+1}(A*G)
  \simeq  T_A(H^d\mathrm{Hom}_{A^e}(P,A^e))*G$.
\end{itemize}

\section{Acknowledgements}

I thank Yann Palu for helpful discussions on cluster tilting theory.
\bibliographystyle{plain}
\bibliography{biblio}

\end{document}